\documentclass{tac}

\usepackage[all, knot, arc]{xy}
\usepackage{amsfonts}
\usepackage{graphicx}
\usepackage{mathtools}
\usepackage{bbm}

\author{Dimitri Chikhladze}

\thanks{The research was supported by the Funda\c{c}\~ao para a Ci\^encia e a Tecnologia (Portugal) Postdoctoral Fellowship SFRH/BPD/79360/2011, by the Shota Rustaveli National Science Foundation (Georgia) grants PG/45/5-113/12 and DI/18/5-133/13, and by the Centre for Mathematics of the University of Coimbra (funded by the European Regional Development Fund program COMPETE and by the FCT project PEst-C/MAT/UI0324/2013).}

\address{Centro de Matem\'atica da Universidade de Coimbra \\
Apartado 3008, 3001-501 Coimbra, Portugal}

\title{Lax formal theory of monads, monoidal approach to bicategorical structures and generalized operads}

\copyrightyear{2015}



\eaddress{d.chikhladze@gmail.com}

\def\beq#1{\begin{equation}\label{#1}}
\def\eeq{\end{equation}}

\def\hiso{\rotatebox{0}{\scalebox{1.2}[1.2]{$\cong$}}}
\def\viso{\rotatebox{270}{\scalebox{1.2}[1.2]{$\cong$}}}
\def\hequ{\rotatebox{0}{\scalebox{1.2}[1.2]{$=$}}}

\def\l{+}
\def\r{-}
\def\o{\circ}
\def\s{\ast}
\def\d{\mathrm{\top}}
\def\u{\mathrm{\bot}}
\def\op{\mathrm{^{op}}}
\def\cat{\mathcal{C}\mathrm{at}}
\def\span{\mathcal{S}\mathrm{pan}}
\def\twocat{2\text{-}\mathcal{C}\mathrm{at}}
\def\set{\mathcal{S}\mathrm{et}}
\def\mat{\mathcal{M}\mathrm{at}}
\def\bicat{\mathcal{B}icat}
\def\twoprof{2\text{-}\mathcal{P}\mathrm{rof}}
\def\mod{\mathcal{M}\mathrm{od}}
\def\prof{\mathcal{P}\mathrm{rof}}
\def\talg{T\text{-}\mathcal{A}\mathrm{lg}}
\def\mon{\mathcal{M}\mathrm{on}}
\def\tmon{T\text{-}\mathcal{M}\mathrm{on}}
\def\B{\mathcal{B}}
\def\C{\mathcal{C}}
\def\T{\mathcal{T}}
\def\L{\mathcal{L}}
\def\M{\mathcal{M}}
\def\K{\mathcal{K}}
\def\MMatV{\mathbb{M}\mathrm{at}(V)}
\def\SSpan{\mathbb{S}\mathrm{pan}}
\def\MMod{\mathbb{M}\mathrm{od}}
\def\PProf{\mathbb{P}\mathrm{rof}}
\def\AA{\mathbb{A}}
\def\BB{\mathbb{B}}
\def\FF{\mathbb{F}}
\def\GG{\mathbb{G}}
\def\TT{\mathbb{T}}
\def\II{\mathbb{I}}
\def\LL{\mathbb{L}}
\def\NN{\mathbb{N}}
\def\MM{\mathbb{M}}
\def\SS{\mathbb{S}}
\def\RR{\mathbb{R}}
\def\one{\mathbbm{1}}
\def\tt{\mathbbm{t}}
\def\ee{\mathbbm{e}}
\def\mm{\mathbbm{m}}
\def\dd{\mathfrak{d}}
\def\fM{\mathfrak{M}}
\def\comp{\mathfrak{Comp}}
\def\id{\mathfrak{Un}}
\def\kl{\mathfrak{Kl}}
\def\Eq{\mathfrak{Eq}}
\def\sEq{^\ast\Eq}
\def\mll{\mathfrak{L}^{\bot\circ\circ}}
\def\mcoll{\mathfrak{L}^{\bot\bullet\circ}}
\def\mcolcol{\mathfrak{L}^{\bot\bullet\bullet}}
\def\mlcol{\mathfrak{L}^{\bot\circ\bullet}}
\def\moll{\mathfrak{L}^{\top\circ\circ}}
\def\mocolcol{\mathfrak{L}^{\top\bullet\bullet}}
\def\mocoll{\mathfrak{L}^{\top\bullet\circ}}
\def\molcol{\mathfrak{L}^{\top\circ\bullet}}

\makeatletter
\def\slashedarrowfill@#1#2#3#4#5{%
  $\m@th\thickmuskip0mu\medmuskip\thickmuskip\thinmuskip\thickmuskip
   \relax#5#1\mkern-7mu%
   \cleaders\hbox{$#5\mkern-2mu#2\mkern-2mu$}\hfill
   \mathclap{#3}\mathclap{#2}%
   \cleaders\hbox{$#5\mkern-2mu#2\mkern-2mu$}\hfill
   \mkern-7mu#4$%
}
\def\rightslashedarrowfill@{%
  \slashedarrowfill@\relbar\relbar\mapstochar\rightarrow}
\newcommand\xslashedrightarrow[2][]{%
  \ext@arrow 0055{\rightslashedarrowfill@}{#1}{#2}}
\makeatother
\def\vectarr{\xslashedrightarrow{}}

\begin{document}

\maketitle

\begin{abstract}
Generalized operads, also called generalized multicategories and $T$-monoids, are defined as monads within a Kleisli bicategory. With or without emphasizing their monoidal nature, generalized operads have been considered by numerous authors in different contexts, with examples including symmetric multicategories, topological spaces, globular operads and Lawvere theories. In this paper we study functoriality of the Kleisli construction, and correspondingly that of generalized operads. Motivated by this problem we develop a lax version of the formal theory of monads, and study its connection to bicategorical structures.
\end{abstract}

\section{Introduction}\label{introduction}
As the title suggests, this paper revolves around three themes. The first of them is developing a new general framework for the theory of generalized multicategories. The second is generalizing the formal theory of monads internal to a bicategory to its lax version internal to a tricategory. The third is an idea of a monoidal approach to bicategorical structures, which also serves as a bridge between the two other themes. 

The concept of a generalized multicategory, also called a generalized operad and a $T$-monoid, involves few steps of abstraction. The basic notion of a multicategory \cite{Her00} is a generalization of a category, in which the domain of a morphism, instead of being a single object, is allowed to be a finite list of objects. A one-object multicategory is a non-symmetric operad of \cite{May72}. At the next step, one observes that the domain of a morphism of a multicategory is an element of the free monoid on the set of its objects, and replaces the free-monoid construction by an arbitrary monad. Furthermore, from the context internal to the category of sets one switches to a more general ambient, so as to allow structures such as enriched multicategories. Numerous works following this paradigm include \cite{Bur71, Kell92, BD98, ClTh03, DS03, Lei04, Che04, Her04, Web05, Sea05, FGHW08, Gar08}, with examples as diverse as symmetric multicategories, topological spaces, metric spaces, globular operads and Lawvere theories. A unifying approach was developed in \cite{CrSh10}, where a comprehensive account of the subject can be also found. We develop a new framework for generalized multicategories which contains abstractly all other contexts. The framework is essentially at the same level of generality as that of \cite{CrSh10}. The difference from the latter is that we took a more structural algebraic approach, which we deemed more appropriate for our purposes. 

More concretely, a generalized multicategory, or a $T$-monoid, is defined internal to a bicategory-like structure $\AA$, and with respect to a monad-like structure on it $\TT$. First, by the Kleisli construction from the data $(\AA, \TT)$ one produces another bicategory-like structure $\kl(\AA, \TT)$, and then defines a $T$-monoid to be a monoid, or a monad, within the the latter. To develop a precise theory, first one needs to formalize the data  $(\AA, \TT)$. We formalize the bicategory-like structure $\AA$ under the name of an equipment, and formalize the data $(\AA, \TT)$ under the name of a $T$-equipment. Further we study functoriality of the Kleisli construction. There are several interesting notions of a morphism between $T$-equipment (in which both $\AA$ and $\TT$ may vary), corresponding to different 2-categories of $T$-equipments. One of them serves as a domain for the Kleisli construction, which becomes a 2-functor $\kl$. We introduce $\ast$-equipments and $T$-$\ast$-equipments, which make equipments closer to the proarrow equipments of \cite{Woo82}, and thus reflect more structure usually present in the examples. It turns out that the Kleisli construction on $T$-$\ast$-equipments has another 2-functorial extension $^{\ast}\kl$. The functoriality of the Kleisli construction can be used to compare to each other the categories of $T$-monoids within different $T$-equipments. Furthermore, it can be used as a technical tool for various constructions on $T$-equipments and $T$-monoids within them. As an application of this technique, we construct a free $T$-algebra functor and the underlying $T$-monoid functor, which are analogues of the free monoidal category functor and the underlying multicategory functor going between the categories of monoidal categories and multicategories. We then generalize the results of \cite{ChClHo14}.   

The second theme of the paper is the lax formal theory of monads within a tricategory. This is a generalization of the formal theory of monads within a bicategory, originally developed in \cite{St72}, through a ``lax categorification'' at the second dimension. By the latter we mean switching to the context internal to a tricategory, and modifying the given theory by replacing the equations between 2-cells by non-invertible 3-cells. We consider lax monads within a tricategory, defined as lax monoids of \cite{DS03} in endohom monoidal bicategories. Furthermore, we consider lax variants of the categories of monads of \cite{St72}, study lax distributive laws, and introduce a construction of composition of a pair of lax monads related by a lax distributive law $\comp$.

The idea of the third theme is a two dimensional analogue of the simple fact that a category is a monad in the bicategory of spans. More specifically, as the composition structure of a category can be encoded by the multiplication structure of a monad in the bicategory of spans, so a horizontal composition of a bicategory-like structure can be encoded by a multiplication structure of a monad-like structure in the tricategory of pseudoprofunctors, which are higher dimensional analogues of spans. The first step here is to define the tricategory of pseudoprofunctors $\mod$. This has categories as its objects and pseudoprofunctors, or modules, as its morphisms. Furthermore, in order to be able to consider functors between bicategory-like structures, one needs to work with an embedding
\[\xymatrix{\cat^{\op} \ar[r] & \mod.}\]

\noindent The tricategory $\mod$ is perhaps well known. We however give an independent outline of its definition. In fact, we define a tricategory whose objects are bicategories, and whose morphisms are biprofunctors $\twoprof$. This itself is done through monad theoretic approach, by observing bicategories to be pseudomonads in a certain tricategory, and biprofunctors to be bipseudomodules of pseudomonads. 

Our equipments are observed to be lax monads in $\mod$. This provides a bridge between the lax formal theory of monads and the theory of generalized multicategories. The constructions of the latter are then established to be expressible by the constructions of the former. In particular, it is shown that the Kleisli construction $\kl$ is an instance of the distributive composition $\comp$. Note also, that $T$-monoids are defined as monads within equipments which themselves are monads. This is a kind of microcosm principle. 

The structure of the paper is the following. In Section~\ref{monads} we review the formal theory of monads. In sections~\ref{equipments}--\ref{starequipments} we study equipments and the theory of generalized multicategories. In Section~\ref{psmonpsmod} we consider pseudomonads and pseudomodules, which we use in Section~\ref{biprofmod} to construct the tricategories $\twoprof$ and $\mod$. In Section~\ref{laxmonads} we develop the lax formal theory of monads, revisiting equipments and the theory of generalized multicategories as an example. In Section~\ref{further} we briefly gives some further perspectives on the subject.

\section{Monads in a bicategory}\label{monads}
\smallskip In this section we recollect the formal theory of monads within a bicategory. Most of the material is essentially from \cite{St72}. We however introduce it under different notation and terminology.

\smallskip A monad $\TT = (X, T) = (X, T, m, e)$ in a bicategory $\B$ consists of an object $X$ of $\B$ and a monoid $(T, m : T^2 \rightarrow T, e : 1_{X} \rightarrow T)$ in the endohom monoidal category $\B(X, X)$. A monad \textbf{upmap} $\FF = (F, u) : (X, T) \rightarrow (Y, S)$ (called a monad map in \cite{St72}) consists of a morphism $F : X \rightarrow Y$ and a 2-cell
\[
\xymatrix{
X \ar[rr]^T_{\quad}="1"  \ar[d]_{F} && X \ar[d]^{F}\\
Y \ar[rr]_S^{\quad}="2"&& Y 
\ar@{=>} "2"+(0,3.5); "2"+(0,8.5) _u
}
\]

\noindent satisfying two axioms, expressing compatibility with the monad multiplication and unit. A monad upmap transformation $\FF \rightarrow \GG$ is a 2-cell $t : F \Rightarrow G$ satisfying one axiom. Monads, monad upmaps and monad upmap transformations form a bicategory which we denote by $\fM^\bot(\B)$. Another bicategory whose objects are monads in the bicategory $\B$ is defined by the formula $\fM^\top(\B) = \fM(\B^\op)^\op$. A monad \textbf{downmap} (called a monad opmap in \cite{St72}) is a morphism of $\fM^\top(\B)$. More explicitly, a monad downmap $\FF = (F, d) : (X, T) \rightarrow (Y, S)$ consists of a morphism $F : X \rightarrow Y$ and a 2-cell 
\[
\xymatrix{
X \ar[rr]^T_{\quad}="1"  \ar[d]_{F}&& X \ar[d]^{F}\\
Y \ar[rr]_S^{\quad}="2"&& Y
\ar@{=>} "1"+(0,-3.5); "1"+(0,-8.5) ^d
}
\]

\noindent satisfying two axioms. A monad downmap transformation is a morphism of $\fM^\top(\B)$. 

For any object $X$ of the bicategory $\B$, there is a trivial monad $\id(X) = (X, 1_X)$. We have two functors $\id_\B^\bot: \B \rightarrow \fM^\bot(\B)$ and $\id_\B^\top : \B \rightarrow \fM^\top(\B)$ extending $\id$ to morphisms and 2-cells of $\B$ in the obvious way. For $\id(X)$ we will shortly write $\one_X$.

A (right) module of a monad $\TT = (X, T)$,  is an object $Z$ together with a monad downmap $\one_Z \rightarrow \TT$. Thus essentially, a module is a 2-cell
\[
\xymatrix{
&Z  \ar[rd]^F="2" \ar[ld]_F="1"  &\\
X \ar[rr]_T& \ar@{=>}"1"+(6,-4);"1"+(11,-4) ^d& X
}
\]

\noindent satisfying two axioms. A module of $\TT$ is an object of the comma category $\id_\B^\top\downarrow\TT$. When it exists, the terminal object in this category is called the EM object of the monad $\TT$, and is denoted by $X^T$. For any module $\one_Z \rightarrow \TT$ the unique map $Z \rightarrow X^T$ in $\id_\B^\top\downarrow\TT$ is called the comparison morphism. 

Take $\B$ to be the 2-category of categories $\cat$. A monad $\TT = (X, T)$ in $\cat$ is a usual monad $T$ on a category $X$. The EM object  $X^T$ is then the category of $T$-algebras. An object of $X^T$, i.e. an algebra  $(x, h : Tx \rightarrow x)$, can be itself identified with a module of $\TT$ given by a 2-cell
\[
\xymatrix{
&I  \ar[rd]^x="2" \ar[ld]_x="1"  &\\
X \ar[rr]_T& \ar@{=>}"1"+(6,-4);"1"+(11,-4) ^h& X.
}
\]
 
\noindent where $I$ denotes the terminal category, and $x : I \rightarrow X$ denotes the functor which chooses the object $x$. 

Suppose now that $\TT = (X, T)$ is a monad in a 2-category $\B$. Fix an object $Z$. Then the functor $\B(Z, T)$ becomes a monad on the category $\B(Z, X)$. The category of EM algebras for this monad is the same as the category $\fM_\B^\top(\one_Z, \TT)$ of modules of $\TT$ with the fixed underlying object $Z$. 

Further recall from \cite{St72}, that a \textbf{distributive pair of monads} $(\SS, \TT, c)$ consists of monads $\TT = (T, X)$ and $\SS = (S, X)$ on the same object $X$, and a 2-cell
\[
\xymatrix{
X \ar[rr]^{S}_{\quad}="1"  \ar[d]_{T} && X \ar[d]^{T}\\
X \ar[rr]_{S}^{\quad}="2"&& X 
\ar@{=>} "2"+(0,3.5); "2"+(0,8.5) ^{c}
}
\]

\noindent such that equivalently:

\begin{itemize}
\item $(T, c)$ is  a monad upmap $\SS \rightarrow \SS$, and the monad multiplication and unit of $\SS$ are monad upmap transformations.
\item $(S, c)$ is a monad downmap $\TT \rightarrow \TT$, and the monad multiplication and unit of $\TT$ are monad downmap transformations.
\end{itemize}

\noindent It follows that a distributive pair of monads determines a monad in $\fM^\bot(\B)$, and also a monad in $\fM^\top(\B)$. There are four bicategories whose objects are distributive pairs of monads:
\[\fM^\bot\fM^\bot(\B),\qquad \fM^\bot\fM^\top(\B),\qquad\fM^\top\fM^\bot(\B),\qquad\fM^\top\fM^\top(\B).\]

\noindent Let us identify the morphisms and 2-cells of each of them: A morphism $(\SS', \TT') \rightarrow (\SS, \TT)$ of $\fM^\bot\fM^\bot(\B)$ consists of a pair of monad upmaps $\FF = (F, u) : \TT' \rightarrow \TT$ and $\GG = (G, u') : \SS' \rightarrow \SS$, such that $F = G$ and the equation
\[
\xymatrix@R=0.8em{
&& X \ar[drr]^{S'}_\;="2"&&\\
X  \ar[rru]^{T'}_{\;}="1" &&&& X  \\
&& X  \ar[uu];[]^{F} \ar[drr]^{S}="4" ="6"&& \\
X  \ar[urr]^{T}="3" \ar[uu];[]_{F}  \ar[drr]_{S}&&&&X \ar[uu];[]^{F}  \\
&&X  \ar[urr]_{T}  ="5" &\ar@{=>} "5"+(-13,3.5); "5"+(-13,8.5) ^c & \ar@{=>}"3"+(0,4.5);"3"+(0,9.5) ^u \ar@{=>}"4"+(0,4.5);"4"+(0,9.5) _{u'}
}
\qquad
\xymatrix@C=1em{\ar@{}[ddd]|-{\hequ} \\ \\ \\  \\}
\qquad
\xymatrix@R=0.8em{
&& X  \ar[drr]^{S'} ="5"&&\\
X  \ar[rru]^{T'}  \ar[drr]^{S'}_\; ="1"&&&& X  \\
&& X \ar[urr]^{T'}_\;="2" && \\
X  \ar[uu];[]_{F} \ar[drr]_{S}^\;="3"&&&&X \ar[uu];[]^{F} \\
&&X \ar[urr]_{T}^{\;}="4"   \ar[uu];[]^{F}&& \ar@{=>}"3"+(0,4.5);"3"+(0,9.5) ^{u'} \ar@{=>}"4"+(0,4.5);"4"+(0,9.5) _u \ar@{=>} "2"+(-12,3); "2"+(-12,8) ^{c'}
}
\]

\noindent holds. A 2-cell $(\GG', \FF') \rightarrow (\GG, \FF)$ in $\fM^\bot \fM^\bot(\B)$ is a 2-cell $F' \rightarrow F$ which becomes an upmap transformation both, between $\FF'$ and $\FF$, and between $\GG'$ and $\GG$. A morphism $(\SS', \TT') \rightarrow (\SS, \TT)$ of $\fM^\top\fM^\bot(\B)$ consists of an upmap  $\FF = (F, u) : \TT' \rightarrow \TT$ and a downmap $\GG = (G, d) : \SS' \rightarrow \SS$ such that $F = G$, and the equation
\[
\xymatrix@R=0.8em{
&& X \ar[drr]^{T'}_\;="2"&&\\
X \ar[rru];[]_{S'}_{\;}="1" &&&& X  \\
&& X  \ar[uu];[]^{F} \ar[drr]_{T}="4" && \\
X  \ar[urr];[]^{S}="3" \ar[uu];[]_{F}  \ar[drr]_{T}&&&&X \ar[uu];[]^{F}  \\
&&X \ar[urr];[]^{S} && \ar@{=>} "1"+(0,-6.5); "1"+(0,-11.5) ^{d} \ar@{=>}"4"+(0,7.5);"4"+(0,12.5) ^{u} \ar@{=>} "1"+(9.5,-22); "1"+(14.5,-22) ^c
}
\qquad
\xymatrix@C=1em{\ar@{}[ddd]|-{\hequ} \\ \\ \\  \\}
\qquad
\xymatrix@R=0.8em{
&& X  \ar[drr]^{T'}&&\\
X  \ar[rru];[]_{S'}  \ar[drr]^{T'}_\;="1"&&&& X  \\
&& X \ar[urr];[]_{S'}_\;="2"  && \\
X  \ar[uu];[]_{F} \ar[drr]_{T}^\;="3"&&&&X \ar[uu];[]^{F} \\
&&X \ar[urr];[]^{S}_{\;}="4"   \ar[uu];[]^{F}&& \ar@{=>}"3"+(0,4.5);"3"+(0,9.5) ^{u} \ar@{=>}"2"+(0,-6.5); "2"+(0,-11.5) ^{d} \ar@{=>} "1"+(9.5,6); "1"+(14.5,6) ^{c'}
}
\]

\noindent holds. A 2-cell in $\fM^\top \fM^\bot(\B)$ $(\GG', \FF') \rightarrow (\GG, \FF)$ is a 2-cell $F' \rightarrow F$ which becomes both, an upmap transformation between $\FF'$ and $\FF$, and a downmap transformation between $\GG'$ and $\GG$. The bicategory $\fM^\bot\fM^\top(\B)$ can be easily shown to be isomorphic to $\fM^\top\fM^\bot(\B)$. Finally, a morphism $(\SS', \TT') \rightarrow (\SS, \TT)$ in $\fM^\bot\fM^\bot(\B)$ consists of a pair of downmaps $\FF = (F, d) : \TT' \rightarrow \TT$ and $\GG = (G, d') : \SS' \rightarrow \SS$ with $F = G$, and satisfying the equation
\[
\xymatrix@R=0.8em{
&& X \ar[drr]^{T'}_\;="2"&&\\
X \ar[rru]^{S'}_{\;}="1" &&&& X  \\
&& X  \ar[uu];[]^{F} \ar[drr]_{T}="4" && \\
X  \ar[urr]_{S}^\;="3" \ar[uu];[]_{F}  \ar[drr]_{T}&&&&X \ar[uu];[]^{F}  \\
&&X \ar[urr];[]^{S}  && \ar@{=>}  "1"+(0,-4.5); "1"+(0,-9.5) ^{d'} \ar@{=>} "2"+(0,-4.5) ; "2"+(0,-9.5)^{d} \ar@{=>}  "3"+(13,-3);  "3"+(13,-8)^c
}
\qquad
\xymatrix@C=1em{\ar@{}[ddd]|-{\hequ} \\ \\ \\  \\}
\qquad
\xymatrix@R=0.8em{
&& X  \ar[drr]^{T'}&&\\
X  \ar[rru]^{S'}  \ar[drr]^{T'}_\;="1"&&&& X  \\
&& X \ar[urr]^{S'}_\;="2"  && \\
X  \ar[uu];[]_{F} \ar[drr]_{T}^\;="3"&&&&X \ar[uu];[]^{F} \\
&&X. \ar[urr];[]^{S}_{\;}="4"   \ar[uu];[]^{F}&& \ar@{=>}"1"+(1,-4.5); "1"+(1,-9.5) ^{d} \ar@{=>}"2"+(0,-4.5); "2"+(0,-9.5) ^{d'} \ar@{=>} "1"+(13,8);  "1"+(13,3) ^{c'}
}
\]

\noindent A 2-cell $(\GG', \FF') \rightarrow (\GG, \FF)$ in $\fM^\top \fM^\top(\B)$ is a 2-cell $F' \rightarrow F$ which becomes a downmap transformation both, between $\FF'$ and $\FF$, and between $\GG'$ and $\GG$. 

The \textbf{composite of a distributive pair of monads} $(\SS, \TT)$, denoted $\comp(\SS, \TT)$, is defined to be the monad $(X, ST)$ with the multiplication: 
\[
\xymatrix{
\ar[d]^{T} X \ar@/_2.2pc/[dd]_{T}="1"\\
X \ar[d]^{T}="4" \ar[r]^{S} \ar@{=>} "1"+(8.5,0); "1"+(3.5,0) _m &X \ar[d]^{T}_{\;}="3" \ar@{=>} "3"+(-1.5,0); "3"+(-6.5,0) _c\\
X\ar[r]^{S}  \ar@/_2pc/[rr]_{S}^{\;}="5"& X \ar[r]^{S} \ar@{=>} "5"+(0,5);"5"^m &X,
}
\]
\noindent and the unit:
\[
\xymatrix{X \ar@/^1pc/[r]^{1_{X}}_{\quad}="1" \ar@/_1pc/[r]_{T}^{\quad}="2"  \ar@{=>}"1" ; "2" ^{e}& X \ar@/^1pc/[r]^{1_{X}}_{\quad}="3" \ar@/_1pc/[r]_{S}^{\quad} ="4" \ar@{=>}"3" ; "4" ^{e}&X.}
\]

\noindent There is a functor 
\[\xymatrix{\comp^\bot_\B : \fM^\bot\fM^\bot(\B)  \ar[r] & \fM^\bot(\B),}\]

\noindent defined on objects as the composite of distributive pairs of monads, and defined on morphisms and 2-cells by $\comp((F_0, h), (F_0, h')) = (F_0, (Sh)(h'T))$ and $\comp(t) = t$ respectively. We will not use this functor itself, but in Section~\ref{laxmonads} we will consider its lax generalization.

We conclude the section by a couple of simple definitions and a simple technical lemma. Suppose that $\TT = (X, T)$ is a monad. Define the \textbf{multiplication upmap} to be the upmap $\TT^\u = (T, m) : \one_X \rightarrow \TT$ consisting of the morphism $T : X \rightarrow X$ and the multiplication 2-cell of $\TT$
\[
\xymatrix{
X \ar[rr]^{1_{X}}_{\quad}="1"  \ar[d]_{T} && X \ar[d]^{T}\\
X \ar[rr]_{T}^{\quad}="2"&& X.
\ar@{=>} "2"+(0,3.5); "2"+(0,8.5) _m
}
\]

\noindent Analogously, define the \textbf{multiplication downmap} to be the downmap $\TT^\d = (T, m) : \TT \rightarrow \one_X$ consisting of the morphism $T : X \rightarrow X$ and the 2-cell 
\[
\xymatrix{
X \ar[rr]^{T}_{\quad}="1"  \ar[d]_{T} && A_0 \ar[d]^{T}\\
X\ar[rr]_{1_{X}}^{\quad}="2"&& X.
\ar@{=>} "2"+(0,8.5); "2"+(0,3.5) _m
}
\]

\begin{lemma}\label{lemma}
For any monad $\TT = (T, X)$, the multiplication downmap $\TT^\d$ and the multiplication upmap $\TT^\u$ determine a morphism in $\fM^\bot\fM^\top(\B)$:
\[(\TT^\u, \TT^\d) : (\TT, \one_X) \rightarrow (\one_X, \TT).\] 
\end{lemma}

\section{Equipments}\label{equipments}
\smallskip Informally, an equipment $\AA$ consists of objects, scalar arrows between objects, vector arrows between objects, and 2-cells between vector arrows, written respectively as:
\[
x 
\qquad\qquad 
\xymatrix{x \ar[r]^f& y} 
\qquad\qquad 
\xymatrix{x\ar|-{\object@{|}}[r]^a & y}  
\qquad\qquad 
\xymatrix{x  \ar@/^1pc/|-{\object@{|}}[r]^a_{\quad}="1" \ar@/_1pc/|-{\object@{|}}[r]_b^{\quad}="2"  \ar@{=>}"1" ; "2" ^{\alpha}& y.}
\]

\noindent Objects and scalar arrows form a category $A_0$. Objects, vectors and 2-cells form a lax bicategory $A$, meaning that, vectors and 2-cells between any fixed pair of objects $x$ and $y$ form a category $A(x, y)$; a multifold composite of vectors
\[\xymatrix{x_0 \ar|-{\object@{|}}[r]^{a_1} & x_1 \ar|-{\object@{|}}[r]^{a_2} &\cdots  \ar|-{\object@{|}}[r]^{a_n} & x_n,}\]

\noindent producing a vector $a_n\hdots a_2a_1 : x_0 \vectarr x_n$ is defined, and extends functorially to 2-cells; for each object $x$, there is a chosen identity vector $i_x : x \vectarr x$; and there are suitably coherent non-invertible 2-cells
\beq{xi}
\xymatrix{(a_{1,1} \hdots a_{1,n_1})(a_{2,1} \hdots a_{2,n_2})\hdots(a_{k,1} \hdots a_{k,n_k}) \ar@{=>}[r]^>>>>>{\xi} & (a_{1,1}\hdots a_{1,n_1} a_{21}\hdots  a_{kn_k}), 
}
\eeq

\noindent wherein strings of zero lengths should be interpreted as identity vectors. Furthermore, scalars act on vectors from left and right, i.e. diagrams such as
\[ 
\xymatrix{ w \ar[r]^f & x \ar[r]|-@{|}^a & y} \qquad\qquad \xymatrix{x \ar[r]|-@{|}^a & y \ar[r]^g & z}
\]

\noindent evaluate to vectors $af : w \vectarr y$ and $ga : x \vectarr z$ respectively; these actions are functorial in their vector argument, and associate with the vector multifold composition in various possible ways. Now we give a more formal definition:

Given categories $X$ and $Y$, by a \textbf{module} $A$ from $X$ to $Y$ we will mean a pseudofunctor $A : X^{\op} \times Y \rightarrow \cat$. The objects of $X$ and $Y$ will be called objects of the module, their morphisms will be called scalars of the module, the objects of $A(x, y)$ will be called vectors of the module, and its morphisms will be called 2-cells of the module.

\begin{definition} An \textbf{equipment} $\AA = (A_0, A, P^A, \xi^A)$ consists of the following data:
\begin{itemize}

\item A category $A_0$.  

\item A module $A$ from $A_0$ to itself.

\item For each $n > 1$, an $n$-fold vector composition $P^A_n$, which is a family of functors 
\begin{equation}\label{Pxxx}
P_{x_0, ..., x_n} : A(x_0, x_1)\times A(x_1, x_2)\times\cdots\times A(x_{n-1}, x_n)  \rightarrow A(x_0, x_n),
\end{equation}

\noindent  pseudonatural in $x_0$ and $x_n$, and pseudo-dinatural in $x_1, ..., x_{n-1}$. 

\item An identity $P^A_0$, which is a family of functors
\[P_x : I  \rightarrow A(x, x),\]  

\noindent where $I$ denotes the terminal category. 

\item A lax associator $\xi^A$, which is a specification for every partition $n = n_1 + n_2 + \cdots n_k$ of a modification with components natural transformations $\xi_{(x_{11}, ..., x_{1n_1}), ... (x_{k1}, ..., x_{kn_k})}$:
\[
\xy
(0,0)*{\big{(}A(x_{00}, x_{01})\times\cdots\times A(x_{0(n_0-1)}, x_{0n_0})\big{)}\times \cdots \times\big{(}A(x_{k0},  x_{k1})\times\cdots\times A(x_{k(n_k-1)}, x_{kn_k})\big{)}};
(-35,-20)*{A(x_{00}, x_{0n_0})\times\cdots\times A(x_{k0}, x_{kn_k})};
(35,-20)*{A(x_{00}, x_{01})\times\cdots\times A(x_{k(n_k-1)}, x_{kn_k})};
(0,-40)*{A(x_{00},\cdots x_{nn_k})};
{\ar_>>>>>>{P_{x_{00}, ..., x_{0n_0}}\times ... \times P_{x_{k0}, ..., x_{kn_k}}} (-5, -3); (-19, -17)}
{\ar^>>>>>>>>{\cong} (5, -3); (19, -17)}
{\ar_{P_{x_{00}, x_{0n_0}, ..., x_{kn_k}}} (-19, -23); (-5, -37)}
{\ar^{P_{x_{00}, x_{12}, ..., x_{kn_k}}} (16, -23); (1, -37)}
{\ar@{=>} (-7.5,-20); (-3.5,-20)}
\endxy
\]

\noindent satisfying the suitable coherence condition. Here $P_{x,y}$ is a family $P_1$ of the identical functors on $A(x, y)$. 
\end{itemize}
\end{definition}

\noindent Following the earlier informal description, instead of $P_{x_0, ..., x_n}(a_1, a_2, \hdots, a_n)$ we write $a_n\hdots a_2a_1$, and for the vector chosen by the functor $P_x$  we write $i_x$. In these notations, the components of $\xi_{(x_{11}, \hdots, x_{1n_1}), \hdots (x_{k1}, ..., x_{kn_k})}$ are 2-cells of $A$ of the form (\ref{xi}).

\begin{example} 
Our equipments are close to the ``virtual double categories with composites'' of \cite{CrSh10}. The latter work starts with a notion of a virtual double category which has a richer structure than our module. The vertical category of a virtual double category corresponds to our category of scalars, and its horizontal arrows correspond to our vectors. Composites of horizontal arrows are defined by a universal property. A virtual double category with specified choice of these composites can be regarded as our equipment. All the examples of virtual double categories considered in \cite{CrSh10} have composites, and hence can be regarded as our equipments.
\end{example}

\begin{example} Suppose that $J : A_0 \rightarrow A$ is a pseudofunctor from a category $A_0$ to a bicategory $A$. Then, the module
\[ 
\xymatrix{A_0^\op \times A_0 \ar[r]^{J \times J} & A^\op \times A \ar[r]^>>>>>{\mathrm{Hom}} & \cat}
\]
\noindent becomes an equipment with the vector composition structure induced by the horizontal composition of $A$ in the obvious way. The associator $\xi^A$ of this equipment is invertible. It is not difficult to see that all equipments with this property arise from a pseudofunctor from a category to a bicategory. 
\end{example}

\smallskip
\begin{definition}\label{eqfunct}
A \textbf{(lax) functor} $\FF = (F_0, F, \kappa^F) : \AA \rightarrow \BB$ between equipments, consists of the following data:

\begin{itemize}
\item A functor between the categories of scalars $F_0 : A_0 \rightarrow B_0$.

\item A family of functors between the vector categories
\begin{equation}\label{Fxy}
F_{x, y} :  A(x, y)  \rightarrow B(Fx, Fy),
\end{equation}

\noindent pseudonatural in both arguments.

\item A lax comparison structure $\kappa^F$, which consists of a family of natural transformations 
\begin{equation}\label{kappa}
\kappa^F_{x_0\cdots x_n} : P_{F_0x_0\cdots F_0x_n}\big{(}F_{x_0,x_1}\times F_{x_1,x_2}\times \cdots \times F_{x_{n-1},x_n}\big{)}  \rightarrow F_{x_0,x_n}P_{x_0\cdots x_n},
\end{equation}

\noindent for each sequence of objects $x_0, x_1,... x_n$ of $A$, compatible in the suitable sense with the vector multifold composition structures of $\AA$ and $\BB$. 
\end{itemize}
\end{definition}

\noindent For an object $x$, instead of $F_{0}(x)$ we write $Fx$, and for a vector $a : x \vectarr y$, instead of $F_{x, y}(a)$ we write $Fa$. Using this notations the components of $\kappa_{x_0,x_1\cdots x_n}$ are the 2-cells of $B$
\[\kappa^F_{a_1, ..., a_n} : Fa_1\cdots Fa_n \Rightarrow F(a_1\cdots a_n).\]

\noindent By changing the direction of $\kappa$ in this definition, we obtain a notion of \textbf{colax functor} between equipments. For the time being we will not use colax functors, so we keep the short name of a functor for lax functors. 

\begin{definition}\label{eqfuncttrans}
A \textbf{(lax) transformation} between functors of equipments $\tt = (t, \nu^t) : \FF \rightarrow \GG : \AA \rightarrow \BB$ consists of 
\begin{itemize}

\item A natural transformation $t : F_0 \rightarrow G_0 : A_0 \rightarrow B_0$

\item A modification with components the pseudonatural transformations
\begin{equation}\label{txy}
\nu^t_{x,y} : A(x, t_{y})F_{x, y}  \rightarrow A(t_{x}, y)G_{x, y} : A(x, y)  \rightarrow B(Fx, Gy),
\end{equation}

\noindent suitably compatible with the functor structures of $\FF$ and $\GG$.
\end{itemize}
\end{definition}

\noindent The components of the natural transformation $\nu^t_{x, y}$ are the 2-cells of $B$ 
\[ 
\xymatrix{
Fx \ar[rr]|-{\object@{|}}^{Fa}_{\quad}="1" \ar[d]_{t_x} && Fy \ar[d]^{t_y}\\
Gx \ar[rr]|-{\object@{|}}_{Ga}^{\quad}="2" &\ar@{=>}"1"+(0,-3.5); "1"+(0,-8.5) ^{\nu^t_a}& Gy.
}
\]

\noindent A \textbf{colax} transformation between lax functors of equipments is defined by reversing the direction of $\nu^t_{x, y}$ in the definition of a lax functor. Similarly, one can define lax and colax transformations between colax functors. However, for the time being we will only work with lax transformations between lax functors, hence we shortly refer to them as transformations of functors of equipments. 

There is an obvious way of defining a composition of functors of equipments, as well as of vertical and horizontal compositions of transformations of functors of equipments. Under these compositions, equipments, functors and transformations form a 2-category, which we denote by $\Eq$. There is a forgetful 2-functor $\Eq \rightarrow \cat$ acting on 0-, 1- and 2-cells as:
\[(A_0, A, P^A, \xi^A) \mapsto A_0,\]
\[(F_0, F, \kappa^F) \mapsto F_0,\]
\[(t_0, \nu^t) \mapsto t_0.\]

\section{T-equipments}\label{Tequipments}
\smallskip A $T$-equipment is a monad $(\AA, \TT) = (\AA, \TT, \mm, \ee)$ in $\Eq$. So defined, a $T$-equipment has an underlying monad $\TT_0 = (A_0, T_0)$ in $\cat$. We state the formal definition in a way that the roles of $\AA$ and $\TT_0$ appear more symmetric:

\begin{definition}\label{defteq}
A \textbf{T-equipment} $(\TT_0, \AA) = (\TT_0, \AA, T, \kappa^T, \nu^m, \nu^e)$ consists of an equipment $\AA$, a monad $\TT_0 = (T_0, m, e)$ on the category $A_0$, and a lifting of this monad in $\cat$ to a monad $(\TT, \mm, \ee)$ in $\Eq$ on the object $\AA$; with $\TT = (T_0, T, \kappa^T)$, $\mm = (m, \nu^m)$ and $\ee = (e, \nu^e)$.
\end{definition}

\noindent Thus, the data of the $T$-equipment $(\TT_0, \AA, T, \kappa^T, \nu^m, \nu^e)$, besides the equipment $\AA$, and the monad $\TT_0$ consists of:

\begin{itemize}
\item A family of functors $T_{x, y} : A(x, y) \rightarrow A(Tx, Ty)$.

\item A family of 2-cells $\kappa_{a_1, ..., a_n} : Ta_1\cdots Ta_n \Rightarrow T(a_1\cdots a_n).$

\item Two families of 2-cells
\begin{equation}\label{em}
\xymatrix{
T^2x \ar[rr]|-{\object@{|}}^{T^2a}_{\quad}="1" \ar[d]_{m_x} && T^2y \ar[d]^{m_y}\\
Tx \ar[rr]|-{\object@{|}}_{Ta}^{\quad}="2" &\ar@{=>}"2"+(0,8.5); "2"+(0,3.5) ^{\nu^m_a}& Ty 
}
\qquad
\qquad
\xymatrix{
x \ar[rr]|-{\object@{|}}^{a}_{\quad}="1" \ar[d]_{e_x} && y \ar[d]^{e_y}\\
Tx \ar[rr]|-{\object@{|}}_{Ta}^{\quad}="2" &\ar@{=>}"2"+(0,8.5); "2"+(0,3.5) ^{\nu^e_a}& Ty.
}
\end{equation}
\end{itemize}

\noindent A set of axioms satisfied by this data can be extracted from the definition with some effort. 

\begin{example} 
We already noticed that our equipment is an analogue of the virtual double category with composites of \cite{CrSh10}. Hence a monad on a virtual double category with composites is an analogue of a $T$-equipment. Correspondingly, numerous examples listed in Table 1 in \cite{CrSh10} are examples of $T$-equipments.
\end{example}

\begin{example} We consider one specific situation, and briefly review a couple of its sub-examples. Let $V$ be a monoidal category with coproducts which distribute over the monoidal product. Recall that the bicategory of matrices $\mat(V)$ has small sets as its objects, while for sets $X$ and $Y$, $\mat(V)(X, Y)$ is the category $[X \times Y, V]$. So, a morphism $X \rightarrow Y$ of $\mat(V)$ is a family of objects $a_{x,y}$ of $V$ indexed by elements of the set $X\times Y$. The horizontal composition is by the usual matrix composition formula
\[(a\circ b)_{x, y} = \coprod_z a_{x, z}\otimes b_{z, y}.\]

\noindent The identity morphisms are the matrices with the monoidal unit at the diagonal and the initial object everywhere else. There is a pseudofunctor $\set \rightarrow \mat(V)$, which takes a set map $f : X \rightarrow Y$ to a matrix whose components are the monoidal unit on pairs of the form $(x, f(x))$, and the initial object otherwise. Corresponding to this pseudofunctor, there is an equipment $\MMatV = (\set, \mat(V))$.

A $T$-equipment $(\TT_0, \MMatV)$ is the same as the ``monad $\TT_0$ with a lax extension $\TT$ to $\mat(V)$'' of \cite{ClTh03}. By varying the monad $\TT_0$ on $\set$ and the monoidal category $V$, together with its lax extension $\TT$ to $\mat(V)$, we obtain various examples. Later, we will return to two of them: In the first case, $\TT_0$ is taken to be the free-monoid monad on $\set$, and $V$ is any monoidal category. In the second case, $\TT_0$ is taken to be the ultrafilter monad on $\set$, and $V$ is taken to be the lattice $2$. The construction of the lax extension $\TT$ in both of these cases can be found in \cite{ClTh03}.
\end{example}

\smallskip Below we will introduce several 2-categories whose objects are $T$-equipments. As the first example, we consider the 2-category $\fM^\top(\Eq)$. Objects of $\fM^\top(\Eq)$ are $T$-equipments since they are monads in $\Eq$. Let us identify the morphisms and the 2-cells. In the notation of Section~\ref{monads}, a downmap between monads $(\BB, \SS) \rightarrow (\AA, \TT)$ in $\Eq$ is a pair $(\FF, \dd)$ where $\FF = (F_0, F) : \AA \rightarrow \BB$ is an equipment functor and $\dd = (d, \nu^d) : \FF\SS \rightarrow \TT\FF$ is a transformation of functors of equipments. The following is a presentation of the same morphism in line with Definition~\ref{defteq}. 

A morphism  in $\fM^\top(\Eq)$ between $T$-equipments $(\SS_0, \BB) \rightarrow (\TT_0, \AA)$ is a triple $(\FF_0, \FF, \nu^d)$ consisting of

\begin{itemize}
\item A functor of equipments $\FF = (F_0, F, \kappa^F) : \BB \rightarrow \AA$.

\item A downmap of monads $\FF_0 = (F_0, d) : \SS_0 \rightarrow \TT_0$ in $\cat$, consisting of the functor $F_0 : B_0 \rightarrow A_0$ and a natural transformation $d : F_0S_0 \rightarrow T_0F_0$.

\item A modification $\nu^d$ with the components natural transformations 
\[\xymatrix{A(x, y) \ar[rr]^{F_{}S_{x,y}}_{\;}="1" \ar[d]_{T_{}F_{x,y}}&& B(FSx, FSy) \ar[d]^{B(FSx, d_y)}\\
B(TFx, TFy) \ar[rr]_{B(d_x, TFy)}^{\;}="2"  && B(FSx, TFy) \ar@{=>} "1"+(0,-3.5); "1"+(0,-8.5) ^{\nu^d_{x,y}}
}\]
\end{itemize}

\noindent satisfying a few axioms. Observe that, the components of $\nu^d_{x,y}$ are 2-cells of $B$
\[ 
\xymatrix{
FSx \ar[rr]|-{\object@{|}}^{FSa}_{\quad}="1" \ar[d]_{d_x} && FSx \ar[d]^{d_y}\\
TFx \ar[rr]|-{\object@{|}}_{TFa}^{\quad}="2" &\ar@{=>}"1"+(0,-3.5); "1"+(0,-8.5) ^{\nu^d_a}& TFy. 
}
\]
 
\noindent A 2-cell  $(\FF, \FF_0, \nu^d) \rightarrow (\GG, \GG_0, \nu^{d'}) : (\SS_0, \BB) \rightarrow (\TT_0, \AA)$ in $\fM^\top(\Eq)$ amounts to a transformation of functors of equipments $(t, \nu^t) : \FF \rightarrow \GG$, such that $t_0 : F_0 \rightarrow G_0$ is a downmap transformation, and a certain additional axiom expressing compatibility of $\nu^t$ with $\nu^d$ and $\nu^{d'}$ is satisfied.

Define now another 2-category $\overline{\fM^\top(\Eq)}$ whose objects are $T$-equipments. A morphism $(\BB, \SS_0) \rightarrow (\AA_0, \TT_0)$ in it is defined to be a triple $(\FF, \FF_0, \overline{\nu}^d)$ where $\FF$ and $\FF_0$ are as in a morphism of $\fM^\top(\Eq)$, while $\overline{\nu}^d$ takes the opposite direction to $\nu^d$, that is it is a modification with the components natural transformations
\[\xymatrix{A(x, y) \ar[rr]^{F_{}S_{x,y}}_{\;}="1" \ar[d]_{T_{}F_{x,y}}&& B(FSx, FSy) \ar[d]^{B(FSx, d_y)}\\
B(TFx, TFy) \ar[rr]_{B(d_x, TFy)}^{\;}="2"  && B(FSx, TFy). \ar@{=>} "2"+(0,3.5); "2"+(0,8.5) _{\overline{\nu}^d_{x,y}}
}\]
 
\noindent The axioms which $\overline{\nu}^d$ should satisfy are obtained from the equations satisfied by the $\nu^d$ component of a morphism of $\fM^\top(\Eq)$ by replacing in them all arrows involving $\nu^d$ by the oppositely directed arrows involving $\overline{\nu}^d$ (commutativity of which still makes sense). The components of $\bar\nu^d$ are 2-cells of $B$
\[ 
\xymatrix{
FSx \ar[rr]|-{\object@{|}}^{FSa}_{\quad}="1" \ar[d]_{d_x} && FSx \ar[d]^{d_y}\\
TFx \ar[rr]|-{\object@{|}}_{TFa}^{\quad}="2" &\ar@{=>}"1"+(0,-8.5); "1"+(0,-3.5) ^{\overline\nu^d_a}& TFy. 
}
\]

\noindent A 2-cell $(\FF, \FF_0, \overline{\nu}^d) \rightarrow (\GG, \GG_0, \overline{\nu}^{d'}) : (\BB, \SS_0) \rightarrow (\AA, \TT_0)$ of $\overline{\fM^\top(\Eq)}$ is defined to be a transformation of functors of equipments $(t, \nu^t) : \FF \rightarrow \GG$, such that $t : F_0 \rightarrow G_0$ is a monad downmap transformation, and a certain additional axiom expressing compatibility of $\nu^t$ with $\overline{\nu}^d$ and $\overline{\nu}^{d'}$ is satisfied.

\smallskip The following construction on $T$-equipments is of principle interest to us.
  
\begin{definition}
The \textbf{Kleisli equipment} of a $T$-equipment $(\TT_0, \AA)$, denoted $\kl(\TT_0, \AA)$, is defined to be an equipment consisting of

\begin{itemize}
\item The category of scalars $A_0$. 

\item The module of vectors $A(-, T-)$, i.e. the pseudofunctor
\[\xymatrix{A_0^{\op}\times A_0 \ar[r]^{1\times T} & A_0^{\op}\times A_0 \ar[r]^{\mathrm{Hom}} &\cat.}\]

\item The $n$-fold composition of vectors defined by
\[\xymatrix{
A(x_0, Tx_1)\times A(x_1, Tx_2)\times\cdots\times A(x_{n-1}, Tx_n) \ar[d]^{1\times T\times\cdots\times T^{n-1}}\\
 A(x_0, Tx_1)\times A(Tx_1, T^2x_2)\times\cdots\times A(T^{n-1}x_{n-1}, T^nx_n) \ar[d]^{P_{x_0,\cdots x_n}}\\
 A(x_0, T^nx_n) \ar[d]^{A(x_0, (m_n)_x)} \\
 A(x_0, Tx_n).
}
\]

\item The identity vectors defined by
\[\xymatrix{I \ar[r]^{P_x} & A(x, x) \ar[r]^{A(x, e_x)} & A(x, Tx)}\] 

\item The lax associator defined from the associator $\xi^A$, using $\nu^m$, $\nu^e$ and $\kappa^T$ (see below).

\end{itemize}
\end{definition}

\noindent A vector in the Kleisli equipment from $x$ to $y$ is a vector $x \vectarr Ty$ of $A$. An $n$-fold composite of Kleisli vectors is formed as a composite of vectors of $A$:
\[\xymatrix{x_0 \ar|-{\object@{|}}[r]^{a_1} & Tx_1 \ar|-{\object@{|}}[r]^{Ta_2} &\cdots  \ar|-{\object@{|}}[r]^>>>>>{T^{n-1}a_n} & T^nx_n \ar[r]^{(m_n)_x} & Tx.}\]

\noindent The identity Kleisli vectors are
\[\xymatrix{x  \ar|-{\object@{|}}[r]^{i_x} & x \ar[r]^{e_x} & Tx}.\]

\noindent The components of the components of the Kleisli lax associator can be defined by hand as follows. For the partitions  $0+1$, $1+ 0$, $2+1$ and $1+2$, they are given by the diagrams:
\[ 
\xymatrix{
&&&&Tx \ar[rrd]^{e_{Ty}}="2"  \ar@{=>} "2"+(-15,0); "2"+(-15,-5) ^{\nu^e_a}&&&&\\
x \ar|-{\object@{|}}@/^1.2pc/[rrrru]^{a}="1"  \ar|-{\object@{|}}[rr]_{i_{x}} && x \ar|-{\object@{|}}[rru]^{a} \ar[rr]_{e_x} \ar@{=>} "1"+(0,-5);  "1"+(0,-10) ^{\xi^A_{(0+1)-}}&& Tx \ar|-{\object@{|}}[rr]_{Ta} && T^2y \ar[rr]_{m_y} && Ty}
\]
\[ 
\xymatrix{
x \ar|-{\object@{|}}@/^2.5pc/[rrrr]^{a}="1" \ar|-{\object@{|}}[rr]_{a} && Tx \ar@{=>} "1"+(0,-4);  "1"+(0,-9) ^{\xi^A_{(1+0)-}} \ar|-{\object@{|}}@/_1.8pc/[rr]_{Ti_x}="3" \ar|-{\object@{|}}[rr]^{i_{Tx}}="2" && Tx \ar[rr]_{Te_x} && T^2y \ar[rr]_{m_y} && Ty \ar@{=>} "2"+(0,-3.5); "2"+(0, -7.5) ^{\kappa^F_{x}}}
\]
\[ 
\xymatrix{
&&&&&&T^3y \ar[rrd]^{m_{Ty}}="2" &&&&\\
v \ar|-{\object@{|}}@/^1.2pc/[rrrrrru]^{T^2cT(b)c}="1" \ar|-{\object@{|}}[rrrr]_{T(b)c} &&&& T^2x \ar|-{\object@{|}}[rru]^{T^2a} \ar[rr]_{m_x} \ar@{=>} "1"+(5, -5.5); "1"+(5, -10.5)^{\xi^A_{(2+1)-}}&& Tx \ar|-{\object@{|}}[rr]_{Ta} && T^2y \ar[rr]_{m_y} && Ty \ar@{=>} "2"+(-15,0); "2"+(-15,-5) ^{\nu^m_a}
}
\]
\[ 
\xymatrix{v \ar|-{\object@{|}}@/^3pc/[rrrrrr]^{T^2aT(b)c}="1"  \ar|-{\object@{|}}[rr]_{c} && Tw \ar|-{\object@{|}}@/_2.5pc/[rrrr]_{T(T(a)b)}^{\;}="3"  \ar|-{\object@{|}}[rrrr]^{T^2aTb}="2" && && T^3y \ar@<-1ex>[rr]_{m_{Ty}} \ar@<1ex>[rr]^{Tm_y} && T^2y \ar[rr]_{m_y} && Ty \ar@{=>} "1"+(-1,-6.5); "1"+(-1,-11.5) ^{\xi^A_{(1+2)-}} \ar@{=>} "2"+(0,-5);  "2"+(0,-10) ^{\kappa^T_{Ta, b}}
}
\]

\noindent where we have ignored structural isomorphisms for scalar actions on vectors. For higher partitions the components of the Kleisli associator involve more complex diagrams, writing out which, albeit requiring some effort, is fairly straightforward. A somewhat more conceptual perspective will be subsequently provided in Section~\ref{laxmonads}, where the Kleisli construction will be observed to be a case of the formal theory developed there.

\begin{example} The Kleisli equipment construction is analogous to the construction of the Kleisli virtual double category of \cite{CrSh10}.
\end{example}

\begin{example}
The composition of vectors in the Kleisli equipment $\kl(\TT_0, \MMatV)$ is the ``Kleisli convolution'' of $(T, V)$-relations.
\end{example}

\begin{remark} Suppose that $\AA$ is an equipment with an invertible associator (that is, it is an equipment coming from a pseudofunctor $A_0 \rightarrow A$). Observe that, given a $T$-equipment $(\TT_0, \AA)$, the lax associator of the Kleisli equipment $\kl(\TT_0, \AA)$ is no longer invertible. In examples the initial input $T$-equipments are usually equipments with an invertible associator. The Kleisli construction however takes us out of this situation. 
\end{remark}

\smallskip A functorial extension of the Kleisi construction exists to the 2-category $\overline{\fM^\top(\Eq)}$. Suppose that $(\FF_0, \FF, \overline{\nu}^d) : (\SS_0, \BB) \rightarrow (\TT_0, \AA)$ is a morphism in $\overline{\fM^\top(\Eq)}$. Define a functor $\kl(\FF_0, \FF) : \kl(\SS_0, \BB) \rightarrow \kl(\TT_0, \AA)$ between the Kleisli equipments to have the following component:

\begin{itemize}
\item The functor of scalars $F_0 : A_0 \rightarrow A_0$.

\item The functors between the categories of vectors defined as
\[
\xymatrix{
B(x, Sy) \ar[r]^{F_{x, y}}& A(Fx, FSy) \ar[r]^{A(Fx, d_y)}& A(Fx, TFy).
}
\]

\item The lax comparison structure given by a family of natural transformations
{\small
\beq{Bx0Sx1}
\xy
(0,100)*{B(x_0, Sx_1)\times B(x_1, Sx_2)\times \cdots \times B(x_{n-1}, Sx_n)}="1";
(-35,90)*{B(Fx_0, FSx_1)\times \cdots \times B(Fx_{n-1}, FSx_n)}= "2";
(35,90)*{B(x_0, Sx_1)\times \cdots \times B(S^{n-1}x_{n-1}, S^nx_n)}="3";
(-35,80)*{B(Fx_0, TFx_1)\times \cdots \times B(Fx_{n-1}, TFx_n)}="4";
(35,80)*{B(x_0, S^nx_n)}="5";
(-35,70)*{B(Fx_0, TFx_1)\times \cdots \times B(T^{n-1}Fx_{n-1}, T^nFx_n)}="6";
(35,70)*{B(x_0, Sx_n)}="7";
(-35,60)*{B(Fx_0, T^nFx_n)}="8";
(35,60)*{B(Fx_0, FSx_n)}="9";
(0,50)*{B(Fx_0, TFx_n)}="10";
{\ar _{F\times F\times\cdots \times F} "1"; "2"} 
{\ar ^{1\times S\times\cdots S^{n-1}} "1"; "3"}
{\ar _{(d-)\times(d-)\times\cdots\times(d-)} "2"; "4"}
{\ar ^{P_n} "3"; "5"}
{\ar _{1\times T\times\cdots T^{n-1}} "4"+(-2,-2); "6"+(-2,2)}
{\ar ^{m-} "5"+(-2,-2); "7"+(-2,2)}
{\ar _{P_n} "6"+(-4,-2); "8"+(-4,2)}
{\ar ^{F} "7"+(-4,-2); "9"+(-4,2)}
{\ar _{m-} "8"+(-4,-2); "10"+(-6,2)}
{\ar ^{(d-)} "9"+(-8,-2); "10"+(-4,2)}
{\ar@{=>} (-10,77); (-5,77)}
\endxy
\eeq
}

\noindent defined using $\overline{\nu}^d$ and the lax comparison structure $\kappa^F$, as outlined on components below.

\end{itemize}

\noindent For $n = 0$ the components of (\ref{Bx0Sx1}) are
\[\xymatrix{Fx \ar|-{\object@{|}}[rr]_{Fi_x}^{\;}="2" \ar|-{\object@{|}}@/^2pc/[rr]^{i_{Fx}}_{\;}="1" && Fx \ar[r]_{Fe_x} \ar@/^2pc/[rr]^{e_{Fx}}& FTx \ar[r]_{d} & SFx. \ar@{=>} "1"+(0,-1); "1"+(0,-5) ^{\kappa^F_x}}\]

\noindent  For $n = 2$, they are
\[ 
\xymatrix{
Fw  \ar|-{\object@{|}}@/_1.3pc/[rrrd]_{F(T(a)b)}="1" \ar|-{\object@{|}}[r]^{Fb} & FTx \ar[rr]^{d} \ar|-{\object@{|}}[rrd]^{FTa} \ar@{=>} "1"+(-2,11) ; "1"+(-2,6) ^{\kappa^F_{a,Ta}}&& SFx \ar|-{\object@{|}}[rr]^{SFa} \ar@{=>} {}+(0,-4);{}+(0,-9) ^{\nu_a^d}&& SFTy \ar[r]^{Sd} & S^2Ty \ar[r]^{m_{Ty}} & STy\\
&&&FT^2y \ar[rru]^{d_{Ty}} \ar[rr]_{Fm_{y}}&&FTy. \ar[rru]_{d}\\
}
\]

\noindent For higher $n$, a similar description works. For example, for $n = 3$, the component are
\[ 
\xymatrix@C=2.2em{
Fv \ar|-{\object@{|}}@/_2.5pc/[rrrrrdd]_{F(T^2(a)T(b)c)}="1" \ar|-{\object@{|}}[r]^{Fc} & FTw \ar|-{\object@{|}}[rrd]^{FTb} \ar[rr]^{d_w} \ar@{=>}{}+(7,-9);{}+(7,-14) ^{\kappa^F_{c,Tb,T^2a}}&& SFw \ar|-{\object@{|}}[rr]^{SFb} \ar@{=>}{}+(0,-4); {}+(0,-9)^{\nu_b^d}&& SFTx \ar[rr]^{Sd_x} \ar|-{\object@{|}}[rrd]_{SFTa} \ar@{=>}{}+(0,-11);{}+(0,-16)^{\nu^d_{Ta}}&& S^2Fx \ar|-{\object@{|}}[rr]^{S^2Fa} \ar@{=>} {}+(0,-4); {}+(0,-9)^{S\nu_a^{d}} && S^2FTy \ar[d]^{S^2d_y} &&&& \\
&&&FT^2x \ar[rru]_{d_{Tx}} \ar|-{\object@{|}}[rrd]^{FT^2a} &&&&SFT^2y \ar[rru]_{Sd_{Ty}}&&S^3y \ar[d]^{m_{3y}}&& \\
&&&&& FT^3y \ar[rru]_{d_{T^2y}} \ar[rr]_{Fm_{3y}}&& FTy \ar[rr]_{d} &&SFy. 
}
\]

\noindent We leave it to the reader to define $\kl$ on the 2-cells of $\fM^\top(\Eq)$, and to conclude that $\kl$ is a 2-functor

\[\xymatrix{\kl : \overline{\fM^\top(\Eq)} \ar[r]&\Eq.}\]

\smallskip Let us also characterize the 2-category $\fM^\bot(\Eq)$. This is another 2-category whose objects are $T$-equipments. A morphism $(\SS_0, \BB) \rightarrow (\TT_0, \AA)$ in it is a triple $(\FF_0, \FF, \nu^u)$ consisting of

\begin{itemize}
\item A functor of equipments $\FF = (F_0, F, \kappa^F) : \BB \rightarrow \AA$.
\item An upmap of monads $\FF_0 = (F_0, u) : \SS_0 \rightarrow \TT_0$ in $\cat$, consisting of the functor $F_0 : B_0 \rightarrow A_0$ and a natural transformation $u : T_0F_0 \rightarrow F_0S_0$.
\item A modification $\nu^u$ with the components natural transformations 
\[\xymatrix{A(x, y) \ar[rr]^{T_{}F_{x,y}}_{\;}="1" \ar[d]_{F_{}S_{x,y}}&& B(TFx, TFy) \ar[d]^{B(x, u_y)}\\
B(FSx, FSy) \ar[rr]_{B(u_x,  FSy)}^{\;}="2" && B(TFx, FSy) \ar@{=>} "2"+(0,3.5); "2"+(0,8.5) _{\nu^u_{x,y}}
}\]

\end{itemize}

\noindent satisfying a few equations. The component of $\nu^u_{x, y}$ are 2-cells of $B$
\[ 
\xymatrix{
FSx \ar[rr]|-{\object@{|}}^{FSa}_{\quad}="1" \ar[d]_{d_x} && FSx \ar[d]^{d_y}\\
TFx \ar[rr]|-{\object@{|}}_{TFa}^{\quad}="2" &\ar@{=>}"1"+(0,-3.5); "1"+(0,-8.5) ^{\nu^u_a}& TFy. 
}
\]

\noindent A 2-cell $(\FF_0, \FF, \nu^u) \rightarrow (\GG_0, \GG, \nu^{u'}) : (\BB, \SS_0) \rightarrow (\AA, \TT_0)$ amounts to a transformation of functors of equipments $(t, \nu^t) : \FF \rightarrow \GG$, such that $t : F_0 \rightarrow G_0$ is a monad upmap transformation, and a certain additional axiom expressing compatibility of $\nu^t$ with $\nu^u$ and $\nu^{u'}$ is satisfied.

Finally, there is a 2-category $\overline{\fM^\bot(\Eq)}$ whose morphisms are like morphisms of $\fM^\top(\Eq)$ except that their $\nu^u$ component takes the opposite direction.

\section{Monoids in an equipment}\label{monoidsinequ} 

\smallskip Slightly changing the previous notation, let $I_0$ stand for the terminal category, and let $\II = (I_0, I)$ be the terminal equipment, its module of vectors $I$ being the constant pseudofunctor $I_0\times I_0 \rightarrow \cat$ at the terminal category.

\begin{definition}\label{monoid}
The category of monoids $\mon(\AA)$ in an equipment $\AA$ is by definition the category $\Eq(\II, \AA)$;  its objects are called \textbf{monoids}, and its morphisms are called \textbf{monoid homomorphisms}.
\end{definition}

\noindent A monoid amounts to a data $(x, a, \mu_a, \eta_a)$, where $x$ is an object of $A$, $a : x \vectarr x$ is a vector, and $\mu_a$ and $\eta_a$ are 2-cells
\[
\xymatrix{
&x \ar[rd]|-{\object@{|}}^{a} \ar@{=>}{}+(0,-4.5); {}+(0,-9.5) ^{\mu_a}  & \\
x \ar[ru]|-{\object@{|}}^{a} \ar[rr]|-{\object@{|}}_{a}^{\quad}="2"&& x 
}
\qquad
\qquad
\qquad
\xymatrix{&&\\
x  \ar@/^2pc/|-{\object@{|}}[rr]^{i_x}_{\quad}="1" \ar|-{\object@{|}}[rr]_a^{\quad}="2" & \ar@{=>}"1" ; "1"+(0,-5) ^{\eta_a}& x,}
\]

\noindent satisfying associativity and unitivity axioms. Indeed, suppose that $\FF : \II \rightarrow \AA$ is an equipment functor. Its scalar functor $F_0 : I_0 \rightarrow A_0$ is determined by the object $F(\emptyset
) = x$ of $A$. The vector functor $F_{\emptyset, \emptyset} : \{\emptyset\} =  I(\emptyset, \emptyset) \rightarrow A(F(\emptyset), F(\emptyset)) = A(x, x)$ is determined by the vector $a = F(\emptyset) : x \vectarr x$. The 2-cells $\mu_a$ and $\eta_a$ are determined by $\kappa_{\emptyset, \emptyset}$ and $\kappa_0$ respectively. Since $\II$ is an equipment with an invertible associator, the comparison structure $\kappa^F$ is completely determined by these two. Similarly, it can be easily seen that a monoid homomorphism $(x, a) \rightarrow (y, b)$, amounts to a pair $(f, \phi_f)$, where $f : x \rightarrow y$ is a scalar and $\phi_f$ is a 2-cell 
\[
\xymatrix{
x \ar[rr]|-{\object@{|}}^{a}_{\quad}="1" \ar[d]_{f} && x \ar[d]^{f}\\
y \ar[rr]|-{\object@{|}}_{b}^{\quad}="2" &\ar@{=>}"1"+(0,-3.5); "1"+(0,-8.5) ^{\phi_f}& y 
}
\]

\noindent satisfying two axioms. Immediately from the definition it follows that taking the category of monoids is a representable 2-functor:
\[\xymatrix{\mon(-) =  \Eq(\II, -) : \Eq \ar[r] & \cat.}\]

\begin{definition}\label{tmonoid}
The category of $T$-monoids $\tmon(\TT_0, \AA)$ in a $T$-equipment $(\TT_0, \AA)$ is by definition the category of monoids in the Kleisli equipment $\mon(\kl(\TT_0, \AA))$; its objects are called \textbf{T-monoids}, and its morphisms are called \textbf{$T$-monoid homomorphisms}.
\end{definition}

\noindent A $T$-monoid consists of a data $(x, a, \mu_a, \eta_a)$, where $x$ is an object of $A$, $a$ is a Kleisli vector $x \vectarr Tx$ and $\mu_a$ and $\eta_a$ are the 2-cells:
\[
\xymatrix{
x \ar[rr]|-{\object@{|}}^{a} \ar@/_1pc/[rrrrd]|-{\object@{|}}_{a}&& Tx \ar[rr]|-{\object@{|}}^{Ta}_{\quad}="1" && T^2x \ar[d]^{m_x}\\
&&\ar@{}[rr]^{\quad}="2"&& Tx \ar@{=>}"1"+(-5,-2.5); "1"+(-5,-7.5) ^{\mu_a}
}
\qquad
\qquad
\xymatrix{&x \ar[rd]^{e_x} \ar@{=>} {}+(0,-5); {}+(0,-10) ^{\eta_a}&\\
x  \ar|-{\object@{|}}[ru]^{i_x}_{\quad}="1" \ar|-{\object@{|}}[rr]_a^{\quad}="2" && Tx,}
\]

\noindent satisfying three axioms. A $T$-monoid homomorphism $(x, a) \rightarrow (x, b)$ is a pair $(f, \phi_f)$ consisting of a scalar $f : x \rightarrow y$ and a 2-cell
\[
\xymatrix{
x \ar[rr]|-{\object@{|}}^{a}_{\quad}="1" \ar[d]_{f} && Tx \ar[d]^{Tf}\\
y \ar[rr]|-{\object@{|}}_{b}^{\quad}="2" &\ar@{=>}"1"+(0,-3.5); "1"+(0,-8,5) ^{\phi_f}& Ty 
}
\]

\noindent satisfying two axioms. Since $\kl$ has a functorial extension to $\overline{\fM^\top(\Eq)}$, it follows that taking  $T$-monoids is a 2-functor:
\[\xymatrix{\tmon(-) = \Eq(\II, \kl(-)) : \overline{\fM^\top(\Eq)} \ar[r]&\cat.}\]

\begin{example} Our monoids and $T$-monoids are essentially the same as monoids and $T$-monoids of \cite{CrSh10}. Numerous examples can be found listed in Table 1 there.
\end{example}

\begin{example} A monoid in $\MMatV$ is a $V$-category. A $T$-monoid in $(\TT_0, \MMatV)$ is  a $(T, V)$-category introduced in \cite{ClTh03}. In particular: When $\TT_0$ is the free-monoid monad and $V$ is an arbitrary monoidal category, a $T$-monoid is a $V$-multicategory. When $\TT_0$ is the ultrafilter monad and $V$ is the lattice $2$, a $T$-monoid is a topological space. 
\end{example}

\smallskip The trivial monad $\id(I_0) = \one_{I_0} = (I_0, 1_{I_0})$ on the terminal category $I_0$ lifts to the trivial monad $\id(\II) = \one_\II$ on the terminal equipment $\II$, giving the $T$-equipment $(\one_{I_0}, \II)$, which is the terminal object in $\fM^\top(\Eq)$.

\begin{definition}\label{talgebra}
The category of $T$-algebras $\talg(\TT_0, \AA)$ in a $T$-equipment $(\TT_0, \AA)$ is by definition the category $\fM^\bot(\Eq)\big{(}(\one_{I_0}, \II), (\AA, \TT_0)\big{)}$; its objects are called \textbf{$T$-algebras} and its morphisms are called \textbf{$T$-algebra homomorphisms}.
\end{definition}

\noindent A $T$-algebra amounts to a monoid $(x, a)$ in $\AA$, together with a scalar $h : Tx \rightarrow x$ and a 2-cell

\[ 
\xymatrix{
Tx \ar[rr]|-{\object@{|}}^{Ta}_{\quad}="1" \ar[d]_{h} && Tx \ar[d]^{h}\\
x \ar[rr]|-{\object@{|}}_{a}^{\quad}="2" &\ar@{=>}"1"+(0,-3.5); "1"+(0,-8.5) ^{\sigma_h}& x 
}
\]

\noindent such that $(x, h)$ is an algebra for the monad $\TT_0$, and $\sigma_h$ satisfies two axioms. Indeed, given a morphism $(\FF_0, \FF, \nu^u) : (\one_{I_0}, \II) \rightarrow (\TT_0, \AA)$ in $\fM^\bot(\Eq)$, the equipment functor $\FF : \II \rightarrow \AA$ amounts to a monoid $(x, a)$ in $\AA$, the monad downmap  $\FF_0 : \one_{I_0} \rightarrow \TT_0$, as observed in Section~\ref{monads}, amounts to an algebra $(x, h)$ of $\TT_0$, while $\sigma_a$ is determined by $\nu^u_{\emptyset, \emptyset}$. A $T$-algebra homomorphism is a monoid map $(f, \phi_f)$ compatible with the algebra structures of the source and the target. Immediately from the definition it follows that taking $T$-algebras is a representable 2-functor:
\[\xymatrix{\talg(-) = \fM^{\bot}(\Eq)\big{(}(\one_{I_0}, \II), -\big{)}: \fM^{\bot}(\Eq) \ar[r]& \cat.}\] 

Given a $T$-equipment $(\TT_0, \AA)$, its defining monad $(\AA, \TT)$ in $\Eq$  is taken by the 2-functor $\mon$ to a monad $(\mon(\AA), \mon(\TT))$  in $\cat$. Let  $\T$ denote the functor $\mon(\TT)$
\[\xymatrix{\T : \mon(\AA) \ar[r] & \mon(\AA).}\]

\noindent Explicitly, the image of a monoid under $\T$ is given by the formula

\[\T(x, a, \mu_a, \eta_a) = (Tx, Ta, T\mu_a\kappa_{a, a}, T\eta_a\kappa_x).\]

 \noindent The monad multiplication $\T^2 \rightarrow \T$ and the unit $1_\T \rightarrow \T$ are the natural transformations with the components on a monoid $(x, a)$ respectively the monoid homomorphisms $(m_x, \nu_a^m)$ and $(e_x, \nu_a^e)$.

The category of $T$-algebras $\talg(\TT_0, \AA)$ is by definition the category of modules $\fM^\bot(\Eq)(\one_\II, (\AA, \TT))$ of the monad $(\AA, \TT)$ with the fixed underlying object $\II$.  So, by the observation made in Section~\ref{monads}, it is the same as the category of EM algebras $\mon(\AA)^\T$ for the monad $\T = \mon(\TT) = \Eq(\II, \TT)$ on the category $\mon(\AA) = \Eq(\II, \AA)$. Consequently, we have a diagram in $\cat$ 
\[\xymatrix@C=0.1em{
&\talg(\TT_0, \AA) \ar[dl]_{}="1" \ar[dr]^{}="2"&\\
\mon(\AA)  \ar[rr]_{\T} && \mon(\AA) \ar@{=>} "1"+(8,-1); "1"+(13,-1)
}\]

\noindent exhibiting $\talg(\TT_0, \AA)$ as the EM category. Thus, a $T$-algebra can be alternatively defined as an algebra of the monad $\T$, from which point of view, it  consists of a monoid $(x, a, \mu_a, \eta_a)$ in $\AA$ and an algebra structure $(h, \sigma_h) : \T(x, a) \rightarrow (x, a)$.

\begin{example} In the case $\AA = \MMatV$, $T$-algebras are exactly the $T$-algebras considered in \cite{ChClHo14}. In particular: When $\TT_0$ is the free monoid monad, and $V$ is an arbitrary monoidal category, a $T$-algebra is a strict monoidal $V$-category. When $\TT_0$ is an ultrafilter monad, and $V$ is the lattice 2, a $T$-algebra is an ordered Compact Hausdorff space.

When $\AA = \MMatV$, and $\TT_0$ is the free monoid monad, the monad $\T$ is the free strict monoidal $V$-category monad on the category of $V$-categories.
\end{example}

\smallskip Further we fix a $T$-equipment $(\TT_0, \AA, T, \kappa^T, \nu^m, \nu^e)$, such that $\nu^m$ is invertible. 

Consider the multiplication downmap (see Section~\ref{monads}) of the monad $(\AA, \TT)$ in $\Eq$. It is a morphism of $\fM^\top(\Eq)$ 
\[(\TT_0^\d, \TT, \nu^m) : (\TT_0, \AA) \rightarrow (\one_{A_0}, \AA),\]

\noindent where $(\one_{A_0}, \AA) = \id(\AA)$ is the trivial $T$-equipment, and $\TT_0^\d : \TT_0 \rightarrow \one_{A_0}$ is the multiplication downmap of the monad $\TT_0$. Replacing $\nu^m$ in the triple by its inverse results in a morphism
\[(\TT_0^\d, \TT, (\nu^m)^{-1}) : (\TT_0, \AA) \rightarrow (\one_{A_0}, \AA)\]

\noindent of $\overline{\fM^\top(\Eq)}$. To this we can apply the Kleisli construction, by which we get a functor of equipments 
\[\LL = \kl(\TT_0^\d, \TT, (\nu^m)^{-1}) : \kl(\TT_0, \AA) \rightarrow \AA.\] 

\noindent In a more explicit way $\LL$ can be defined to be a functor of equipments which consists of the following components:

\begin{itemize}
\item The functor $L_0 = T_0 : A_0 \rightarrow A_0$.
\item The family of functors $L_{x, y}$ between the categories of vectors defined by
\[\xymatrix{A(x, Ty) \ar[r]^{F_{x, y}}& A(Tx, T^2y) \ar[r]^{A(Tx, m_y)} & A(Tx, Tx)} \]
\item The lax comparison maps $\kappa^L$ defined using $\kappa$ and $(\nu^m)^{-1}$.
\end{itemize}

\noindent Still more explicitly, for an object $x$, $L(x) = T(x)$, for a Kleisli vector $a : x \vectarr Ty$, $L(a)$ is the composite
\[\xymatrix{Tx \ar|{\object@{|}}[r]^{Ta} & T^2y \ar[r]^{m_y}& Ty.} \]

\noindent The components of $\kappa^L_0$ are:
\[\xymatrix{Tx \ar|{\object@{|}}[rr]_{Ti_x}^{\;}="2" \ar|{\object@{|}}@/^2pc/[rr]^{i_{Tx}}_{\;}="1" && Tx \ar[r]_{Te_x} \ar@/^2pc/[rr]^{1_{Tx}}& T^2x \ar[r]_{m_x} & Tx \ar@{=>} "1"; "1"+(0,-5) ^{\kappa^T_{x}}}\]

\noindent The components of $\kappa^L_2$ are:
\[ 
\xymatrix{
Tw  \ar|-{\object@{|}}@/_1.4pc/[rrrrd]_{T(T(a)b)}="1" \ar|-{\object@{|}}[rr]^{Tb} && T^2x \ar[rr]^{m_x} \ar|-{\object@{|}}[rrd]^{T^2a} \ar@{=>} {}+(0,-4.5); {}+(0,-9.5) ^{\kappa^T_{b,Ta}}&& Tx \ar|-{\object@{|}}[rr]^{Ta} \ar@{=>}{}+(0,-4);{}+(0,-9) ^{(\nu^m_a)^{-1}}&& T^2y \ar[rr]^{m_y} && Ty \\
&&&&T^3y \ar@<0.5ex>[rru]^>>>>>>>>>{m_{Ty}} \ar@<-1ex>[rru]_{Tm_{y}}&&\\
}
\]

\noindent The higher components of $\kappa^L$ are defined similarly, e.g. for $n = 3$ the components are
\[ 
\xymatrix{
Tv \ar|-{\object@{|}}@/_2.3pc/[rrrrrdd]_{T(T^2(a)T(b)c)}="1" \ar|-{\object@{|}}[r]^{Tc} & T^2w \ar|-{\object@{|}}[rrd]^{T^2b} \ar[rr]^{m_w} \ar@{=>}{}+(7,-9);{}+(7,-14) ^{\kappa^T_{c,Tb,T^2a}}&& Tw \ar|-{\object@{|}}[rr]^{Tb} \ar@{=>}{}+(0,-4);{}+(0,-9) ^{(\nu^m_b)^{-1}}&& T^2x \ar[rr]^{m_x} \ar|-{\object@{|}}[rrd]_{T^2a} \ar@{=>}{}+(0,-11);{}+(0,-16)^{(\nu^m_{Ta})^{-1}}&& Tx \ar|-{\object@{|}}[rr]^{Ta} \ar@{=>}{}+(0,-4);{}+(0,-9) ^{(\nu^m_a)^{-1}}&& T^2y \ar[r]^{m_y} & Ty\\
&&&T^3x \ar[rru]_{m_{Tx}} \ar|-{\object@{|}}[rrd]^{T^3a} &&&&T^3y \ar[rru]_{m_{Ty}}&& \\
&&&&& T^4y \ar[rru]_{m_{T^2y}} \ar@/_1.8pc/[rrrrruu]_{m_{4y}}&&&
}
\]

Applying the 2-functor $\mon$ to $\LL$, we obtain a functor between the categories of monoids:
\[\xymatrix{\L =  \mon(\LL) : \tmon(\AA, \TT) \ar[r]& \mon(\AA).}\]

\noindent In more details, on objects $\L$ is determined by the formula 

\[\L(x, a, \mu_a, \nu_a) = (Tx, m_xTa, \mu_{m_xTa}, \nu_{m_xTa}),\]

\noindent where $\mu_{m_xTa}$ and $\nu_{m_xTa}$ are defined by the diagrams
\[
\xymatrix{
Tx \ar[rr]|-{\object@{|}}^{Ta} \ar@/_1.1pc/[rrrrdd]|-{\object@{|}}_{Ta}="1" && T^2x \ar[rr]^{m_x}  \ar[rrd]|-{\object@{|}}^{T^2a} \ar@{=>} {}+(0,-7); {}+(0,-12) ^{T(\mu_a)\kappa_{-}}&& Tx \ar[rr]|-{\object@{|}}^{Ta}  \ar@{=>}{}+(0,-4); {}+(0,-9) ^{(\nu^m_a)^{-1}}&& T^2x \ar[rr]^{m_x}="3" && T^2x \\
&&&& T^3x \ar[d]^{Tm_x} \ar[rru]_{m_{Tx}} &&&&  \\
&&&& T^2x \ar@/_1.1pc/[uurrrr]_{m_x}&&&& 
}
\]
\noindent and
\[
\xymatrix{&& Tx  \ar[rrd]^{Te_x}  \ar@/^1pc/[rrrrd]^{1_{Tx}} \ar@{=>} {}+(0,-5); {}+(0,-10) ^{T(\eta_a)\kappa_{-}} &&&& \\
Tx  \ar[rrrr]|-{\object@{|}}_{Ta}^{\;}="3" \ar[rru]|-{\object@{|}}_{Ti_{x}}^{\;}="2" \ar@/^2.5pc/[rru]|-{\object@{|}}^{i_{Tx}}_{\;}="1" &&&& T^2x  \ar[rr]_{m_x}&& Tx. \ar@{=>} "1"; "2" ^{\kappa_x} 
}
\]

\noindent To a homomorphism of $T$-monoids $(f, \phi_f)$ $\L$ assigns a homomorphism of $T$-algebras $(Tf, \phi_{Tf})$ where $\phi_{Tf}$ is defined by 
\[
\xymatrix{
Tx \ar|-{\object@{|}}[rr]^{Ta}_{\;}="1" \ar[d]_{Tf}&& T^2x \ar[d]^{T^2f} \ar[rr]^{m_x} && Tx \ar[d]^{Tf}\\
Tx \ar|-{\object@{|}}[rr]_{Ta}^{\;}="2" && T^2y \ar[rr]_{m_y} && T^2y. \ar@{=>} "1"+(0,-3.5); "1"+(0,-8.5) ^{\phi_f}
}
\]

Consider now the multiplication upmap of the monad $(\AA, \TT)$. It is a morphism of $\fM^\bot(\Eq)$:
\[(\TT_0^\u, \TT, \nu^m) : (\one_{A_0}, \AA) \rightarrow (\TT, \AA).\]

\noindent By Lemma~\ref{lemma} the multiplication upmap and the multiplication downmap define a morphism of distributive pairs of monads. So the pair 
\beq{Td0}\big{(}(\TT_0^\u, \TT, \nu^m), (\TT_0^\d, \TT, \nu^m)\big{)}\eeq

\noindent determines a morphism of $\fM^\bot(\fM^\top(\Eq))$. We proceed relying on the obvious functorial nature of the constructions involved. Replacing $\nu^m$ in the second component of (\ref{Td0}) by its inverse, we get a pair
\[\big{(}(\TT_0^\u, \TT, \nu^m), (\TT_0^\d, \TT, (\nu^m)^{-1})\big{)}\]

\noindent which determines a morphism in $\fM^\bot(\overline{\fM^\top(\Eq)})$. Taking the Kleisli construction of the second component of this, we get a morphism of $\fM^\bot(\Eq)$
\beq{TdLnu}
(\TT_0^\u, \LL, \nu^{\prime m}) : (\one_{A_0}, \kl(\TT_0, \AA)) \rightarrow (\TT_0, \AA).
\eeq

\noindent where $(\one_{A_0}, \kl(\TT_0, \AA)) = \id\kl(\TT_0, \AA)$ is the trivial $T$-equipment. The components of the modification $\nu^{\prime m}$ can be verified to be the 2-cells of $A$:
\begin{equation}\label{numdash}
\xymatrix{
T^2x \ar|-{\object@{|}}[rr]^{T^2a}_{\;}="1" \ar[d]_{m_{x}}&& T^3y \ar[d]^{m_{Ty}} \ar[rr]^{Tm_y} && T^2y \ar[d]^{m_y}\\
Tx \ar|-{\object@{|}}[rr]_{Ta}^{\;}="2" && T^2y \ar[rr]_{m_y} && T^2y. \ar@{=>} "1"+(0,-3.5); "1"+(0,-8.5) ^{\nu^m_a}
}
\end{equation}

\smallskip
\begin{definition} The \textbf{free $T$-algebra functor}
\[\xymatrix{\M : \tmon(\AA) \ar[r]& \talg(\AA)}\] 

\noindent is defined as the composite 
\[
\xymatrix{\Eq(\II, \kl(\TT_0, \AA)) \ar[d]^{\id^\bot_{\Eq}} \\
\fM^\bot(\Eq)(\id(\II), \id\kl(\TT_0, \AA)) \ar[d]^{[\id(\II), (\TT_0^\u, \LL, {\nu^{\prime m}})]}\\
\Eq\big{(}(\one_{I_0}, \II), (\TT_0, \AA)\big{)}.}
\]
\end{definition}

Observe that, (\ref{TdLnu}) exhibits $\LL$ as a module of the monad $(\AA, \TT)$ in $\Eq$
\[\xymatrix{
&\kl(\TT_0, \AA) \ar[dl]_{\LL}="1" \ar[dr]^{\LL}="2"&\\
\AA  \ar[rr]_{\TT} && \AA. \ar@{=>} "1"+(8,-3); "1"+(13,-3)
}\]

\noindent Applying the 2-functor $\mon$ we get a module of the monad $(\mon(\AA), \T)$ in $\cat$
\[\xymatrix{
&\tmon(\AA) \ar[dl]_{\L}="1" \ar[dr]^{\L}="2"&\\
\mon(\AA)  \ar[rr]_{\T} && \mon(\AA).  \ar@{=>} "1"+(10,-3); "1"+(15,-3)
}\]

\noindent Then, alternatively, the free $T$-algebra functor $\M$ is the comparison functor 
\[\xymatrix{\tmon(\TT_0, \AA) \ar[r]&  \mon(\AA)^{\T} = \talg(\TT_0, \AA).}\] 

\noindent In more details, for a $T$-monoid $(x, a)$, the $T$-algebra $\M(x, a)$ consists of the monoid $\L(x, a)$ and the algebra structure $(m, \sigma_{m_xTa}) : \T\L(x, a) \rightarrow \L(x, a)$, where $\sigma_{m_xTa}$ is defined by (\ref{numdash}). So, on objects we have a formula for $\M$

\[\M(x, a, \mu_a, \nu_a) = (Tx, m_xTa, \mu_{m_xTa}, \nu_{m_xTa}, m, \sigma_{m_xTa}).\]

\begin{example} In the case $\AA = \MMatV$, the free $T$-algebra functor is the functor from the category of $(T, V)$-categories to the category of $T$-algebras constructed in \cite{ChClHo14}. In particular, when $\TT_0$ is taken to be the free monoid monad, then the free $T$-algebra functor becomes the free monoidal $V$-category functor on the category of $V$-multicategories. Furthermore, when $V = \set$ we obtain the free monoidal category functor from the category of multicategories to the category of strict monoidal categories. 
\end{example}

\smallskip Let us look at the details of the general constructions developed until now in the basic multicategory-monoidal category case. This case is captured by the context of the $T$-equipment $(\TT_0, \AA)$, where $\AA = (\set, \mat(\set))$ and $\TT_0$ is the free monoid monad on $\set$.

The bicategory $\mat(\set)$ is the same as the bicategory of spans $\span$. Its objects are sets. A morphism in it $a: x \vectarr y$ is a diagram of the form

\beq{span}
\xymatrix{&a_o \ar[rd]^{a^\r} \ar[ld]_{a^\l}&\\ x&&y.
}
\eeq

\noindent A 2-cell $a \Rightarrow b$ is a map $a_o \rightarrow b_o$ which respects the two span legs. The composite of spans is defined by taking a pullback, as shown here:

\[\xymatrix{
&&a_o\times_{y} b_o\ar[rd]^{a^\r} \ar[ld]_{a^\l}&&\\
&a_o \ar[rd]^{a^\r} \ar[ld]_{a^\l}&&b_o\ar[rd]^{b^\r} \ar[ld]_{b^\l}&\\ 
x&&y&&z.
}\]

\noindent The identities are those spans whose both legs are identities. The equipment $\SSpan = (\set, \span)$ arises from the pseudofunctor $\set \rightarrow \span$ which takes a set map $f : x \rightarrow y$ to the span with the left leg the identity and the right leg $f$.

Let now $T_0$ be an endofunctor on the category of sets.  Define functors $T_{x, y} : \span(x, y) \rightarrow \span(Tx, Ty)$ by setting for any span $a$, $T(a)$ to be the span
\[\xymatrix{&Ta_o \ar[rd]^{Ta^\r} \ar[ld]_{Ta^\l}&\\
Tx&&Ty.
}\]

\noindent Define 2-cells $\kappa_{a, b} : T(ab) \rightarrow T(a)T(b)$ in $\span$ by the obvious maps

\beq{pb}
T(a_o\times_{y} b_o) \rightarrow T(a_o)\times_{Ty}T(b_o).
\eeq

\noindent For $n > 2$, define 2-cells $\kappa_{a_1, a_0, \cdots a_n} : T(a_1, a_2, \cdots a_n) \rightarrow T(a_1)T(a_2)\cdots T(a_n)$ analogously. This gives a colax functor $\TT : \SSpan \rightarrow \SSpan$. If $T$ preserves pullbacks then (\ref{pb}) are invertible. 

Suppose next that $(T_0, m, e)$ is a monad on $\set$, such that $T_0$ preserves pullbacks. This monad has a lift to a monad on the equipment $(\set, \SSpan)$. It consists of the (lax) functor $\TT : \SSpan \rightarrow \SSpan$ whose lax comparison structure is given by the inverse of (\ref{pb}) and the monad multiplication and the monad unit which are the lax transformations $\mm : \TT^2 \rightarrow \TT$ and $\ee : 1_{\SSpan} \rightarrow \TT$ with the 2-cells (\ref{em}) needed to define them given by the obvious maps:

\beq{nuspan}
T^2a_o \rightarrow T^2x\times_{Tx}Ta^1,\qquad\qquad
a_o \rightarrow x\times_{Tx}Ta^1.
\eeq
 
\noindent  Hence we have a $T$-equipment $(\TT_0, \SSpan)$. Note that if the naturality squares of the natural transformations $e$ and $m$ are pullback squares, then (\ref{nuspan}) are invertible.

Further in this section we take $\TT_0$ to be the free monoid monad. Thus, for any set $x$, $Tx$ is the set of lists of elements of $x$. A component of the monad multiplication $m_x : T^2x \rightarrow Tx$ concatenates lists of lists into lists. While, a component of the monad unit $e_x : x \rightarrow Tx$ sends an element to the singleton list. It is a well know fact that the free monoid monad is Cartesian, i.e. its underlying endofunctor preserves pullbacks, and the naturality squares of its multiplication and unit are pullback squares. Therefore, we get a $T$-equipment $(\TT_0, \SSpan)$, for which $\kappa^T$, $\nu^m$ and $\nu^e$ are invertible.

A monoid in $\SSpan$ is the same as a monad in the bicategory $\span$, which are well known to be ordinary categories. A $T$-monoid in $(\TT_0, \SSpan)$ is a multicategory (otherwise known as a planar colored operad) \cite{Her00}. In more details, a $T$-monoid $(x, a)$ is a  multicategory whose set of objects (colors) is $x$, and whose set of multimorphisms is $a_o$. The span $a$ itself is a diagram

\[\xymatrix{&a_o \ar[rd]^{a^\r} \ar[ld]_{a^\l}&\\ x&&Tx.
}\]

\noindent in which $a^\r$ takes a multimorphism to its target, which is an element of $x$, while $a^r$ takes a multimorphism to its source, which is a list of elements of $x$, i.e. an element of $Tx$. The $T$-monoid multiplication and unit are 2-cells of $\span$ which amount to maps

\[a_o \times_{Tx} Ta_o \rightarrow a_o,\qquad\qquad x  \rightarrow a_o\]

\noindent respectively. These correspond to the operations of multimorphism composition and taking identity multimorphisms.

Recall that $\T = \mon(\TT) : \mon(\SSpan) \rightarrow \mon(\SSpan)$ takes a monoid $(x, a)$ to a monoid $(Tx, Ta)$. The latter is easily identified to be the free strict monoidal category on the category $(x, a)$. In fact it is easily seen that $\T$ is the free strict monoidal category monad. Correspondingly, the $T$-algebras in the current context are the strict monoidal categories.

Let us compute the free $T$-algebra $\M(x, a)$, showing it to be the free strict monoidal category on the multicategory $(x, a)$ (see \cite{Her00}). First, let us identify its underlying monoid $\L(x, a)$. Recall that this is given by the data $(Tx, m_xTa, \mu_{m_xTa}, \nu_{m_xTa})$. Thus, as a category, its set of objects is $Tx$, i.e. the set of lists of objects of the original multicategory. Its set of morphisms is determined by the span $m_xTa$ which is 

\[\xymatrix{
&Ta_o \ar[rd]^{Ta^\r} \ar[ld]_{Ta^\l}&&\\
Tx&&T^2x \ar[rd]^{m_x}&\\
&&&Tx.
}\]

\noindent We can see that, a morphism of $\L(x, a)$ is a list of multimorphisms of the original multicategory. The target of such a list is the list consisting of the targets of the multimorphisms in the list. The source is the concatenation of the sources (themselves lists) of the multimorphisms in the list. This agrees precisely with the construction of the free strict monoidal category on a multicategory. That the multiplication and the unit of  $\L(x, a)$ coincide with the composition and the identities of the strict free monoidal category is left to the reader.

Finally let us work out the details of the strict monoidal structure on $\L(x, a)$. Recall, that the latter can be given as the algebra $(m, \sigma_{m_xTa}) : \T\L(x, a) \rightarrow \L(x, a)$. This is the functor whose object part is $m_x : T^2x \rightarrow Tx$, i.e. concatenation of lists of lists of objects to lists of objects. Its morphism part, given by (\ref{numdash}), amounts to concatenation of lists of lists of multimorphisms to lists of multimorphisms. We have arrived precisely at the strict monoidal structure on the free strict monoidal category.

Another basic example of operads captured as $T$-monoids are symmetric multicategories. For this see Section~\ref{further}.

\section{Star equipments}\label{starequipments}

Informally, a $\ast$-equipment is an equipment for which, in addition to the scalar actions, scalar opactions are defined, and the actions and the opactions are adjoint to each other. This means that, for scalars $f : x \rightarrow w$ and $g : z \rightarrow y$, we can draw diagrams 
\[ 
\xymatrix{w \ar[r]^{f^\ast} & x \ar[r]|-@{|}^{a}& y} \qquad\qquad \xymatrix{x \ar[r]|-@{|}^{a}& y \ar[r]^{g^\ast}& z}
\]

\noindent which evaluate to vectors $af^\ast : w \vectarr y$ and $g^\ast a : x \vectarr z$ respectively, and there are universal 2-cells
\[
\xymatrix{x \ar|-{\object@{|}}[rr]^a_{\quad}="1" \ar[d]_{f} && y \\
w \ar[rr]_{f^\ast}^{\quad}="2" &  \ar@{=>} "1"+(0,-3.5);"1"+(0,-8.5)& x \ar|-{\object@{|}}[u]_a} 
\qquad\qquad\qquad\qquad
\xymatrix{x \ar|-{\object@{|}}[rr]^a_{\quad}="1"  \ar|-{\object@{|}}[d]_{a} && y \\
y \ar[rr]_{g^\ast}^{\quad}="2"&\ar@{=>} "1"+(0,-3.5);"1"+(0,-8.5)& z \ar[u]_{g}.}
\]

\noindent We turn to formal definitions. First we introduce a notion of a $\ast$-module, which appears in \cite{Ca} under the name of ``starred module''.

\begin{definition} A $\ast$-module  $A$ from a category $X$ to a category $Y$ is a module together with a choice, for each morphism $f : x \rightarrow y$ and an object $w$ of $X$, of a right adjoint $A(w, f^\ast)$ of the functor $A(w, f)$, and of a left adjoint $A(f^\ast, w)$ of the functor $A(f, w)$
\[\xymatrix{A(w, x) \ar@/_1pc/[rr]_{A(w, f)}\ar@{}[rr]|-\top&&A(w, y) \ar@/_1pc/[ll]_{A(w, f^\ast)}}, \qquad\qquad \xymatrix{A(y, w) \ar@/_1pc/[rr]_{A(f, w)}\ar@{}[rr]|-\bot&&A(x, w)\ar@/_1pc/[ll]_{A(f^\ast, w)}}\]

\noindent such that the natural transformations
\[
\xymatrix{
A(x, y) \ar[rr]^{A(f^\s, y)}_{\;}="1" \ar[d]_{A(x, k)}&& A(w, y) \ar[d]^{A(w, k)}\\ 
A(x, z)  \ar[rr]_{A(f^\s, z)}^{\;}="2" && A(w, z)  \ar@{=>} "1"+(0,-8.5);"1"+(0,-3.5)
}
\qquad\qquad
\xymatrix{
A(x, y) \ar[rr]^{A(x, g^\s)}_{\;}="1" \ar[d]_{A(l, y)} && A(x, z) \ar[d]^{A(l, x)}\\ 
A(w, y)  \ar[rr]_{A(w, g^\s)}^{\;}="2" && A(w, z) \ar@{=>} "1"+(0,-3.5);"1"+(0,-8.5)
}
\]

\noindent defined as mates of the structural isomorphisms 
\[A(x, k)A(f, y) \cong A(f, z)A(w, k), \qquad A(w, g)A(l, x) \cong A(l, y)A(x, g),\] 

\noindent are invertible, thus giving isomorphisms

\begin{equation}\label{opactiso1}
A(f^\s, z)A(x, k) \cong A(w, k)A(f^\s, y), \qquad\qquad  A(l, x)A(x, g^\s)\cong A(w, g^\s)A(l, y).
\end{equation}

\noindent Furthermore the natural transformations
\[
\xymatrix{
A(x, z) \ar[rr]^{A(f^\s, z)}_{\;}="1" && A(w, z) \\ 
A(x, y)  \ar[rr]_{A(f^\s, y)}^{\;}="2" \ar[u]^{A(x, g^\s)}&& A(w, y) \ar[u]_{A(w, g^\s)} \ar@{=>} "1"+(0,-3.5);"1"+(0,-8.5)
}
\qquad\qquad
\xymatrix{
A(w, y) \ar[rr]^{A(x, g^\s)}_{\;}="1"&& A(x, z)\\ 
A(x, y)  \ar[rr]_{A(x, g^\s)}^{\;}="2"  \ar[u]^{A(f^\s, y)} && A(x, z)  \ar[u]_{A(f^\s, z)} \ar@{=>} "1"+(0,-8.5);"1"+(0,-3.5)
}
\]

\noindent obtained as the mates of isomorphisms (\ref{opactiso1}) (with $k = g$ and $l = f$) should be equal to each other and be invertible, thus giving an isomorphism

\begin{equation}\label{opactiso2}
A(f^\s, z)A(x, g^\s) \cong A(f^\s, y)A(w, g^\s).
\end{equation}

\end{definition}

\noindent Let us shortly write  $af^\s$ for $A(f^\s, a)$, and $f^\s a$ for $A(a, f^\s)$. Since composites of adjoints are adjoints, there are isomorphism  $(fg)^\s a \cong f^\s (g^\s a)$ and $(fg)^\s a \cong f^\s (g^\s a)$. The isomorphisms (\ref{opactiso1}) and (\ref{opactiso2}) say that $(fa)g^\s \cong f(ag^\s)$, $f^\s(a g) \cong (f^\s ag)$ and $(f^\s a)g^\s \cong f^\s(ag^\s)$. Via all these, $A$ can be regarded as a module in three other ways: as a module from $X^\op$ to $Y$, as a module from $X$ to $Y^\op$, and as a module from $X^\op$ to $Y^\op$. 

\smallskip For any category $X$ there exists a 2-category $\Pi(X)$ which freely adjoins right adjoints to all morphisms of $X$.  More precisely, there is a functor $X \rightarrow \Pi(X)$ with the codomain a 2-category, which has the property that the image of any morphism under it has a right adjoint, and it is universal among functors with this property. $\Pi(X)$ exists for formal reasons. An explicit description, which we now recount briefly, can be found in \cite{DPP}. Objects of $\Pi(X)$ are the same as objects of $X$. Its morphisms can be represented by chains of arrows labeled by morphisms of $X$, such as 
\[\xymatrix{\ar[r]^{f_1} & \ar[r];[]_{g_1}& \ar[r]^{f_2} & \ar@{}[r]^{\hdots} & \ar[r];[]_{g_{n-1}} &\ar[r]^{f_n} &}.\]

\noindent The inclusion $X \rightarrow \Pi(X)$ is identical on objects, and takes an arrow $f$ to the chain of a single morphism directed to the left and labeled with $f$. Its right adjoint is given by the chain of a single arrow directed to the right and labeled with $f$. It was shown in \cite{DPP}, that $X \rightarrow \Pi(X)$ is faithful and locally fully faithful. 

It is not difficult to see that, a $\ast$-module $A$ from $X$ to $Y$ is essentially the same as a pseudofunctor
\beq{APiPi}
\xymatrix{A : \Pi(X)^\op \times \Pi(Y) \ar[r]& \cat.}
\eeq

\smallskip Now we come to a formal definition of $\ast$-equipments as well as their functors and transformations.

\begin{definition}
A \textbf{$\ast$-equipment} $\AA = (A_0, A, P, \xi)$ is an equipment where $A$ has a $\ast$-module structure, such that the family of $n$-fold composition functors (\ref{Pxxx}) is pseudonatural in the outermost arguments and pseudo-dinatural in inner arguments when $A$ is considered as a module from $A_0^\op$ to itself via the left and the right opactions. 

A functor of $\ast$-equipments $\AA \rightarrow \BB$ is a functor between the underlying equipments, such that  (\ref{Fxy}) becomes a pseudonatural family when $A$ is considered as a module from $A_0^\op$ to itself via the left and the right opactions. 

A transformation of functors between equipments is by definition a transformation between the underlying equipment functors. 

$\ast$-equipments, functors of $\ast$-equipment and transformations form a 2-category, which we denote by $\sEq$. 
\end{definition}

\begin{example} Our $\ast$-equipments are closely related to the equipments of \cite{CrSh10}. Like the compositions, the opactions in \cite{CrSh10} are defined by universal properties rather than an extra structure, which it is in our setting.
\end{example} 

\begin{example} Suppose that $\AA$ is an equipment which comes from a pseudofunctor $J : A_0 \rightarrow A$ from a category to a bicategory. Then $\AA$ becomes a $\ast$-equipment as soon as for each morphism $f$ of $A_0$, $J(f)$ has a right adjoint. So, the proarrow equipments (\cite{Woo82}) can be regarded as $\ast$-equipments.
\end{example}

\begin{example}
The pseudofunctor $\set \rightarrow \mat(V)$ considered earlier is a proarrow equipment. An image of a set map $f$ has a right adjoint given by the matrix which is the monoidal unit at pairs of the form $(f(x), x)$, and the initial object otherwise. Hence $\MMatV$ becomes a $\ast$-equipment. 
\end{example}

\smallskip A $T$-$\ast$-equipment is a monad in $\sEq$. Or in line with Definition~\ref{defteq}: 

\begin{definition} A \textbf{$T$-$\ast$-equipment} is a $T$-equipment $(\TT_0,  \AA, \kappa, \nu)$ with a $\ast$-equipment structure on $\AA$, such that $\TT : \AA \rightarrow \AA$ is a functor of $\ast$-equipments. 
\end{definition}

We have 2-categories of monads $\fM^\top(\sEq)$ and $\fM^{\bot}(\sEq)$, as well as 2-categories $\overline{\fM^\top(\sEq)}$ and $\overline{\fM^\bot(\sEq)}$, defined analogously to their non-star versions. All of these have $T$-$\ast$-equipments as their objects. Their morphisms are respectively morphisms of $\fM^\top(\Eq)$, $\fM^\bot(\Eq)$, $\overline{\fM^\bot(\Eq)}$ and $\overline{\fM^\bot(\Eq)}$ whose underlying equipment functors are $\ast$-equipment functors.

The Kleisli construction on a $T$-$\ast$-equipment is by definition the Kleisli construction on the underlying equipment. The result however is only an equipment, because the obvious $\ast$-structure of the module of Kleisli vectors $A(-, T-)$ inherited from the $\ast$-structure of $A$ is not compatible with the Kleisli composition. So the Kleisli construction 2-functor on $T$-$\ast$-equipments lands in the category of equipments

\[\xymatrix{\kl : \overline{\fM^\top(\sEq)}  \ar[r] & \Eq.}\]

Suppose that $A$ is a $\ast$-module from $Z$ to $X$. Suppose that $t : F \rightarrow G : Y \rightarrow X$ is a natural transformation. We will say that the opaction of $t$ on $A$ is \textbf{Cartesian} if the natural transformation
\[
\xymatrix{
A(z, Gy) \ar[rr]^{A(Gx, t_y^\ast)}_{\;}="1" \ar[d]_{A(z, Gg)} && A(Gx, Fy) \ar[d]^{A(Gx, Fg)}\\ 
A(z, Gv)  \ar[rr]_{A(Gx, t_v^\ast)}^{\;}="2" && A(z, Fv) \ar@{=>}  "1"+(0,-3.5); "1"+(0,-8.5)
}
\]

\noindent defined as a mate of the structural isomorphism $A(Gx, Gg)A(Gx, t_y) \cong A(Gx, t_v)A(Gx, Fg)$ is invertible. This means that there are invertible 2-cells $Fgt_y^\ast a \cong t_v^\ast Gga$. A transformation $\tt = (t, \nu^t) : \FF \rightarrow \GG : \BB \rightarrow \AA$ between functors of equipments will be said to be Cartesian if the opaction of $t$ on $A$ is Cartesian.

Let $\sEq_{\mathrm{C}}$ denote the sub 2-category of $\sEq$ with 2-cells restricted to the Cartesian transformations. A monad in $\sEq_{\mathrm{C}}$ is a $T$-$\ast$-equipment $(\TT_0, \AA)$ for which $\mm : \TT^2 \rightarrow \TT$ and $\ee : 1_{\AA} \rightarrow \TT$ are Cartesian transformations. A Kleisli equipment of such a $T$-equipment inherits a $\ast$-equipment structure. Moreover, we have a 2-functor 

\[\xymatrix{\kl : \overline{\fM^\top(\sEq_{\mathrm{C}})}  \ar[r] & \sEq_{\mathrm{C}}.}\]

Changing the discussion from the downmaps to the upmaps, let $\fM_\mathrm{C}^\bot(\sEq)$ be the sub 2-category of $\fM^\bot(\sEq)$ whose objects are all $T$-$\ast$-equipments, but whose morphisms are restricted to those $(\FF_0, \FF, \nu^u) : (\SS_0, \BB) \rightarrow (\TT_0, \AA)$ for which the transformation $(u, \nu^u) : \FF\SS \rightarrow \TT\FF$ is Cartesian. The opaction structure allows a functorial extension of the Kleisli construction to $\fM_\mathrm{C}^\bot(\sEq)$. Given a morphism $(\FF, \FF_0, u): (\SS_0, \BB) \rightarrow (\TT_0, \BB)$ in $\fM_\mathrm{C}^\bot(\sEq)$, a functor between the Kleisli equipments $\kl^\ast(\FF_0, \FF) : \kl^\ast(\SS_0, \BB) \rightarrow \kl^\ast(\TT_0, \AA)$ is defined to have the following components:

\begin{itemize}
\item The scalar functor $F_0$.

\item The family of functors between the categories of vectors
\[
\xymatrix{
B(x, Sy) \ar[r]^{F_{x, y}}& A(Fx, FSy) \ar[r]^{A(Fx, u_y^\ast)}& A(Fx, TFy),
}
\]

\noindent which is pseudonatural by the virtue of the Cartesian property of $(u, \nu^u)$.  

\item The lax comparison structure defined from $\kappa^F$, a certain transform of an equality satisfied by the upmap $(F_0, u)$, and a certain transform of $\nu^u$, as described below.
\end{itemize}

\noindent $\kl^\ast$ extends in the obvious way to the 2-cells, and so gives a 2-functor:

\[\xymatrix{\kl^\ast : \fM_\mathrm{C}^\bot(\sEq) \ar[r]& \Eq.}\]

\noindent To give a more concrete description, we invoke the fact that the $\ast$-module $A$ extends to a pseudofunctor $\Pi(A_0)\times \Pi(A_0) \rightarrow \cat$. This means that the 2-cells of $A$ are acted upon by morphisms and 2-cells of $\Pi(A_0)$. We use this to form pasting diagrams involving 2-cells of $A$ and those of $\Pi(X)$. We are also enabled to talk about mates of 2-cells of $A$ under adjunctions of morphisms of $\Pi(X)$. So: The $0$ component of the lax comparison structure of $\kl^\ast(\FF_0, \FF)$ is defined by
\[\xymatrix{Fx \ar|{\object@{|}}[rr]_{Fi_x}^{\;} \ar|{\object@{|}}@/^2pc/[rr]^{i_{Fx}}_{\;}="1" && Fx \ar[r]_{Fe_x} \ar@/^2pc/[rr]^{e_{Fx}}_{\;}="2"& FTx \ar[r]_{u_x^\ast} & SFx \ar@{=>} "1"+(0,-1); "1"+(0,-5) ^{\kappa^F_x} \ar@{=>} "2"+(0,-1); "2"+(0,-4) ^{\tilde1_{Fe_x}} }\]

\noindent where $\tilde1_{Fe_x}$ is the mate of the identity $Fe_x = ue_{Fx}$. The $2$ component is defined by
\[ 
\xymatrix{
Fw  \ar|-{\object@{|}}@/_1.3pc/[rrrd]_{F(T(a)b)}="1" \ar|-{\object@{|}}[r]^{Fb} & FTx \ar[rr]^{u_x^\ast} \ar|-{\object@{|}}[rrd]^{FTa} \ar@{=>} "1"+(-2,11) ; "1"+(-2,6) ^{\kappa^F_{a,Ta}}&& SFx \ar|-{\object@{|}}[rr]^{SFa} \ar@{=>} {}+(0,-4);{}+(0,-9) ^{\tilde\nu_a^u}&& SFTy \ar@{=>} {}+(0,-4.5);{}+(0,-9.5) ^{\tilde1_{u_ym_{Ty}}} \ar[r]^{Su_y^\ast} & S^2Ty \ar[r]^{m_{Ty}} & STy\\
&&&FT^2y \ar[rru]^{u^\ast_{Ty}} \ar[rr]_{Fm_{y}}&&FTy \ar[rru]_{u_y^\ast}\\
}
\]

\noindent where $\tilde\nu_a^u$ is the mate of $\nu_a^u$, and $\tilde1_{u_ym_{Ty}}$ f the identity $u_ym_{Ty} = Fm_yu_{Ty}Su_y$. The higher components of the lax comparison structure are defined similarly. We leave it to the reader to define $\kl^\ast$ on the 2-cells of $\fM_{\mathrm{C}}^\bot(\sEq)$, completing the construction of the 2-functor. A more conceptual insight on the definition of $\kl^\ast$ will be given at the end of Section~\ref{laxmonads}. Observe that, if we restrict the domain of $\kl^\ast$ to the Cartesian $T$-$\ast$-equipments, then we get a 2-functor landing in $\ast$-equipments

\[\xymatrix{\kl^\ast : \fM^\bot(\sEq_{\mathrm{C}})  \ar[r] & \sEq_{\mathrm{C}}.}\]

\smallskip
\begin{definition} The \textbf{underlying $T$-monoid functor}
\[\xymatrix{\K : \talg(\TT_0, \AA) \ar[r]& \tmon(\TT_0, \AA)}\]

\noindent is defined by the $\ast$-Kleisli construction 
\[\xymatrix{\kl^\ast : \fM^{\bot}(\sEq)\big{(}(\one_{I_0}, \II), (\TT_0, \AA)\big{)} \ar[r]& \Eq\big{(}\II, \kl(\TT_0, \AA)\big{)}}.\]

\noindent ($\kl^\ast$ is defined on the left hand side since, trivially, every transformation between functors $\II \rightarrow \AA$ is Cartesian.) 
\end{definition}

\noindent Here is a more explicit description of $\K$. Suppose that $(x, b, \mu_b, \nu_b, h, \sigma_h)$ is a data for a $T$-algebra. Then $\K$ gives a $T$-monoid with the formula

\[\K(x, b, \mu_b, \nu_b, h, \sigma_h) = (x, h^\ast b, \mu_{h^\ast b}, \eta_{h^\ast b}),\] 

\noindent where $\mu_{h^\ast b}$ is defined by
\[
\xymatrix{
x \ar[rr]|-{\object@{|}}^b \ar@/_1.1pc/[rrrrd]|-{\object@{|}}_b^{\;}="1" && x \ar[rr]^{h^\s}  \ar[rrd]|-{\object@{|}}^b \ar@{=>} {}+(0,-3.5); {}+(0,-8.5) ^{\mu_b} && Tx \ar[rr]|-{\object@{|}}^{Tb} \ar@{=>} {}+(0,-3.5); {}+(0,-8.5) ^{\tilde\sigma^h}&& Tx \ar[rr]^{Th^\s}_{\;}="2"  && T^2x \ar[d]^{m_x}\\
&&&& x \ar[rru]^{h^\s} \ar[rrrr]_{h^\s}^{\;}="3" &&&& Tx \ar@{=>} "2"+(-5,-3.5); "2"+(-5,-8.5)^{\tilde1_{hm}}
}
\]

\noindent where $\tilde\sigma^h$ is the mate of $\sigma^h$, and $\tilde1_{hm}$ is the mate of the identity $hTh = hm$, and $\eta_{h^\ast b}$ is defined by
\[
\xymatrix{x  \ar[rr]|-{\object@{|}}_b^{\;}="1" \ar@/^2.pc/[rr]|-{\object@{|}}^{i_x}_{\;}="3" && x  \ar[rr]_{h^\s}^{\;}="2" \ar@/^2.pc/[rr]^{e_x}_{\;}="4"&& Tx \ar@{=>} "3"; "3"+(0,-5) ^{\eta_b} \ar@{=>} "4"; "4"+(0,-5) ^{\tilde1_{1_x}}}
\]

\noindent where $\tilde1_{1_x}$ is the mate of the identity $he_x = 1_x$.  To a homomorphism of $T$-algebras $(f, \phi_f) : (y, b, h) \rightarrow (y', b', h')$ $\K$ assigns a homomorphism of $T$-monoids $(f, \tilde\phi_f)$, where $\tilde\phi_f$ is defined by 
\[
\xymatrix{
y \ar|-{\object@{|}}[rr]^{b}_{\;}="1" \ar[d]_{f}&& x \ar[d]^{f} \ar[rr]^{h^\s}_{\;}="2" && Ty \ar[d]^{Tf}\\
y' \ar|-{\object@{|}}[rr]_{b'}^{\;}="3" && y' \ar[rr]_{h^{\prime \s}}^{\;}="4" && Ty'
 \ar@{=>} "1"+(0,-3.5); "1"+(0,-8.5) ^{\phi_f} \ar@{=>} "2"+(0,-3.5); "2"+(0,-8.5)^{\tilde1_{fh}} }
\]

\noindent wherein $\tilde1_{fh}$ is the mate of the identity $Tfh = hg$.

\begin{example} In the case $\AA = \MMatV$, the underlying $T$-monoid functor is the functor from the category of $T$-algebras to the category of $(T, V)$-categories constructed in \cite{ChClHo14}. In particular, when $\TT_0$ is taken to be the free monoid monad, then the underlying $T$-monoid functor becomes the underlying $V$-multicategory functor on the category of $V$-categories.
\end{example}

\smallskip The following is a generalization of Theorem 5.4 of \cite{ChClHo14}:

\begin{theorem}\label{adjthm}
Given a $T$-$\ast$-equipment $(\TT_0, \AA, \kappa^T, \nu^m, \nu^e)$, such that $\nu^m$ is invertible, the  free $T$-algebra functor $\M$ is a left adjoint to the underlying $T$-monoid functor $\K$.
\[\xymatrix{\talg(\TT_0, \AA)  \ar@/_1em/[rr]_{\K}\ar@{}[rr]|-\top&&\tmon(\TT_0, \AA)  \ar@/_1em/[ll]_{\M}}.\]
\end{theorem}

\noindent The component of the unit of this adjunction at a $T$-monoid $(x, a)$ is the homomorphism of $T$-monoids $(e_x, \phi_{e_x}): (x, a) \rightarrow \K\M(x, a)$, where $\phi_{e_x}$ is defined by
\[
\xymatrix{
x  \ar[rrrrrr]|-{\object@{|}}^a_{\;}="1"  \ar[d]_{e_x}&&&&&& Tx  \ar[d]^{Te_x} \ar[dllll]_{e_{Tx}}\\
Tx  \ar[rr]|-{\object@{|}}_a^{\;}="2"&& T^2x \ar[rr]_{m_x}&& Tx \ar[rr]_{m^\s_x}^{\;}="3"&& T^2x \ar@{=>} "2"+(0,8.5); "2"+(0,3.5) ^{e_a} \ar@{=>} "3"+(0,6); "3"+(0,1) ^{\tilde1_{me_{Tx}}}
}
\]

\noindent where $\tilde1_{me_{Tx}}$ is the mate of the identity map $m_xTe_x = m_xe_{Tx}$. The component of the counit at a $T$-algebra $((y, b), (h, \sigma_h))$ is the $T$-algebra homomorphism $(h, \phi_h) : \M\K((x, b), (h, \sigma_h)) \rightarrow ((x, b), (h, \sigma_h))$, where $\phi_h$ is defined by
\[
\xymatrix{
Ty \ar[rr]|-{\object@{|}}^{Tb}_{\;}="1" \ar[d]_{h} && Ty \ar[rr]^{Th^\s} \ar[drrrr]_{h}&& T^2y \ar[rr]^{m_y}_{\;}="2"&& Ty  \ar[d]^{Th}   \\
y  \ar[rrrrrr]|-{\object@{|}}_b^{\;}="3"  &&&&&& Ty \ar@{=>} "1"+(0,-3.5); "1"+(0,-8.5) ^{\sigma_f} \ar@{=>} "2"+(0,-1.5); "2"+(0,-6.5) ^{\tilde1_{hTh}}
}
\]

\noindent where $\tilde1_{hTh}$ is the mate of the identity $Thm_x = hTh$.

Below we state a more general result. Varying the $T$-equipment $(\TT_0, \AA)$, (\ref{TdLnu}) becomes a family of morphisms of $\fM^\bot(\Eq)$. This family of morphisms can be given a structure of a lax natural transformation $\id^\bot_{\Eq}\kl^\ast \rightarrow 1_{\fM_\mathrm{C}(^\ast\Eq)}$. It can be shown that:

\begin{theorem}
The $\ast$-Kleisli 2-functor $\kl^\ast : \fM_\mathrm{C}^\bot(\sEq) \rightarrow \Eq$ is a lax right 2-adjoint to the 2-functor $\id^\bot_{\Eq} :  \Eq \rightarrow \fM_\mathrm{C}^\bot(^\ast\Eq)$, with the counit the lax natural family (\ref{TdLnu}), and the unit the trivial lax natural transformation. In particular, for any $T$-$\ast$-equipment $(\TT_0, \AA)$ with an invertible $\nu^m$ and any $\ast-$equipment $\BB$, there is an adjudication
\beq{ladj}
\xymatrix{\sEq_\mathrm{C}\big{(}(\one_{B_0}, \BB), (\TT_0, \AA)\big{)} \ar@/_1.5em/[r]\ar@{}[r]|-\top& \fM_\mathrm{C}^\bot(\sEq)\big{(}\BB, \kl^\ast(\TT_0, \AA)\big{)}. \ar@/_1.5em/[l]}
\eeq

\end{theorem}

\noindent A lax 2-adjunction between 2-categories (\cite{Gr74}), is like an adjunction except that its unit and counit are lax natural transformations and the triangle identities are replaced by non-invertible 2-cells themselves satisfying triangle-type identities. Instead of the usual isomorphisms between the homsets of an adjunction, a lax 2-adjunction gives rise to a family of adjunction between homcategories as (\ref{ladj}). Theorem~\ref{adjthm} is a special case of (\ref{ladj}) with $\BB$ taken to be the terminal equipment $\II$.

\section{Pseudomonads and pseudomodules}\label{psmonpsmod}
The purpose of this section is to provide a background for the subsequent section where we define the tricategory of modules. The theory outlined here also provides a generalization of the formal theory of monads from the context of the bicategory to the context of a tricategory, and from the strict version to the weak version in the sense of weakening equalities to isomorphisms. A pseudomonad is defined within a tricategory as a pseudomonoid in an endohom bicategory. They have been introduced in \cite{Mar99}, and further studied in \cite{Lack00}. Here we consider pseudomonads within a tricategory whose homs are strict 2-categories. In this situation we introduce pseudomodules between pseudomonads, which also form a tricategory with strict 2-categories as homs. A version of a tricategory whose homs are strict 2-categories is an enriched bicategory of \cite{GSh13} with the enrichment in the 2-category of categories. An alternative is an unbiased version of this notion, that is, a $\twocat$ enriched bicategory which instead of binary compositions has specified multifold compositions of any length, associative in the suitable sense up to isomorphisms. 

Further we assume that $\C$ is an enriched bicategory either in the sense of \cite{GSh13}, or in the sense of its unbiased analogue.  We will write as if it were a Gray category. A pseudomonad $\TT = (X, T)$ in $\C$ consists of an object $X$, an endomorphism $T : X \rightarrow X$, 2-cells $T^2 \rightarrow T$ and $1_X \rightarrow T$ and invertible 3-cells expressing associativity and unitivity, satisfying the usual coherence axioms. A pseudomodule $\MM$ from a pseudomonad $\TT = (X, T)$ to a pseudomonad $\SS = (Y, S)$ is defined to consist of a morphism $M : X \rightarrow Y$, left and right pseudoaction 2-cells
\[
\xymatrix{
&X  \ar[rd]^M="2" \ar[ld]_M="1"  &\\
Y \ar[rr]_S& \ar@{=>}"1"+(6,-4);"1"+(11,-4) ^l& Y
}
\qquad\qquad\qquad
\xymatrix{
&Y  \ar[rd];[]_M="2" \ar[ld];[]^M="1"  &\\
X \ar[rr]_T& \ar@{=>}"1"+(11,-4); "1"+(6,-4) _r& X
}
\]

\noindent and pseudoaction isomorphisms, the invertible 3-cells
\[
\xymatrix{TTM \ar@{=>}[rr]^{Tr} \ar@{=>}[d]_{mM} \ar@{}[rrd]|{\hiso}&& TM \ar@{=>}[d]^{r}\\
TM \ar@{=>}[rr]_{r} && M}
\qquad\qquad
\xymatrix{&TM \ar@{=>}[rd]^{r} \ar@{}[d]|-{\viso} &\\
M \ar@{=>}[rr]_{1_M} \ar@{=>}[ru]^{eM} && M
}
\]
\[
\xymatrix{MSS \ar@{=>}[rr]^{lS} \ar@{=>}[d]_{Mm} \ar@{}[rrd]|{\hiso}&& MS \ar@{=>}[d]^{l}\\
MS \ar@{=>}[rr]_{l} && M}
\qquad\qquad
\xymatrix{&MS \ar@{=>}[rd]^{l} \ar@{}[d]|-{\viso} &\\
M \ar@{=>}[rr]_{1_M} \ar@{=>}[ru]^{Me} && M
}
\]
\[
\xymatrix{TMS \ar@{=>}[rr]^{Tl} \ar@{=>}[d]_{rS} \ar@{}[rrd]|{\hiso}&&  TM \ar@{=>}[d]^{r}\\
MS \ar@{=>}[rr]_{l} && M}
\]

\noindent satisfying coherence conditions, which can be quite obviously understood in the ``all diagrams commute" way, or obtained in a finitary form from the known coherence theorems. A map between modules $\MM$ and $\NN$ consists of a 2-cell $t : M \rightarrow N$ and suitably coherent invertible 2-cells:
\[
\xymatrix{TM \ar@{=>}[rr]^{r} \ar@{=>}[d]_{Tt} \ar@{}[rrd]|{\hiso}&&  M \ar@{=>}[d]^{t}\\
TN \ar@{=>}[rr]_{r} && N}
\qquad\qquad 
\xymatrix{MS \ar@{=>}[rr]^l \ar@{=>}[d]_{tS} \ar@{}[rrd]|{\hiso}&&  M \ar@{=>}[d]^{t}\\
NS \ar@{=>}[rr]_{l} && N.}
\]

\noindent A morphism between pseudomodule maps $t$ and $s$ consists of a 3-cell $\alpha : t \rightarrow s : M \rightarrow N$ satisfying the obvious conditions. Pseudomodules between any pair of pseudomonads $\TT$ and $\SS$, maps between them, and their morphisms form a strict 2-category $\mod(\C)(\TT, \SS)$. Pseudomodules can be horizontally composed in a fairly standard way once $\C$ has certain cocompletness properties. We give a somewhat heuristic description of the process. Suppose that hom-2-categories of $\C$ have strict 2-coequalizers of split pairs. In other words 1-categorical coequalizers of such pairs exist, and they are taken into equalizers in $\twocat$ by the contravariant hom functors. Suppose also that these coequalizers  are preserved by left and right compositions with any fixed 1-cell. Under these conditions, if $\MM$ is a pseudomodule from $\RR$ to $\TT$, and  $\NN$ is a pseudomodule from $\TT$ to $\SS$,  a composite module $\NN\circ \MM$ from $\RR$ to $\SS$ is defined to consists of the coequalizer $N\circ M$ as in  
\[\xymatrix{
&X \ar[rr]^{T}_{\;}="2"&& X \ar[rd]^N&\\
Z \ar[ru]^{M} \ar[rr]^M \ar@/_2pc/[rrrr]_{N\circ M}^{\;}="1" && X \ar[rr]^{N} \ar@{=>} "1"+(0,5); "1"+(0,1) \ar@{=>} "2"+(-2,-3); "2"+(-2,-8)  _{rM}  \ar@{=>} "2"+(2,-3); "2"+(2,-8) ^{Nl} && Y.} 
\]

\noindent  with pseudoaction 2-cells $R(N\circ M) \rightarrow (N\circ M)$ and $(N\circ M)T \rightarrow (N\circ M)$  and the pseudoaction isomorphisms for them induced by the pseudoaction 2-cells $MR \rightarrow M$ and $SN \rightarrow N$, and their pseudoaction isomorphisms. Moreover, $\circ$ extends to 2-functors:
\[\xymatrix{\mod(\C)(\RR, \TT)\times \mod(\C)(\TT, \SS)  \ar[r] & \mod(\C)(\RR, \SS).}\]

\noindent Associativity and unitivity isomorphisms of $\C$ induce invertible module maps expressing associativity and unitivity for the operation $\circ$. These module maps are functorial, and satisfy coherence conditions. It follows that, with $\circ$ as a horizontal composition, pseudomonads and the 2-categories of pseudomodules between them form a $\twocat$-enriched bicategory. We denote this by $\mod(\C)$. Alternatively, it is possible to define a multifold version of the operation $\circ$, which will lead to an unbiased version of $\mod(\C)$. 

\section{Biprofunctors and modules}\label{biprofmod} 

\smallskip We consider a special case of the previous section, taking $\C$ to be the $\twocat$-enriched bicategory $\mat(\cat)$ of $\cat$-valued matrices. Its objects are small sets. For sets $X$ and $Y$, the 2-category  $\mat(\cat)(X, Y)$ is defined as $[X\times Y, \cat]$. In more details, a morphism $M : X \rightarrow Y$ of $\mat(\cat)$ consists of a collection of categories $M(x, y)$ indexed by elements of $X\times Y$. While, morphisms and 2-cells are given by indexed collections of functors and natural transformations. The horizontal composition is defined by the usual matrix multiplication formula:
\[NM(x, y) = \coprod_z(M(x, z)\times N(z, y)).\]

\noindent Identity morphisms are matrices with the terminal category at the diagonal and the empty category everywhere else. 

A pseudomonad $\AA = (X, A)$ in $\mat(\cat)(X, Y)$ is the same as a bicategory with set of objects $X$, and the homcategories the components of the matrix $A$. We set
\[\twoprof = \mod\big{(}\mat(\cat)\big{)}.\]

\noindent Then,
\[\twoprof(A, B) = \bicat(A^{op}\times B, \cat).\]

\noindent The morphisms of $\twoprof$ can be called biprofunctors. $\twoprof$ itself can be called the \textbf{tricategory of biprofunctors}.

\smallskip Let $\mod$ denote the full sub- 2-cat enriched bicategory of $\twoprof$ whose objects are categories. A morphism in $\mod$ between categories $X$ and $Y$ is the same as a pseudofunctor $X^\op\times Y \rightarrow \cat$, i.e. a module from $X$ to $Y$ in the sense of Section~\ref{equipments}. In this way we have organized modules into a $\twocat$ enriched bicategory. Further, we outline how various data that we have used previously can be shortly described internal to $\mod$. 

A 2-cell of $\mod$
\[\xymatrix{X \ar@/^1.5em/[r]^{M}="1"  \ar@/_1.5em/[r]_{N}& Y \ar@{=>}"1"+(0,-5); "1"+(0,-10) ^F}\]

\noindent is a pseudonatural family of functors $F_{x, y} : M(x, y) \rightarrow N(x, y)$. A 3-cell
\[\xymatrix{X \ar@/^1.5em/[r]^{M}="1"  \ar@/_1.5em/[r]_{N}& Y \ar@{=>}"1"+(0,-5); "1"+(0,-10) ^F}
\quad
\xymatrix{\ar@3{>}[r]^\tau&}
\quad
\xymatrix{X \ar@/^1.5em/[r]^{M}="1"  \ar@/_1.5em/[r]_{N}& Y \ar@{=>}"1"+(0,-5); "1"+(0,-10) ^G}
\]

\noindent  is a modification with components natural transformations $\tau_{x, y} : F_{x, y} \rightarrow G_{x, y}$. A 2-cell $A^{\circ n} \Rightarrow A : X \rightarrow X$ in $\mod$
\[\xymatrix{&X\ar[r]^A&X \ar@{}[r]|-{\hdots}="1" & X \ar[r]^A& X \ar[rd]^A\\
X \ar[rrrrr]_A \ar[ru]^A&&&&&X \ar@{=>} "1"+(0,-5); "1"+(0,-10)
}\]

\noindent amounts to a family of functors
\[\xymatrix{A(x_0, x_1)\times A(x_1, x_2)\times\cdots\times A(x_{n-1}, x_n) \ar[r]& A(x_0, x_n),}\]

\noindent pseudonatural in the first and the last argument, and pseudo-dinatural in all other arguments. 

Consider the embedding:
\beq{cattops}
\xymatrix{\cat^\op  \ar[r] & \mod}
\eeq

\noindent which takes a functor $F : X \rightarrow Y$ to a module $Y(-, F-)$, i.e. the pseudofunctor
\[\xymatrix{Y^\op\times X \ar[r]^{1_{Y^{op}}\times F} & Y^\op \times Y \ar[r]^{\mathrm{Hom}} & \cat}.\]

\noindent Denote this module by $F^\o$. To a natural transformation 
\beq{t}
\xymatrix{X \ar@/^1.5em/[r]^{F}="1"  \ar@/_1.5em/[r]_{G}& Y \ar@{=>}"1"+(0,-5); "1"+(0,-10) ^t}
\eeq

\noindent  (\ref{cattops}) assigns a pseudonatural family of functors $Y(x, t_{y}) : Y(x, Fy) \rightarrow Y(x, Gy)$. We denote this by $t$ again, so the image of (\ref{t}) in $\mod$ becomes
\[\xymatrix{Y \ar@/^1.5em/[r]^{F^\o}="1"  \ar@/_1.5em/[r]_{G^\o}& X \ar@{=>}"1"+(0,-5); "1"+(0,-10) ^t}\]

\noindent For a module $M : Z \rightarrow Y$, and a functor $F : X \rightarrow Y$, the composite $F^\o\circ M$ is the module $Z \rightarrow X$ with categories of vectors $M(z, Fx)$. Furthermore, for a natural transformation $t : F \Rightarrow G : X \rightarrow Y$, $t\circ M$ is a 2-cell of $\mod$ given by the pseudonatural family of functors $M(z, t_x) : M(z, Fx) \rightarrow  M(z, Fy)$. 

$F^\o$ has a (strict) right adjoint $F_\o$ in $\mod$ given by the module $Y(F-, -)$. A natural transformation $t : F \Rightarrow G$ gives rise to a 2-cell $t : F_\o \Rightarrow G_\o$ of $\mod$ given by the pseudonatural family of functors $Y(t_{x}, y)$. This leads to another embedding
\[\xymatrix{\cat^{\mathrm{co}}  \ar[r] & \mod.}\]

\noindent For $M : Z \rightarrow Y$ a module, and $F : X \rightarrow Z$ a functor, the composite  $M\circ G_\o$ is the module $X \rightarrow Y$ with categories of vectors $M(Fx, y)$. For a natural transformation $t : F \Rightarrow G : X \rightarrow Z$, $t\circ M$ is a 2-cell of $\mod$ given by the pseudonatural family of functors $M(t_x, x) : M(Gz, x) \rightarrow M(Fz, x)$. Giving a 2-cell in $\mod$
\beq{F}
\xymatrix{
X \ar[rr]^A="1" \ar[d];[]^{F^\o}&& X \ar[d];[]_{F^\o}\\
Y \ar[rr]_B&& Y \ar@{=>} "1"+(0,-6); "1"+(0,-11)\\
}
\eeq

\noindent is the same as giving its transpose along the adjunction $F^\o \dashv F_\o$
\[
\xymatrix{
X \ar[rr]^A="1" \ar[d]_{F_\o}&& X \ar[d];[]_{F^\o}\\
Y \ar[rr]_B&& Y \ar@{=>} "1"+(0,-6); "1"+(0,-11)\\
}
\]

\noindent Hence a 2-cell of the form (\ref{F}) amounts to a family of functors $F_{x, y} : A(x, y) \rightarrow B(Fx, Fy)$ pseudonatural in both arguments. 

\section{Lax Monads in a 3-category}\label{laxmonads}

\smallskip Suppose that $\C$ is an arbitrary tricategory. We work as if $\C$ were a Gray category. For the composition of 1-cells we use the dot symbol. So for the $n$-fold composite of 
\[
\xymatrix{X_n \ar@/^1pc/[rr]^{A_n}_{\quad}="1" \ar@/_1pc/[rr]_{B_n}^{\quad}="2" & \ar@{=>}"1"+(0,-1); "1"+(0,-6) ^{F_1}& X_{n-1} \ar@{}[r]|{\cdots} & X_2 \ar@/^1pc/[rr]^{A_2}_{\quad}="3" \ar@/_1pc/[rr]_{B_2}^{\quad} ="4" \ar@{=>}"3"+(0,-1); "3"+(0,-6)^{F_2}&& X_{1} \ar@/^1pc/[rr]^{A_1}_{\quad}="5" \ar@/_1pc/[rr]_{B_1}^{\quad}="6" & \ar@{=>}"5"+(0,-1); "5"+(0,-6) ^{F_n}& X_0}
\]
\noindent we write
\[\xymatrix{A_1.A_2.\cdots.A_n \ar@{=>}[rr]^{F_1.F_2.\cdots.F_n} && B_1.B_2.\cdots.B_n.}\]

\noindent For vertical compositions of 2-cells we use concatenation, so the $n$-fold composite
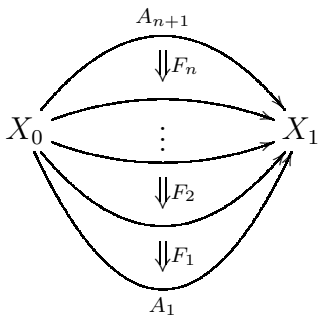
\[
\xymatrix{X_0 \ar@{}[rrr]|-{\vdots} \ar@/_5pc/[rrr]_{A_1}="1" \ar@/_3pc/[rrr]^{}="2" \ar@/_1pc/[rrr]_{}="3" \ar@/^3pc/[rrr]^{A_{n+1}}_{}="4" \ar@/^1pc/[rrr]^{}="5"  \ar@{=>}"2"+(0,-2) ; "2"+(0,-6) ^{F_1}  \ar@{=>}"3"+(0,-2) ; "3"+(0, -6) ^{F_2} \ar@{=>}"4"+(0,-2) ; "4"+(0, -6) ^{F_n}&&& X_1} 
\]

\noindent is denoted by $F_1F_2\cdots F_n : A_{n+1} \Rightarrow A_1$. All pasting composites of 2-cells that we form are obtained through consecutive application of the $n$-fold horizontal and vertical composites. A 2-dimensional pasting diagram may be evaluated to a 2-cell in several possible ways. Between any two values of the same diagram there is a unique structural invertible 3-cell. These structural 3-cells will be denoted by $\cong$, or in some definitions they will be omitted altogether. A 2-cell of the form
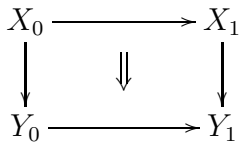
\[
\xymatrix{X_0 \ar[rr]_{}="1" \ar[d] && X_1  \ar[d]\\
Y_0 \ar[rr]^{}="2" && Y_1 \ar@{=>} "1"+(0, -4); "1"+(0, -9)
}
\]

\noindent will be called a square. Given squares in either of the following configurations
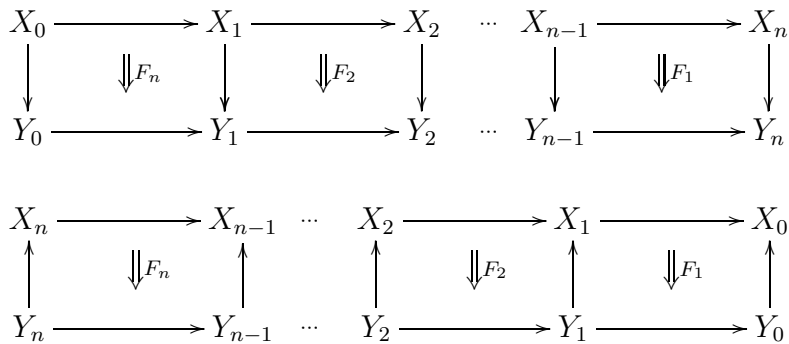
\[
\xymatrix{X_0 \ar[rr]_{}="1" \ar[d] && X_{1}  \ar[d] \ar[rr]_{}="3" \ar[d]  && X_{2} \ar[d] \ar@{}[r]|{\cdots} & X_{n-1} \ar[d] \ar[rr]_{}="5" \ar[d] && X_n \ar[d]\\
Y_0 \ar[rr]^{}="2" && Y_{1}   \ar[rr]^{}="4" & & Y_2 \ar@{}[r]|{\cdots}& Y_{n-1} \ar[rr]^{}="6" && Y_n \ar@{=>} "1"+(0, -4); "1"+(0, -9) ^{F_n} \ar@{=>} "3"+(0, -4); "3"+(0, -9) ^{F_2} \ar@{=>} "5"+(0, -4); "5"+(0, -9) ^{F_1}
}
\]
\[
\xymatrix{X_n \ar[rr]_{}="1" \ar[d];[] && X_{n-1}  \ar[d];[] \ar@{}[r]|{\cdots}&  X_2  \ar[d];[]   \ar[rr]_{}="3" && X_1  \ar[d];[] \ar[rr]_{}="5" && X_0 \ar[d];[]\\
Y_n \ar[rr]^{}="2" && Y_{n-1} \ar@{}[r]|{\cdots}&  Y_2  \ar[rr]^{}="4" && Y_1 \ar[rr]^{}="6" && Y_0 \ar@{=>} "1"+(0, -4); "1"+(0, -9) ^{F_n} \ar@{=>}   "3"+(0, -4); "3"+(0, -9) ^{F_2} \ar@{=>} "5"+(0, -4); "5"+(0, -9) ^{F_1}
}
\]

\noindent the composite, will be denoted by $F_1F_{2}\cdots F_n$. This notation is conflicting with the earlier adoption of concatenation for the vertical composition of 2-cells, but we will use it only when it is clear from the context what is meant. We refer to homomorphisms of tricategories as functors.

\smallskip A lax monad in a tricategory is a lax monoid in the sense of \cite{DS03} in an endohom monoidal bicategory. In  \cite{DS03} a lax monoid was defined in a packaged form, as a strict monoidal lax functor with the domain the simplicial category.  We state the definition in an unpacked form.

\begin{definition} A (normal) \textbf{lax monad} or an \textbf{l-monad} $\AA = (A_0, A, P^A, \xi^A)$ in $\C$ consists of an object $A_0$, a 1-cell $A : A_0 \rightarrow A_0$,  for every $n > 0$, a 2-cell $P^A_n : A^{.n} \Rightarrow A$, with $P^A_1$ an identity, and for every partition $m = n_1+ \cdots+ n_j$ a 3-cell 

\begin{equation}\label{xi}
\xymatrix{\xi^A_{ n_1, ..., n_j} : P^A_m \ar@3{->}[r] & P^A_j(P^A_{n_1}. \cdots .P^A_{n_j}) : A^{.n_1}... A^{.n_j} \ar@{=>}[r] & A,
}
\end{equation}

\noindent satisfying the coherence condition:

\begin{equation}\label{lmonad}
\xymatrix{P^A_l 
\ar@3{->}[r]^{\xi^A_{m_1, ..., m_k}} 
\ar@3{->}[dd]_{\xi^A_{n_{11}, ..., n_{kj_{k}}}} 
& P^A_k(P^A_{m_1}. \cdots .P^A_{m_j})
\ar@3{->}[d]^{-(\xi^A_{ n_{11}, ..., n_{1j_1}}. \cdots .\xi^A_{ n_{k1}, ..., n_{kj_k}})} 
\\ 
& P^A_k((P^A_{j_1}(P^A_{n_{11}}. \cdots .P^A_{n_{1j_1}})). \cdots .(P^A_{j_k}(P^A_{n_{k1}}. \cdots .P^A_{n_{kj_k}})))
\ar@3{->}[d]^{\cong}
\\
P^A_h(P^A_{n_{11}}. \cdots .P^A_{n_{kj_k}})
\ar@3{->}[r]_{\xi^A_{j_1, ..., j_k}-} 
&
P^A_k(P^A_{j_1}. \cdots .P^A_{j_k})(P^A_{n_{11}}. \cdots .P^A_{n_{kj_k}})
}
\end{equation}

\noindent for $l = m_1 +\cdots+m_k$, $m_i = n_{i1}+\cdots+n_{ij_i}$, $h = j_1+\cdots+j_k$, as well as $\xi_{1, ..., 1}$ and $\xi_n$ required to be identities. 
\end{definition}

An l-monad in a tricategory is a lax version of a monad in a 2-category. Generally, given a notion defined within a 2-category which consists of a data of 0-, 1-, and 2-cells, satisfying axioms which are equalities between 2-cells, by the \textit{lax version of the 2-categorical notion within a tricategory} we mean a structure defined within a tricategory which consists of the same data of 0-, 1- and 2-cells, together with non-invertible 3-cells which replace the 2-cell equations of the 2-categorical notion, and are required to satisfy coherence axioms. Our lax monad is obtained in this way from the unbiased presentation of a monad in a 2-category, i.e. the one in which the $n$-fold multiplication 2-cells are regarded as a part of the data. A \textbf{colax monad} is another lax version of an unbiased monad in which the associator 3-cells take the opposite direction. We will not use this notion in this paper.

\begin{example} In the light of the observations made in Section~\ref{biprofmod} it is easy to verify:
\begin{theorem} An l-monad in $\mod$ is the same as an equipment. \end{theorem}
\end{example}

In the following definition we collect definitions of the lax versions of monad upmaps and their transformations.

\begin{definition}
An \textbf{l-upmap} of l-monads  $\FF = (F_0, F, \kappa^F) : \BB \rightarrow \AA$ consists of a morphism $F_0 : B_0 \rightarrow A_0$, a square $F : A.F_0 \Rightarrow F_0.B$, and for every $n \geq 0$ a 3-cell 
\[
\xymatrix{
B_0 \ar@/_3pc/[dddd]_B="7" \ar[rr]^{F_0} \ar[d]^{B}="2" && A_0  \ar[d]^A_{}="1"\\
B_0 \ar[rr]^{F_0} \ar[d]^{B}_{}="4" && A_0  \ar[d]^A_{}="3"\\
B_0 \ar[rr]_{F_0} \ar@{}[d]|{\vdots} \ar@{=>} {}+(-5,0) ; {}+(-10,0) _{P^B_m}&& A_0 \ar@{}[d]|{\vdots}\\
B_0 \ar[rr]^{F_0} \ar[d]^B_{}="6" && A_0  \ar[d]^A_{}="5"\\
B_0 \ar[rr]_{F_0} && A_0 \ar@{=>} "1"+(-10,0); "1"+(-15,0)_{F}  \ar@{=>} "3"+(-10,0); "3"+(-15,0) _F   \ar@{=>} "5"+(-10,0); "5"+(-15,0)_{F} 
}
\qquad
\xymatrix@R=1.4em{\\  \\  \\ \\ \ar@3{->}[r] ^{\kappa^F_n}& }
\qquad
\xymatrix{
B_0  \ar[rr]^{F_0} \ar[dddd]_{B}^{}="2" && A_0  \ar@/_3pc/[dddd]_A="1" \ar[d]^A\\
&& A_0  \ar[d]^A\\
&& A_0 \ar@{}[d]|{\vdots} \ar@{=>} {}+(-5,0) ; {}+(-10,0) _{P^A_m}\\
&& A_0  \ar[d]^A\\
B_0 \ar[rr]_{F_0} && A_0 \ar@{=>} "1"+(-3,0); "1"+(-8,0) _{F}
}
\]

\noindent called a lax comparison map, which for every $m = n_1+\cdots+n_i$ satisfy the equation

\begin{equation}\label{lmap}
\xy 
(0,0)*{F(P^A_m.F_0)}="1"; 
(15,15)*{(F_0.P^B_m)(F^{n_1}\cdots F^{n_j})}="2"; 
(95, 15)*{(F_0.(P^B_j(P^B_{n_1}. \cdots .P^B_{n_j})))(F^{n_1}\cdots F^{n_j})}="3"; 
(105, 0)*{(F_0.P^B_j)(((F_0.P^B_{n_1})F^{n_{1}})\cdots ((F_0.P^B_{n_j})F^{n_{j}}))}="4"; 
(0,-15)*{F((P^A_j(P^A_{n_1}. \cdots .P^A_{n_j})).F)}="5"; 
(15,-30)*{F(P^A_j.F)(P^A_{n_1}.\cdots.P^A_{n_j}.F_0)}="6"; 
(90,-30)*{(F_0.P^B_j)F^{j}(P^A_{n_1}.\cdots.P^A_{n_j}.F_0)}="7"; 
(105,-15)*{(F_0.P^B_j)((F(P^A_{n_1}.F_0)) \cdots (F(P^A_{n_j}.F_0)))}="8";  
{\ar@3{->}_>>>>>>{\kappa^F_m} "2";"1"+(5,3)} 
{\ar@3{->}^>>>>>>>>>>>>>{(-.(\xi^B_{n_1, ..., n_j}))-} "2";"3"} 
{\ar@3{->}^>>>>>>{\cong} "3";"4"}
{\ar@3{->}_{-(\xi^A_{n_1, .., n_j}.-)} "1";"5"}  
{\ar@3{->}^>>>>>{\cong} "5";"6"}
{\ar@3{->}^<<<<<<<<<<<<{\kappa^F_j-} "7"+(-31,0);"6"+(18,0)}
{\ar@3{->}^<<<<<{\cong} "8"+(-7,-3);"7"} 
{\ar@3{->}^{-(\kappa^F_{n_1} \cdots \kappa^F_{n_j})} "4"+(-7,-3);"8"+(-7,3)} 
\endxy 
\end{equation}

\noindent while $\kappa^F_1$ is required to be an identity. A \textbf{cl-upmap} between l-monads is defined similarly except that its lax comparison maps $\kappa^F_n$-s take the opposite direction and satisfy the axiom obtained from (\ref{lmap}) by reversing in them the arrows involving $\kappa$-s. 

An \textbf{l-transformation of l-upmaps} of l-monads $\NN = (N, \nu^N): \FF \rightarrow \GG : \BB \rightarrow \AA$ consists of a 2-cell  $N : F_0 \Rightarrow G_0$, and a 3-cell 
\[
\xymatrix{
B_0  \ar[rr]^B \ar[dd]_{F_0}="3" && B_0  \ar@/^2.5pc/[dd]^{G_0}="1" \ar[dd]_{F_0}="2"\\
&& \\
A_0 \ar[rr]^A && A_0 \ar@{=>}  "2"+(5,0); "2"+(10,0) ^{N} \ar@{=>}  "3"+(13,0); "3"+(18,0) ^{F}  
}
\quad
\xymatrix@R=1em{\\ \\ \ar@3{->}[r] ^{\nu^N}& }
\quad
\xymatrix{
B_0 \ar@/_2.5pc/[dd]_{F_0}="3" \ar[rr]^B \ar[dd]^{G_0}="2" && B_0  \ar[dd]^{G_0}="1"\\
&&\\
A_0 \ar[rr]_A && A_0 \ar@{=>} "2"+(8,0); "2"+(13,0) ^{G} \ar@{=>} "3"+(5,0); "3"+(10,0) ^{N}  
}
\]

\noindent satisfying the axiom

\begin{equation}\label{ltrans}
\xy 
(0,0)*{(G_0.P^B_n)(N.B^{.n})F^{n}}="1"; 
(0, 20)*{(N.B)(F_0.P^B_n)F^{n}}="2"; 
(50, 20)*{(N.B)F(P^A_n.F_0)}="3"; 
(100, 20)*{G(A.N)(P^A_n.F_0)}="4"; 
(50, 0)*{(G_0.P^B_n)G^{n}(A^{.n}.N)}="5"; 
(100, 0)*{G(P^A_n.G_0)(A^{.n}.N)}="6"; 
{\ar@3{->}_{\cong} "2";"1"} 
{\ar@3{->}^{\cong} "4";"6"} 
{\ar@3{->}^<<<<<{-\kappa^F_n} "2"+(15,0);"3"+(-15,0)} 
{\ar@3{->}^{\nu^N-} "3"+(11,0);"4"+(-17,0)} 
{\ar@3{->}^{(-\nu^N)^{``n"}} "1"+(13,0);"5"+(-22,0)} 
{\ar@3{->}^{\kappa^G_n-} "5"+(11,0);"6"+(-22,0)} 
\endxy 
\end{equation}

\noindent A \textbf{cl-transformation of l-upmaps} is defined in the same way except that the 3-cell $\nu$ takes the opposite direction, and satisfies an axiom obtained from (\ref{ltrans}) by reversing in it all arrows involving $\nu$-s. In the same diagram, leaving the direction of $\nu$-s unchanged, but reversing the directions of $\kappa$-s, we obtain the axiom for an \textbf{l-transformation of cl-upmaps}. Changing both, directions of $\nu$-s and $\kappa$-s we obtain a notion of \textbf{cl-transformation of cl-upmaps}. 

A \textbf{modification between l-transformations of l-upmaps} $\NN \rightarrow \SS : \FF \rightarrow \GG : \BB \rightarrow \AA$ is a 3-cell $\lambda : N \Rightarrow S$ satisfying the equation

\begin{equation}\label{modification}
\xy 
(0,0)*{(N.B)F}="1"; 
(40,0)*{G(A.N)}="2"; 
(0,-15)*{(S.B)F}="3"; 
(40,-15)*{G(A.S)}="4"; 
{\ar@3{->}^{\nu^N} "1"+(8,0);"2"+(-8,0)}
{\ar@3{->}^{(\lambda.-)-} "2";"4"} 
{\ar@3{->}_{(\lambda.-)-} "1";"3"} 
{\ar@3{->}_{\nu^S} "3"+(5,0);"4"+(-11,0)}  
\endxy 
\end{equation}

\noindent Modifications between all other kinds of transformations are defined similarly. 

There are four tricategories all of which have l-monads as their objects
\[\mll(\C), \quad \mlcol(\C), \quad \mcoll(\C), \quad \mcolcol(\C).\] 

\noindent  The morphism of the first two are l-upmaps, and morphism of the last two are cl-upmaps. The 2-cells for the first and the third are l-transformations, and the 2-cells of the second and the fourth are cl-transformations. The 3-cells for all of them are modifications. 
\end{definition}

Consider next lax versions of monad downmaps and their transformations. Let $(-)^\op$ denote the functor from $\C$ to itself, which inverts the direction of 1-cells and leaves the directions of 2- and 3-cells unchanged.

\begin{definition}
The tricategories 
\[\moll(\C), \quad \molcol(\C), \quad \mocoll(\C), \quad \mocolcol(\C)\]

\noindent are defined by the scheme
\[\mathfrak{L}^{\top-}(\C) = (\mathfrak{L}^{\bot-}(\C^\op))^\op\]

\noindent The morphisms of the first two are called \textbf{l-downmaps of l-monads}. The morphisms of the last two are called \textbf{cl-downmaps of l-monads}. The 2-cells of the first and fourth are called l-transformations of l- and respectively cl-downmaps, and the 2-cells of the second and the fourth are called cl-transformations of l- and respectively cl-downmaps. The 3-cells of all of them are called modifications.
\end{definition}

\smallskip Suppose that $\mathcal{I} : \K \rightarrow \C$ is a functor of tricategories. Define a tricategory $\mll(\K, \C)$ to have the following components: An object consists of an object $A_0$ of $\K$ and an l-monad $(\mathcal{I}(A_0), A, P, \xi)$ in the tricategory $\C$. A morphism consists of a morphism $F_0$ of $\K$ and an l-downmap of l-monads of the form $(\mathcal{I}(F_0), F, \kappa)$. A 2-cell consists of a 2-cell $N$ of $\K$ and an l-transformation of the form $(\mathcal{I}(N), \nu)$. A 3-cell consists of a 3-cell $\lambda$ of $\K$, such that $P(\lambda)$ is a modification. Analogously one defines tricategories $\mlcol(\K, \C)$, $\mcoll(\K, \C)$, $\mcolcol(\K, \C)$, $\molcol(\K, \C)$, $\mocoll(\K, \C)$ and $\mocolcol(\K, \C)$. We will say that this is a context \textit{relative} to the functor $\mathcal{I}$.

\begin{example} Taking $\mathcal{I}$ to be the functor $\cat^\op \rightarrow \mod$ of Section~\ref{biprofmod}, we get a 2-category $\mll(\cat^\op, \mod)$. Using the observations made in Section~\ref{biprofmod} it is not difficult to see that:

\begin{theorem}\label{LK}
The 2-category $\mll(\cat^\op, \mod)^\op$ is isomorphic to the 2-category $\Eq$.
\end{theorem}

\noindent To be more explicit, an equipment functor $\FF = (F_0, F, \kappa^F) : \AA \rightarrow \BB$, as in Definition~\ref{eqfunct}, becomes an l-upmap $(F_0^\o, F, \kappa^F) : \BB \rightarrow \AA$ of l-monads in $\mod$ consisting of the morphism $F_0^\o : B_0 \rightarrow A_0$, the 2-cell $F : F^\o_0\circ B \rightarrow A\circ F^\o_0$, and the lax comparison 3-cells $\kappa^F$. A transformation of equipment functors $\tt = (t, \nu^t) : \FF \rightarrow \GG$, as in Definition~\ref{eqfuncttrans}, becomes an l-transformation of l-upmaps of l-monads in $\mod$ consisting of the 2-cell $t : F_0^\o \rightarrow G_0^\o$ and the 3-cell $\nu^t$. Since $\cat$ is a 2-category, modifications in the given situation are trivial.

Furthermore, the formula $\Eq^{-} = \mathfrak{L}^{\bot-}(\cat^\op, \mod)$ defines 2-categories
\[\Eq^{\circ\circ}, \quad \Eq^{\circ\bullet}, \quad \Eq^{\bullet\circ}, \quad \Eq^{\bullet\bullet}.\]

\noindent The first of them is just another name for $\Eq$. $\Eq^{\circ\bullet}$ is the 2-category of equipments, lax functors and colax transformations of lax functors. $\Eq^{\bullet\circ}$ is the 2-category of equipments, colax functors and lax transformation of colax functors. $\Eq^{\bullet\bullet}$ is the 2-category of equipments, colax functors and colax transformations of colax functors.
\end{example}

\smallskip Consider an l-monad in $\mll(\C)$. It consists of the data $(\AA, \BB, \mathbb{P}^\BB, \xi^\BB)$, where $\AA$ is an l-monad in $\C$, $\BB = (B_0, B, \kappa^B) : \AA \rightarrow \AA$ is an l-upmap of l-monads, $\mathbb{P}^\BB$ is a family of l-transformations $\mathbb{P}^\BB_n = (P^B_n, \nu^{P^B_n}) : \BB^n \rightarrow \BB$, and $\xi^\BB$ is a family of modifications of l-transformations determined by a family of 3-cells $\xi^B$ of $\C$. The four-tuple $(A_0, B_0, P^B_n, \xi^B)$ defines an l-monad in $\C$, which we denote by $\BB$. Furthermore, the data $(A, B, \nu^{P^B_n})$ becomes an l-downmap of l-monads $\BB \rightarrow \BB$, which via $\mathbb{P}^A$, $\kappa^A$ and $\xi^A$ becomes an l-monad in $\moll(\C)$. Through this correspondence, an l-monad in $\mll(\C)$ is the same as an l-monad in $\moll(\C)$. Renaming $B_0$ into $B$, renaming the old $B$ into $D$, and renaming $\nu^{P^B_n}$ into $\kappa^A$, in the next paragraph we present the same structure under a new name.

An \textbf{ll-distributive pair of l-monads} consists of the data $(\BB, \AA, D, \kappa^B, \kappa^A)$, where $\AA$ and $\BB$ are l-monads in $\C$ whose base objects are equal $A_0 = B_0$, $D$ is a square $A.B \rightarrow B.A$, and $\kappa^B$ and $\kappa^A$ are 3-cells

\begin{equation}\label{dsk}
\xymatrix{
A_0 \ar@/_3pc/[dddd]_A="7" \ar[rr]^{B} \ar[d]^{A}="2" && A_0  \ar[d]^A_{}="1"\\
A_0 \ar[rr]^{B} \ar[d]^{A}_{}="4" && A_0  \ar[d]^A_{}="3"\\
A_0 \ar[rr]_{B} \ar@{}[d]|{\vdots} \ar@{=>} {}+(-5,0) ; {}+(-10,0) _{P^A_m}&& A_0 \ar@{}[d]|{\vdots}\\
A_0 \ar[rr]^{B} \ar[d]^A_{}="6" && A_0  \ar[d]^A_{}="5"\\
A_0 \ar[rr]_{B} && A_0 \ar@{=>} "1"+(-10,0); "1"+(-15,0)_{D}  \ar@{=>} "3"+(-10,0); "3"+(-15,0)_{D}   \ar@{=>} "5"+(-10,0); "5"+(-15,0)_{D} 
}
\qquad
\xymatrix@R=1.4em{\\  \\  \\ \\ \ar@3{->}[r] ^{\kappa^B_n}& }
\qquad
\xymatrix{
A_0  \ar[rr]^{B} \ar[dddd]_{A}^{}="2" && A_0  \ar@/_3pc/[dddd]_A="1" \ar[d]^A\\
&& A_0  \ar[d]^A\\
&& A_0 \ar@{}[d]|{\vdots} \ar@{=>} {}+(-5,0) ; {}+(-10,0) _{P^A_m}\\
&& A_0  \ar[d]^A\\
A_0 \ar[rr]_{B} && A_0 \ar@{=>} "1"+(-3,0); "1"+(-8,0) _{D}
}
\end{equation}

\begin{equation}\label{dsn}
\xymatrix{
A_0 \ar[rr]^{A} \ar[d]_{B}^{}="2" && A_0  \ar[d]_B_{}="1" \ar@/^3pc/[dddd]^B="7"\\
A_0 \ar[rr]^{A} \ar[d]_{B}_{}="4" && A_0  \ar[d]_B_{}="3"\\
A_0 \ar[rr]_{A} \ar@{}[d]|{\vdots} && A_0 \ar@{}[d]|{\vdots} \ar@{=>} {}+(5,0) ; {}+(10,0) ^{P^B_m}\\
A_0 \ar[rr]^{A} \ar[d]_B_{}="6" && A_0  \ar[d]_B_{}="5"\\
A_0 \ar[rr]_{A} && A_0 \ar@{=>} "2"+(10,0); "2"+(15,0)  ^{D} \ar@{=>} "4"+(10,0); "4"+(15,0)^{D}  \ar@{=>}  "6"+(10,0); "6"+(15,0) ^{D} 
}
\qquad
\xymatrix@R=1.4em{\\  \\  \\ \\ \ar@3{->}[r] ^{\kappa^A_n}& }
\qquad
\xymatrix{
A_0  \ar[rr]^{A}  \ar@/^3pc/[dddd]^{B}="1" \ar[d]_B && A_0  \ar[dddd]^{B}="2" \\
A_0   \ar[d]_B&&\\
A_0  \ar@{}[d]|{\vdots} \ar@{=>} {}+(5,0) ; {}+(10,0) ^{P^B_m} &&\\
A_0  \ar[d]_B&&\\
A_0 \ar[rr]_{A} && A_0. \ar@{=>} "1"+(3,0); "1"+(8,0) ^{D}
}
\end{equation}

\noindent satisfying axioms expressing that the data defines an $l$-monad in $\mll(\C)$, or equivalently an $l$-monad in $\moll(\C)$, in the way pointed out above.

An ll-distributive pair is a lax version of the 2-categorical notion of distributive pair of monads recounted in Section~\ref{monads}. Other lax distributive pairs of l-monads are obtained by varying the directions of $\kappa^A$ and $\kappa^B$. These correspond to lax monads in $\mll(\C)$, $\mcoll(\C)$ and $\mlcol(\C)$, or equivalently in $\moll(\C)$, $\molcol(\C)$ and $\mocoll(\C)$ respectively. 

An l-monad in $\mll(\K, \C)$ can be thought of as an ll-distributive pair of an l-monad in $\K$ and an l-monad in $\C$.

\begin{example} Since $\mll(\cat^\op, \mod) = \Eq^\op$ is a 2-category, an l-monad in it is the same as a monad in it, which by definition is a $T$-equipment. Thus, we can regard a $T$-equipment $(\TT_0, \AA)$ as an ll-distributive pair consisting of the monad $\TT_0$ in $\cat$ and the l-monad $\AA$ in $\mod$. The $T$-equipment given by the data $(\TT_0, \AA, T, \kappa^T, \nu^m, \nu^e)$, as in Definition~\ref{Tequipments}, besides the l-monad $\AA$ in $\mod$ and the monad $\TT_0$ in $\cat$, consists of a square $T : A\circ T_0^\o \Rightarrow T_0^\o\circ A$ in $\mod$, and the 3-cells of $\mod$
\[
\xymatrix{
A_0 \ar@/_3pc/[dddd]_A="7" \ar[rr]^{T_0^\o} \ar[d]^{A}="2" && A_0  \ar[d]^A_{}="1"\\
A_0 \ar[rr]^{T_0^\o} \ar[d]^{A}_{}="4" && A_0  \ar[d]^A_{}="3"\\
A_0 \ar[rr]_{T_0^\o} \ar@{}[d]|{\vdots} \ar@{=>} {}+(-5,0) ; {}+(-10,0) _{P^A_m}&& A_0 \ar@{}[d]|{\vdots}\\
A_0 \ar[rr]^{T_0^\o} \ar[d]^A_{}="6" && A_0  \ar[d]^A_{}="5"\\
A_0 \ar[rr]_{T_0^\o} && A_0 \ar@{=>} "1"+(-10,0); "1"+(-15,0)_{T}  \ar@{=>} "3"+(-10,0); "3"+(-15,0)_{T}   \ar@{=>} "5"+(-10,0); "5"+(-15,0)_{T} 
}
\qquad
\xymatrix@R=1.4em{\\  \\  \\ \\ \ar@3{->}[r] ^{\kappa^T_n}& }
\qquad
\xymatrix{
A_0  \ar[rr]^{T_0^\o} \ar[dddd]_{A}^{}="2" && A_0  \ar@/_3pc/[dddd]_A="1" \ar[d]^A\\
&& A_0  \ar[d]^A\\
&& A_0 \ar@{}[d]|{\vdots} \ar@{=>} {}+(-5,0) ; {}+(-10,0) _{P^A_m}\\
&& A_0  \ar[d]^A\\
A_0 \ar[rr]_{T_0^\o} && A_0 \ar@{=>} "1"+(-3,0); "1"+(-8,0) _{T}
}
\]

\begin{equation}
\xymatrix{
A_0 \ar[rr]^{A} \ar[d]_{T_0^\o}^{}="2" && A_0  \ar[d]_{T_0^\o}_{}="1" \ar@/^3pc/[dd]^{T^\o_0}="7"\\
A_0 \ar[rr]^{A} \ar[d]_{T_0^\o}_{}="4" && A_0  \ar[d]_{T_0^\o}_{}="3" \ar@{=>} {}+(5,0) ; {}+(10,0) ^{m}\\
A_0 \ar[rr]_{A} && A_0 \ar@{=>} "2"+(10,0); "2"+(15,0)  ^{T} \ar@{=>} "4"+(10,0); "4"+(15,0)^{T}\\
}
\quad
\xymatrix@R=1.4em{\\  \\  \ar@3{->}[r] ^{\nu^m}& }
\quad
\xymatrix{
A_0  \ar[rr]^{A}  \ar@/^3pc/[dd]^{T_0^\o}="1" \ar[d]_{T_0^\o} && A_0  \ar[dd]^{T_0^\o}="2" \\
A_0   \ar[d]_{T_0^\o} \ar@{=>} {}+(5,0) ; {}+(10,0) ^{m}&&\\
A_0  \ar[rr]^A&& A_0\\
\ar@{=>} "1"+(3,0); "1"+(8,0) ^{T}
}
\end{equation}

\begin{equation}
\xymatrix{
A_0 \ar[rr]^{A} && A_0  \ar[dd]_{T_0^\o}_{}="1" \ar@/^3pc/[dd]^{T_0^\o}="7"\\
&&\\
&& A_0 \ar@{=>} "1"+(4,0) ; "1"+(9,0) ^{e}
}
\quad
\xymatrix@R=2em{\\ \ar@3{->}[r] ^{\nu^e}& }
\quad
\xymatrix{
A_0  \ar[rr]^{A}  \ar@/^3pc/[dd]^{T_0^\o}="1" \ar[dd]_{T_0^\o} && A_0  \ar[dd]^{T_0^\o}="2" \\
\ar@{=>} {}+(4,0) ; {}+(9,0) ^{e}&&\\
A_0  \ar[rr]_A && A_0
\ar@{=>} "1"+(3,0); "1"+(8,0) ^{T}
}
\end{equation}

Furthermore, lax monads in $\Eq^{\circ\bullet}$, $\Eq^{\bullet\circ}$ and $\Eq^{\circ\bullet}$ give new notions of $T$-equipments, in which varying laxities for the equipment functor $\TT$ on the one hand, and the equipment functor transformations $\mm$ and $\ee$ on the other hand are allowed.

\end{example}

\smallskip Now we will look at the lax versions of morphisms of distributive pairs of monads as well as their transformations. In contrast with the 2-categorical context, in the lax situation, not all of these notions correspond to lax down/upmaps of l-monads and lax transformations in a tricategory of l-monads. However, for those which do, such as morphisms and 2-cells of the tricategory $\mll\mll(\C)$, the axioms are readily available.

Let us describe $\mll\mll(\C)$. Its objects are ll-distributive pairs of l-monads. A morphism $(\BB', \AA', D, \kappa^{B'}, \kappa^{A'}) \rightarrow (\BB, \AA, D, \kappa^{B}, \kappa^{A})$ amounts to a triple $(\FF, \GG, \delta)$, where $\FF : \AA' \rightarrow \AA$ and $\GG : \BB' \rightarrow \BB$ are l-upmaps with $F_0 = G_0$, and $\delta$ is a 3-cell
\beq{delta}
\xymatrix@R=0.8em{
&& A'_0 \ar[drr]^{B'}_\;="2"&&\\
A'_0  \ar[rru]^{A'}_{\;}="1" &&&& A'_0  \\
&& A_0  \ar[uu];[]^{F_0} \ar[drr]^{B}="4" && \\
A_0  \ar[urr]^{A}="3" \ar[uu];[]_{F_0}  \ar[drr]_{B}&&&&A_0 \ar[uu];[]^{F_0}  \\
&&A_0 \ar[urr]_{A} \ar@{=>} {}+(0,6); {}+(0,11) ^D && \ar@{=>}"3"+(0,5);"3"+(0,10) ^{F} \ar@{=>}"4"+(0,5);"4"+(0,10) ^{G}
}
\quad
\xymatrix@R=1em{\\\\\\ \ar@3{->}[r] ^\delta & \\ \\}
\quad
\xymatrix@R=0.8em{
&& A'_0  \ar[drr]^{B'}&&\\
A'_0  \ar[rru]^{A'}  \ar[drr]^{B'}_\;="1"&&&& A'_0  \\
&& A'_0 \ar[urr]^{A'}_\;="2" \ar@{=>}  {}+(0,6); {}+(0,11) ^{D'} && \\
A_0  \ar[uu];[]_{F_0} \ar[drr]_{B}^\;="3"&&&&A_0 \ar[uu];[]^{F_0} \\
&&A_0 \ar[urr]_{A}^{\;}="4"   \ar[uu];[]^{F_0}&& \ar@{=>} "3"+(0,5);"3"+(0,10) ^{G} \ar@{=>} "4"+(0,5);"4"+(0,10)  ^F
}
\eeq

\noindent satisfying two equation expressing compatibility with the multiplication structures of $\AA$ and $\BB$:
{\small
\beq{axdelta}
\xymatrix@R=0.8em@C=0.1em{
&&(F_0.B'.P_n)(F_0.{D'}^n)(F^n.B')(A^{.n}.G)\ar[drr]^>>>>>>>{-\kappa^B_n-}&&\\
(G.P_n)(B.F^n)(D^n.F_0) \ar[rru]^<<<<<{``{\delta^n}"} &&&&  (F_0.D')(F_0.P_n.B')(F^n.B')(A^{.n}.G)\\
&&&& \\
(GF)(B.P_n.F_0)(D^n.F_0)  \ar[uu];[]_{-\kappa^F_n-}  \ar[drr]_{-\kappa^B_n-}&&&&(F_0.D')(FG)(F_0.P_n.B) \ar[uu];[]^{-\kappa^F_n-}  \\
&& (GF)(D.F_0)(P_n.B.F_0)  \ar[urr]_{-\delta-}
}
\eeq
}
{\small
\[
\xymatrix@R=0.8em@C=0.1em{
&& (F_0.P_n)(F_0.{D'}^n)(F.B^n)(A.G^n) \ar[drr]^>>>>>>>{-\kappa^A_n-}&&\\
(P_n.A')(G^n.A')(B^n.F)(D^n.F_0) \ar[rru]^<<<<<<<<<{``\delta^n"} &&&&  (F_0.D')(A'.P_n)(F.B^n)(A.G^n) \\
&&&& \\
(GF)(P_n.A.F_0)(D^n.F_0)  \ar[uu];[]_{-\kappa^G_n-}  \ar[drr]_{-\kappa^A_n-}&&&& (F_0.D')(FG)(A.P_n.F_0)  \ar[uu];[]^{-\kappa^G_n-}  \\
&& (GF)(D.F_0)(A.P_n.F_0)  \ar[urr]_{-\delta-}
}
\]
}

\noindent A 2-cell $(\FF', \GG', \delta') \rightarrow (\FF, \GG, \delta) : (A', B') \rightarrow (A, B)$ in $\mll\mll(\C)$ amounts to a pair $(\SS, \NN)$ of l-transformations $\NN : \FF' \rightarrow \FF$ and $\SS : \GG' \rightarrow \GG$, such that $N = S : F'_0 \rightarrow F_0$, and the following equation holds:
{\small
\beq{axnunu}
\xymatrix@R=0.8em@C=0.1em{
&&(N.B'.A')(G_0'.D')(F.B')(A.G)\ar[drr]^>>>>>>>{-\nu^N-}&&\\
(N.B'.A')(G.A')(B.F)(D.F'_0) \ar[rru]^<<<<<<<<{-\delta-} &&&& (G_0.D')(F.B')(A.N.B')(A.G)\\
&&&& \\
(G.A')(B.N.A')(B.F)(D.F'_0)  \ar[uu];[]_{-\nu^S}  \ar[drr]_{-\nu^N-}&&&&(G_0.D')(F.B')(A.G)(A.B.N) \ar[uu];[]^{-\nu^S-}  \\
&& (G.A')(B.F)(D.F_0)(B.A.N)\ar[urr]_{-\delta'-}
}
\eeq
}

\noindent A 3-cell in $\mll\mll(\B)$ amounts to a 3-cell which is a modification between two pairs of l-transformations. 

In the context relative to a functor $\mathcal{I} : \K \rightarrow \C$, we have a tricategory $\mll\mll(K, C)$ with objects ll-distributive pairs of an l-monad in $\K$ and an l-monad in $\C$.

Besides $\mll\mll(\B)$, there are other tricategories of lax distributive pairs, lax morphisms of distributive pairs and their lax transformations. Here is a generic definition of such a tricategory:

\begin{itemize} 
\item For objects there are four possibilities corresponding to lax monads in tricategories $\mll(\C)$, $\mlcol(\C)$, $\mlcol(\C)$ and $\mcolcol(\C)$. 
\item Morphisms are triples $(\FF, \GG, \delta)$, where $\FF$ and $\GG$ are up- or downmaps of monads of arbitrary laxity, going between the two corresponding component l-monads of lax distributive pairs, and $\delta$ is a 3-cell of the form (\ref{delta}), taking any direction such that the commutativity of the diagrams (\ref{axnunu}) makes sense and holds. The 3-cell $\delta$ can be thought of as a distributivity law between the lax up/downmaps $\FF$ and $\GG$. 

\item 2-cells are pairs $(\NN, \SS)$, where $\NN$ and $\SS$ are transformations of arbitrary laxity going between the two corresponding component lax up/downmaps of morphisms previously defined, such that commutativity of (\ref{axnunu}) makes sense and holds. 

\item A 3-cell amounts to a 3-cell of $\C$ which becomes a modification for the two component transformations of a 2-cell previously defined.
\end{itemize}

Following this generic definition, a morphism $(\FF, \GG, \delta)$ of $\mll\mll(\B)$ can be thought of as a pair of l-upmaps $\FF$ and $\GG$ related by a distributivity law $\delta$. Looking at the diagrams (\ref{axdelta}) we observe that, if we replace $\delta$ by a morphism with an opposite direction $\bar\delta$, the commutativity will still make sense. Such a 3-cell $\bar\delta$ gives another type of a distributivity law between the l-upmaps $\FF$ and $\GG$. Now define a tricategory $\overline{\mll\mll}(\C)$. Its objects are ll-distributive pairs. A morphism $(\BB', \AA') \rightarrow (\BB, \AA)$ consists of a triple $(\FF, \GG, \bar\delta)$ where $\FF$ and $\GG$ are l-upmaps of l-monads as in a morphism of $\mll\mll(\C)$, but $\bar\delta$ takes the opposite direction to $\delta$, so it is a 3-cell
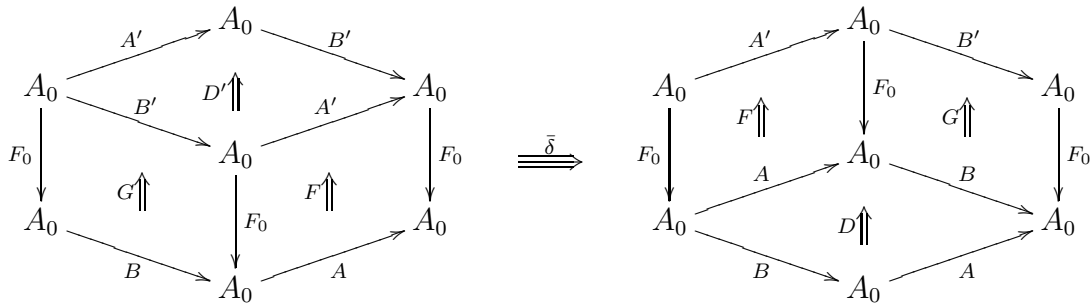
\[
\xymatrix@R=0.8em{
&& A_0  \ar[drr]^{B'}&&\\
A_0  \ar[rru]^{A'}  \ar[drr]^{B'}_\;="1"&&&& A_0  \\
&& A_0 \ar[urr]^{A'}_\;="2" \ar@{=>}  {}+(0,6); {}+(0,11) ^{D'} && \\
A_0  \ar[uu];[]_{F_0} \ar[drr]_{B}^\;="3"&&&&A_0 \ar[uu];[]^{F_0} \\
&&A_0 \ar[urr]_{A}^{\;}="4"   \ar[uu];[]^{F_0}&& \ar@{=>} "3"+(0,5);"3"+(0,10) ^{G} \ar@{=>} "4"+(0,5);"4"+(0,10)  ^F
}
\quad
\xymatrix@R=1em{\\\\\\ \ar@3{->}[r] ^{\bar\delta} & \\ \\}
\quad
\xymatrix@R=0.8em{
&& A_0 \ar[drr]^{B'}_\;="2"&&\\
A_0  \ar[rru]^{A'}_{\;}="1" &&&& A_0  \\
&& A_0  \ar[uu];[]^{F_0} \ar[drr]^{B}="4" && \\
A_0  \ar[urr]^{A}="3" \ar[uu];[]_{F_0}  \ar[drr]_{B}&&&&A_0 \ar[uu];[]^{F_0}  \\
&&A_0 \ar[urr]_{A} \ar@{=>} {}+(0,6); {}+(0,11) ^D && \ar@{=>}"3"+(0,5);"3"+(0,10) ^{F} \ar@{=>}"4"+(0,5);"4"+(0,10) ^{G}
}
\]

\noindent satisfying the equations obtained by modifying (\ref{axdelta}) in the obvious way. A 2-cell of $\overline{\mll\mll}(\C)$ consists of a pair of l-transformations $\NN$ and $\SS$, such that the obvious modification of the equation (\ref{axnunu}) holds. 3-cells are straightforward. 

In the relative context, the tricategory $\overline{\mll\mll}(\K, \C)$ can be defined much like its non-relative version. 

\begin{remark}
Note that, if $(\FF, \GG, \bar\delta)$ is a morphism of $\overline{\mll\mll}(\C)$, then, the pair $(G, \bar\delta)$ is a cl-transformation of l-upmaps $\FF\BB' \rightarrow \BB\FF : \AA' \rightarrow \AA$ (where $\BB$ is considered as an l-upmap $\AA \rightarrow \AA$). So $\bar\delta$ has a different laxity from $\kappa^A$, when the latter is considered as part of the l-transformations $\mathbb{P}_n^{\BB} : \BB^n \rightarrow \BB$. For this reason a morphism of $\overline{\mll\mll}(\C)$ is not an l-upmap in a tricategory of l-monads. Another example of a morphism of ll-distributive pairs which is not a lax map in any tricategory is obtained by allowing $\kappa^F$ to have the opposite laxity to $\kappa^B$, which will result in the upmaps $\FF$ and $\BB$ living in different tricategories. 
\end{remark}

\begin{example} Since $\mll(\cat^\op, \mod)$ is a 2-category, for $\mll\mll(\cat^\op, \mod)$ and $\overline{\mll\mll}(\cat, \mod)$ we should rather write $\fM^\bot\mll(\cat^\op, \mod)$ and $\overline{\fM^\bot\mll}(\cat^\op, \mod)$. In the view of Theorem~\ref{LK}, these are the opposites of the 2-categories of $T$-equipments:

\begin{theorem}
The 2-category $\fM^\bot\mll(\cat^\op, \mod)$ is isomorphic to $\fM^\top(\Eq)^\op$, and the 2-category $\overline{\fM^\bot\mll}(\cat^\op, \mod)$ is isomorphic to $\overline{\fM^\top(\Eq)}^\op$.
\end{theorem}

\noindent In more details, a morphism $(\FF_0, \FF, \nu^d) : (\SS_0, \BB) \rightarrow (\TT_0, \AA)$ in $\fM^\top(\Eq)$ in the monadic terms consists of a downmap $\FF_0 : \SS_0 \rightarrow \TT_0$ of monads in $\cat$, an l-upmap $\FF : \BB \rightarrow \AA$ of l-monads in $\mod$, and a distribution 3-cell of $\mod$

\[
\xymatrix@R=0.8em{
&& A_0 \ar[drr]^{T^\o_0}_\;="2"&&\\
A_0  \ar[rru]^{A}_{\;}="1" &&&& A_0  \\
&& B_0  \ar[uu];[]^{F^\o_0} \ar[drr]^{S^\o_0}="4" && \\
B_0  \ar[urr]^{B}="3" \ar[uu];[]_{F^\o_0}  \ar[drr]_{S^\o_0}&&&&B_0 \ar[uu];[]^{F^\o_0}  \\
&&B_0 \ar[urr]_{B} \ar@{=>} {}+(0,6); {}+(0,11) ^S && \ar@{=>}"3"+(0,5);"3"+(0,10) ^{F} \ar@{=>}"4"+(0,5);"4"+(0,10) ^{d}
}
\quad
\xymatrix@R=1em{\\\\\\ \ar@3{->}[r] ^{\nu^d} & \\ \\}
\quad
\xymatrix@R=0.8em{
&& A'_0  \ar[drr]^{T^\o_0}&&\\
A'_0  \ar[rru]^{A}  \ar[drr]^{T^\o_0}_\;="1"&&&& A'_0  \\
&& A'_0 \ar[urr]^{A}_\;="2" \ar@{=>}  {}+(0,6); {}+(0,11) ^{T} && \\
A_0  \ar[uu];[]_{F^\o_0} \ar[drr]_{S^\o_0}^\;="3"&&&&A_0 \ar[uu];[]^{F^\o_0} \\
&&A_0 \ar[urr]_{B}^{\;}="4"   \ar[uu];[]^{F^\o_0}&& \ar@{=>} "3"+(0,5);"3"+(0,10) ^{d} \ar@{=>} "4"+(0,5);"4"+(0,10)  ^F
}
\]

\noindent A morphism of $\overline{\fM^\top(\Eq)}$ has a similar characterization with the distribution 3-cell taking the opposite direction. Besides of these, there are few different 2-categories for each notion of $T$-equipment corresponding to lax monads in $\Eq^{\circ\bullet}$, $\Eq^{\bullet\circ}$ and $\Eq^{\circ\bullet}$.
\end{example}

\smallskip Define a (left) \textbf{(cl/)l-comodule} of an l-monad $\AA$ in $\C$ to consists of an object $Z$ and a (cl/)l-upmap $\AA \rightarrow \id(Z)$, where $\id(Z)$ denotes the trivial l-monad on $Z$. Further define a (left) \textbf{(cl/)l-module} as a (cl/)l-downmap $\AA \rightarrow \id(Z)$.

\begin{example} A monoid $(x, a, \mu_a, \eta_a)$ in an equipment $\AA$, as per Definition~\ref{monoid}, is an l-comodule of the l-monad $\AA$ in $\mod$, consisting of the morphism $x^\o : A_0 \rightarrow I_0$, the 2-cell $a : x^\o \Rightarrow x^\o\circ A$, and the 3-cells
\[
\xymatrix@R=0.8em{
&&I_0 \ar[lldd];[]^{x^\o}="1" \ar[dd];[]_{x^\o} \ar[rrdd];[]_{x^\o}&&\\\\
A_0 \ar@/_2.5em/[rrrr]_A \ar[rr]_A&& A_0 \ar[rr]_A \ar@{=>} {}+(0,-3); {}+(0,-8) ^{P_2^A}&& A_0 \ar@{=>} "1"+(4,-5); "1"+(9,-5) ^a \ar@{=>} "1"+(21,-5); "1"+(26,-5) ^{a}
}
\quad
\xymatrix@R=0.8em{
\\
\ar@3{->}[r]^{\mu_a}&
}
\quad
\xymatrix@R=0.8em{
&&I_0 \ar[lldd];[]^{x^\o}="1" \ar[rrdd];[]_{x^\o}&&\\\\
A_0 \ar[rrrr]_A&& \ar@{=>} "1"+(12,-4); "1"+(17,-4) ^a&& A_0 
}
\]
\[
\xymatrix@R=0.8em{
&&I_0 \ar[lldd];[]^{x^\o}="1" \ar[rrdd];[]_{x^\o}&&\\\\
A_0 \ar@/_2.5em/[rrrr]_A \ar[rrrr]^{1_{A_0}}&& \ar@{=>} {}+(0,-3); {}+(0,-8) ^{P_0^A} \ar@{=>} "1"+(12,-4); "1"+(17,-4) ^a&& A_0 
}
\quad
\xymatrix@R=0.8em{
\\
\ar@3{->}[r]^{\eta_a}&
}
\quad
\xymatrix@R=0.8em{
&&I_0 \ar[lldd];[]^{x^\o}="1" \ar[rrdd];[]_{x^\o}&&\\\\
A_0 \ar[rrrr]_A&& \ar@{=>} "1"+(12,-4); "1"+(17,-4) ^a&& A_0 
}
\]

\noindent A cl-comodule with the underlying morphism $x^\o : A_0 \rightarrow I_0$ defines a \textbf{comonoid in an equipment}. A $T$-monoid $(x, a, \mu_a, \eta_a)$ in a $T$-equipment $(\TT_0, \AA)$, as per Definition~\ref{tmonoid}, is an l-comodule of the l-monad $\comp(\TT_0, \AA)$, consisting of the morphism $x : A_0 \rightarrow I_0$, the 2-cell $a : x^\o \Rightarrow x^\o\circ A\circ T_0^\o$, and 3-cells
\[
\xymatrix@R=0.8em{
&&I_0 \ar[lldd];[]^{x^\o}="1" \ar[dd];[]_{x^\o} \ar[rrdd];[]_{x^\o}&&\\\\
A_0 \ar@/_1.3em/[rrdd]_A \ar[r]_A& A_0 \ar@{=>} {}+(0,-4); {}+(0,-9) ^{P_2^A} \ar[rdd]_A \ar[r]_{T_0^\o}& A_0 \ar[r]_A \ar@{=>} {}+(0,-4); {}+(0,-9) ^{T}& A_0 \ar[r]_{T_0^\o} \ar@{=>} {}+(0,-4); {}+(0,-9) ^{m} & A_0 \ar@{=>} "1"+(5,-5); "1"+(10,-5) ^a \ar@{=>} "1"+(23,-5); "1"+(28,-5) ^{a}\\\\
&&A_0 \ar[ruu]_{T^\o_0} \ar@/_1.3em/[rruu]_{T^\o_0}&&  
}
\quad
\xymatrix@R=0.8em{
\\
\ar@3{->}[r]^{\mu_a}&
}
\quad
\xymatrix@R=0.8em{
&&I_0 \ar[lldd];[]^{x^\o}="1" \ar[rrdd];[]_{x^\o}&&\\\\
A_0 \ar[rr]_A&& A_0 \ar[rr]_{T^\o_0} \ar@{=>} "1"+(12,-4); "1"+(17,-4) ^a&& A_0 
}
\]
\[
\xymatrix@R=0.8em{
&&I_0 \ar@{}[dd]|-{\hequ}\ar[lldd];[]^{x^\o} \ar[rrdd];[]_{x^\o}&&\\\\
A_0 \ar[rr]^{1_{A_0}}="1" \ar@/_2.3em/[rr]_A&& A_0 \ar[rr]^{1_{A_0}}="2" \ar@/_2.3em/[rr]_{T^\o_0} && A_0 \ar@{=>} "1"+(0,-4); "1"+(0,-9) ^{P_0^A} \ar@{=>} "2"+(0,-4); "2"+(0,-9) ^e
}
\quad
\xymatrix@R=0.8em{
\\
\ar@3{->}[r]^{\eta_a}&
}
\quad
\xymatrix@R=0.8em{
&&I_0 \ar[lldd];[]^{x^\o}="1" \ar[rrdd];[]_{x^\o}&&\\\\
A_0 \ar[rr]_A&& A_0 \ar[rr]_{T^\o_0} \ar@{=>} "1"+(12,-4); "1"+(17,-4) ^a&& A_0. 
}
\]
\end{example}

An \textbf{ll-distributive pair of l-comodules} over an ll-distributive pair of l-monads $(\BB, \AA)$ is a morphism $(\BB, \AA) \rightarrow \id^2(Z)$ in $\mll\mll(\C)$, where $\id^2(Z)$ is a trivial ll-distributive pair of l-modules on the object $Z$. This amounts to a pair of l-comodules respectively in $\AA$ and $\BB$ related by an additional distributivity data. Furthermore, define generically a lax distributive pair of cl/l-(co)modules as a morphism $(\BB, \AA) \rightarrow \id^2(Z)$ in any of the tricategories of ll-distributive pairs of l-monads.

\begin{example} For any equipment $\AA$, consider the ll-distributive pair of l-monads $(\AA, \AA)$  with the trivial distribution structure. An ll-distributive pair of l-comodules over $(\AA, \AA)$ with an underlying morphism $x^\o$, should be regarded as a distributive pair of monoids in $\AA$, similar to a distributive pair of monads in a bicategory. A lax distributive pair of an l-comodule and a cl-comodule with an underlying morphism $x^\o$ is a mixed distributive pair of a monoid and a comonoid in the equipment. 
\end{example}

\begin{example} A $T$-algebra $(x, a, \mu_a, \eta_a, h, \sigma_h)$ in a $T$-equipment $(\TT_0, \AA)$, as in Definition~\ref{talgebra}, is a morphism  $(\TT_0, \AA) \rightarrow \id^2(I_0)$ in $\fM^\bot\mll(\cat^\op, \mod)$. So, it is a lax distributive pair consisting of the l-comodule $(x^\o, a, \mu_a, \eta_a)$ of the l-monad $\AA$ in $\mod$ and a module $(x, h)$ of the monad $\TT_0$ in $\cat$ related by the distributivity data
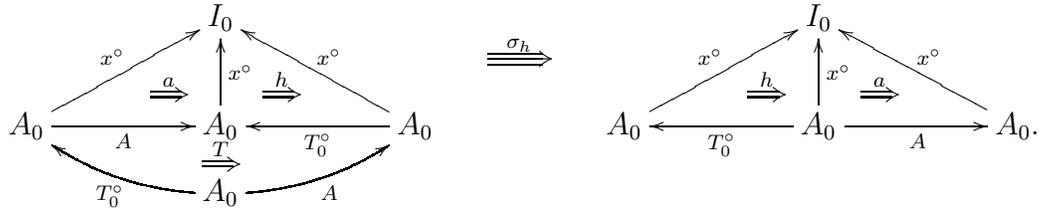
\[
\xymatrix@R=0.8em{
&&I_0 \ar[lldd];[]^{x^\o}="1" \ar[dd];[]_{x^\o} \ar[rrdd];[]_{x^\o}&&\\\\
A_0 \ar[rr]_{A}&& A_0 \ar@{=>} {}+(-2.5,-5); {}+(2.5,-5) ^T \ar[rr];[]^{T^\o_0} && A_0 \ar@{=>} "1"+(5,-5); "1"+(10,-5) ^a \ar@{=>} "1"+(20,-5); "1"+(25,-5) ^{h}\\
&&A_0 \ar@/^.6em/[llu]^{T^\o_0} \ar@/_.6em/[rru]_{A}&&
}
\quad
\xymatrix@R=0.8em{
\\
\ar@3{->}[r]^{\sigma_h}&
}
\quad
\xymatrix@R=0.8em{
&&I_0 \ar[lldd];[]^{x^\o}="1" \ar[dd];[]_{x^\o} \ar[rrdd];[]_{x^\o}&&\\\\
A_0 \ar[rr];[]^{T^\o_0}&& A_0 \ar[rr]_A && A_0. \ar@{=>} "1"+(5,-5); "1"+(10,-5) ^h \ar@{=>} "1"+(20,-5); "1"+(25,-5) ^{a}
}
\]
\end{example}

\smallskip Now we introduce the lax counterpart of the composite of a distributive pair of monads: 

\begin{definition} A \textbf{composite of an ll-distributive pair of l-monads} $(\BB, \AA)$, denoted $\comp(\BB, \AA)$, is an l-monad $(A_0, B.A, P^{BA}, \xi^{BA})$ with the following components:

\begin{itemize} 
\item An object $A_0$.

\item A morphism $B.A : A_0 \rightarrow A_0$.

\item The multiplication 2-cells $P^{BA}$ define by
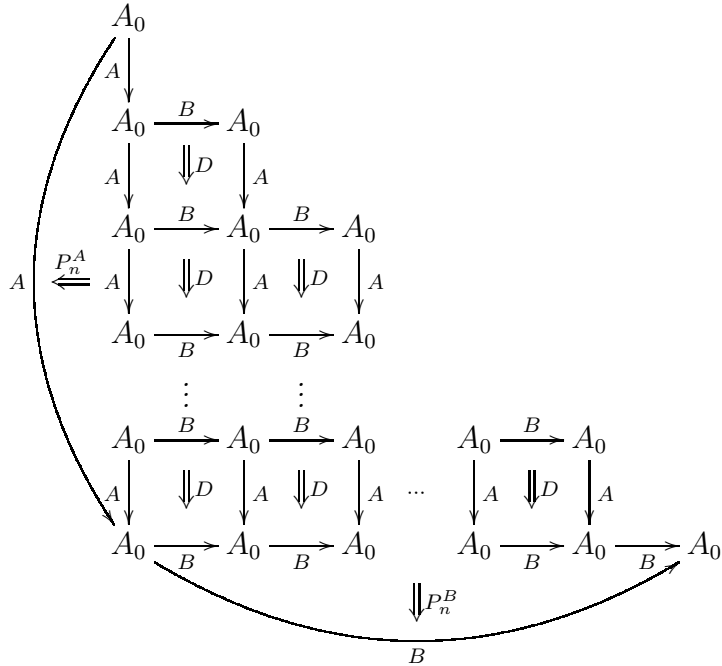
\[
\xymatrix{
A_0 \ar[d]_A \ar@/_3pc/[ddddd]_A="0"&&&&\\
A_0\ar[d]_A \ar[r]^B="1"&A_0\ar[d]^A&&&&\\
A_0\ar[d]_A="12" \ar[r]^B="2"&A_0 \ar[d]^A \ar[r]^B="3"&A_0 \ar[d]^A&&&\\
A_0\ar[r]_B="4" \ar@{}[dr]|-{\vdots}&A_0\ar@{}[dr]|-{\vdots}\ar[r]_B="5"& A_0&&\\
A_0\ar[d]_A \ar[r]^B="6"&A_0 \ar[d]^A \ar[r]^B="7"&A_0\ar@{}[dr]|-{\cdots}\ar[d]^A&A_0 \ar[r]^B="8" \ar[d]^A&A_0\ar[d]^A&\\
A_0\ar[r]_B="9"  \ar@/_3pc/[rrrrr]_B="14"&A_0\ar[r]_B="10"&A_0 \ar@{}[r]_{}="13"&A_0\ar[r]_B="11"&A_0\ar[r]_B& A_0
\ar@{=>} "1"+(0,-5); "1"+(0,-10) ^D \ar@{=>} "2"+(0,-6); "2"+(0,-11) ^D \ar@{=>} "3"+(0,-6); "3"+(0,-11) ^D \ar@{=>}"6"+(0,-6); "6"+(0,-11)^D \ar@{=>} "8"+(0,-6); "8"+(0,-11) ^D \ar@{=>} "7"+(0,-6); "7"+(0,-11) ^D \ar@{=>} "13"+(0,-5); "13"+(0,-10) ^{P^B_n} \ar@{=>} "12"+(-3,0); "12"+(-8,0) _{P^A_n}
}
\]

\noindent Note that, at $n = 0$ this becomes 
\[
\xymatrix{A_0 \ar@/^1.2pc/[r]^{1_{A_0}}_{\quad}="1" \ar@/_1.2pc/[r]_{A}^{\quad}="2"  \ar@{=>}"1" ; "2" ^{P^A_0}& A_0 \ar@/^1.2pc/[r]^{1_{A_0}}_{\quad}="3" \ar@/_1.2pc/[r]_{B}^{\quad} ="4" \ar@{=>}"3" ; "4" ^{P^B_0}&A_0.}
\]

\item The associator 3-cells $\xi^{BA}$ built by repetitive applications of the 3-cells $\kappa^A$ and $\kappa^B$ and the associator 3-cells $\xi^A$ and $\xi^B$ of the l-monads $\AA$ and $\BB$.

\end{itemize}

Similarly, relative to a functor $\K \rightarrow \C$, defined is the composite of an ll-distributive pair of an l-monad $\BB$ in $\K$ and an l-monad $\AA$ in $\C$
\end{definition}

The construction of composite of ll-distributive pairs of l-monads extends to a functor. Given a morphism $(\FF, \GG, \delta) : (\BB', \AA') \rightarrow (\BB, \AA)$ of $\overline{\mll\mll}(\C)$, the l-upmap of l-monads $\comp(\FF, \GG, \delta) : \comp(\BB', \AA') \rightarrow \comp(\BB, \AA)$ is defined by the data $(F_0, GF, \kappa^{GF})$ consisting of:

\begin{itemize}
\item The morphism $F_0 : A_0 \rightarrow A_0$.
\item The square $GF$: 
\[
\xymatrix{A'_0 \ar[rr]^{A'}_{\;}="1" \ar[d]_{F_0} && A'_{0}  \ar[d]^{F_0} \ar[rr]^{B'}_{\;}="3" && A'_{0} \ar[d]^{F_0} \\
Y_0 \ar[rr]_{A}="2" && Y_{1} \ar[rr]_B & & Y_2 \ar@{=>} "1"+(0, -8.5); "1"+(0, -3.5) _{F} \ar@{=>} "3"+(0, -8.5); "3"+(0, -3.5) _{G},
}
\]
\item The comparison 3-cells $\kappa^{GF}$ defined by:
\[
\xymatrix{((GF)^n)(``D^{(n-1)n/2}")(P_n^B.P_n^A) \ar@3{->}[d]^{``\delta^{(n-1)n/2}"}\\(``D'^{(n-1)n/2}")((GF)^n)(P_n^B.P_n^A) \ar@3{->}[d]^{-\kappa^G_n.\kappa^F_n-}\\ (``D'^{(n-1)n/2}")(P_n^{B'}.P_n^{A'})(GF).
} 
\]
\end{itemize}

\noindent We leave it to the reader to define $\comp$ on 2- and 3-cells, and conclude that it is a functor

\[\xymatrix{\comp_\C^\bot: \overline{\mll\mll}(\C) \ar[r]& \mll(C).}\]

\noindent In the relative context, the composite of ll-distributive pairs becomes a functor
\[\xymatrix{\comp_{\K, \C}^\bot : \overline{\mll\mll}(\K, \C) \ar[r]& \mll(\K, \C).}\]

\begin{example} In the view of the fact that a $T$-equipment is an ll-distributive pair, and observations made in Section~\ref{biprofmod} we can easily deduce that:

\begin{theorem}
For a $T$-equipment $(\TT_0, \AA)$, we have $\kl(\TT_0, \AA) = \comp(\TT_0, \AA)$. Moreover, the Kleisli functor
\[\xymatrix{\kl : \overline{\fM^\top(\Eq)} \ar[r]& \Eq}\]

\noindent is the same as the distributive composition functor
\[\xymatrix{\comp_{\cat^\op, \mod}^\bot : \overline{\fM^\bot\mll}(\cat^\op, \mod)^\op \ar[r] & \mll(\cat^\op, \mod)^\op.}\]

\end{theorem}
\end{example}

\smallskip In our current setting we do not have a non ad hoc description of the 2-category of $\ast$-equipments $\sEq$, or either of the 2-categories of $T$-$\ast$-equipments. However the construction that we describe next is closely related to the $\ast$-Kleisli 2-functor $\kl^\ast: \fM_{\mathrm{C}}^\bot(\sEq) \rightarrow \Eq.$  

Consider the tricategory $\mocolcol\mll(\C)$. Its objects are ll-distributive pairs of l-monads. A morphism $(\BB', \AA') \rightarrow (\BB, \AA)$ in it is a triple $(\GG, \FF, \delta)$, where $\FF$ is an l-upmap $\AA' \rightarrow \AA$, $\GG$ is a cl-downmap $\BB' \rightarrow \BB$, and $\delta$ is distribution 3-cell
\[
\xymatrix@R=0.8em{
&& A'_0  \ar[drr]^{A'}&&\\
A'_0  \ar[rru];[]_{B'}  \ar[drr]^{A'}_\;="1" \ar@{=>}{}+(23,0);{}+(28,0)^{D'}&&&& A'_0  \\
&& A'_0 \ar[urr];[]_{B'}^\;="2" && \\
A_0  \ar[uu];[]_{F_0} \ar[drr]_{A}^\;="3"&&&&A_0 \ar[uu];[]^{F_0} \\
&&A_0 \ar[urr];[]^{B}_{\;}="4"   \ar[uu];[]^{F_0}&& \ar@{=>}"3"+(0,5); "3"+(0,10) ^{F} \ar@{=>}"2"+(0,-5); "2"+(0,-10) ^G
}
\quad
\xymatrix@R=1em{\\\\\\ \ar@3{->}[r]^\delta & \\ \\}
\quad
\xymatrix@R=0.8em{
&& A'_0 \ar[drr]^{A'}_\;="2"&&\\
A'_0  \ar[rru];[]_{B'}_{\;}="1" &&&& A'_0  \\
&& A_0  \ar[uu];[]^{F_0} \ar[drr]_{A}="4" && \\
A_0  \ar[urr];[]^{B}="3" \ar[uu];[]_{F_0}  \ar[drr]_{A} \ar@{=>}{}+(23,0);{}+(28,0)^{D} &&&&A_0 \ar[uu];[]^{F_0}  \\
&&A_0 \ar[urr];[]^{B} && \ar@{=>} "1"+(0,-7); "1"+(0,-12) ^G \ar@{=>}"4"+(0,8);"4"+(0,13) ^{F}
}
\]

\noindent satisfying some equations. Suppose that all 2-cells in $\C$ have right adjoints. Under this condition, we will build a functor

\begin{equation}\label{cltol}
\xymatrix{\mocolcol\mll(\C) \ar[r]& \overline{\mll\mll}(\C).}
\end{equation}

\noindent First, construct a functor 
\[\xymatrix{\mocolcol(\C) \ar[r]& \mll(\C).}\]

\noindent On objects set it to be identical. On 1-, 2-, 3- cells define it by the following correspondences. Suppose that $\GG = (G_0, G, \kappa^G)$ is a cl-downmap $\BB \rightarrow \AA$. Let $G^\ast : G_0B \Rightarrow AG_0$ be the right adjoint 2-cell to $G : AG_0 \Rightarrow  G_0B$. Let $\kappa^{G^\ast}_n : (F_0.P_n^{B'})G^{\ast n} \Rightarrow G^\ast (P^B_n.F_0)$ be the mate of $\kappa^G_n : G(F_0.P_n^{B'}) \Rightarrow  (P^B_n.F_0)G^{n}$ under the adjunctions $G \dashv G^\ast$ and $G^n \dashv  G^{\ast n}$. Then, $\GG^\ast = (G_0, F^\ast, \kappa^{G^\ast})$ is an l-upmap $\AA \rightarrow \BB$. Now suppose that $(N, \nu)$ is a l-transformation of cl-upmaps $\GG' \rightarrow \GG$. Let $\tilde\nu^N : (N.B')G^\ast \Rightarrow G^{\prime \ast}(B.N)$ be the mate of $\nu : G'(N.B') \Rightarrow (B.N)G$ under the adjunctions $G \dashv G^\ast$ and $G' \dashv G^{\prime\ast}$. Then $N^\ast = (N, \tilde\nu^N)$ is an l-transformation of l-upmaps $\GG^{\prime\ast} \rightarrow \GG^\ast$. Finally if a 3-cell is a modification $\NN \rightarrow \SS$ of cl-transformation of cl-downmaps, then the same 3-cell becomes a modification $\NN^\ast \rightarrow \SS^\ast$ of l-transformations of l-downmaps.  Now we can describe the functor (\ref{cltol}). On objects it is identical. To a morphism $(\GG, \FF, \delta)$ of $\mocolcol(\C)$ it assigns the morphism of  $\overline{\mll\mll}(\C)$ consisting of the triple $(\GG^\ast, \FF, \tilde\delta)$ where 
\[\xymatrix{(F_0.D')(F.B')(A.G^\ast) \ar@{=>}[r]^{\tilde\delta} & (G^\ast.A')(B.F)(D.F_0)}\] 

\noindent is the mate of 
\[\xymatrix{(G.A')(F_0.D')(F.B')  \ar@{=>}[r]^\delta & (B.F)(D.F_0)(A.G)}\]

\noindent under the adjunctions $A.G \dashv A.G^\ast$ and $G.A'\dashv G^\ast.A'$. To a 2-cell $(\NN, \SS)$ it assigns the 2-cell $(\NN^\ast, \SS)$. On 3-cells it is trivial.

Precomposing (\ref{cltol}) with the functor $\comp_\C^\bot$, we get another functorial extension of the ll-distributive composition 
\[\xymatrix{\comp_\C^{\bot\ast} : \mocolcol\mll(\C) \ar[r]& \mll(\C).}\]

\noindent On objects $\comp^{\bot\ast}$ coincides with $\comp^\bot$. On morphisms, $\comp^{\bot\ast}(\GG, \FF, \delta)$ is an l-map of l-monads $(F_0, G^\ast.F, \kappa^{G^\ast F})$ consisting of:

\begin{itemize}
\item The morphism $F_0$.
\item The square $G^\ast.F$:
\[
\xymatrix{A'_0 \ar[rr]^{A'}_{\;}="1" \ar[d]_{F_0} && A'_{0}  \ar[d]^{F_0} \ar[rr]^{B'}_{\;}="3" && A'_{0} \ar[d]^{F_0} \\
Y_0 \ar[rr]_{A}="2" && Y_{1} \ar[rr]_B & & Y_2 \ar@{=>} "1"+(0, -8.5); "1"+(0, -3.5) _{F} \ar@{=>} "3"+(0, -8.5); "3"+(0, -3.5) _{G^\ast},
}
\]

\item A lax comparison $\kappa^{G^\ast F}$ defined from $\kappa^F$, $\kappa^{G^\ast}$ and $\tilde{\delta}$ by
\beq{kgastf}
\xymatrix{((GF)^n)(``D^{(n-1)n/2}")(P_n^B.P_n^A) \ar@3{->}[d]^{``\tilde\delta^{(n-1)n/2}"}\\(``D'^{(n-1)n/2}")((GF)^n)(P_n^B.P_n^A) \ar@3{->}[d]^{-\kappa^{G^\ast}.\kappa^F_n-}\\ (``D'^{(n-1)n/2}")(P_n^{B'}.P_n^{A'})(GF).
} 
\eeq
\end{itemize}

\noindent The relative counterpart is the functor
\beq{compast}
\xymatrix{\comp_{\K, \C}^{\bot\ast} : \mocolcol\mll(\K, \C) \ar[r] & \mll(\K, \C)}.\eeq

\begin{example} Obviously, $\mocolcol\mll(\cat^\op, \mod) =\fM^\top\mll(\cat^\op, \mod) = \fM^\bot(\Eq)^\op$. However, we can not use (\ref{compast}) because 2-cells in $\cat$ do not have right adjoints, neither do their images in $\mod$. The following, which we do not try to make completely precise, rectifies this. 

Suppose that $A$ is a $\ast$-module from $A_0$ to itself. Suppose that $t : F \rightarrow G$ is a natural transformation whose opaction on $A$ is Cartesian in the sense of Section~\ref{starequipments}. Then, $t\circ A : F\circ A, \rightarrow G\circ A$, which is a pseudonatural transformation with the components $A(x, t_x) : A(y, Fx) \rightarrow A(y, Gy)$, has a right adjoint  in $\mod$ given by the family of functors $A(x, t^\ast_x) : A(y, Gx) \rightarrow A(y, Fy)$, which is pseudonatural by the virtue of the Cartesian property. Denote this right adjoint by $t^\ast\circ A$. This structure on $\cat^\op \rightarrow \mod$ allows a partial definition of $\comp^\bot$. Given a morphism $(\FF_0, \FF, \nu^d) : (\SS_0, \BB) \rightarrow (\TT_0, \AA)$ of $\fM^\top\mll(\cat^\op, \mod)$, with the module $A$ given a $\ast$-structure, and the opaction of $d : F_0S_0 \rightarrow T_0F_0$ on $A$ being Cartesian, define  $\comp^\bot(\FF_0, \FF, \nu^d)$ by the data $(F_0, d^\ast\circ F, \kappa^{d^\ast\circ F})$ consisting of 

\begin{itemize}
\item The morphism $F_0$ of $\cat$.
\item A 2-cell of $\mod$ $d^\ast\circ F$ defined as $(d^\ast\circ A)(S_0\circ F)$.
\item The lax associator $\kappa^{d^\ast\circ F}$, defined as (\ref{kgastf}), which makes use of $\kappa^F$, a transform of the equality satisfied by $d$, and a transform of $\nu^d$. 
\end{itemize}

\noindent $\comp^\bot$ defined in this way coincides with the $\ast$-Kleisli 2-functor $\kl^\ast$.

\end{example}

\begin{remark} Recall from Section~\ref{starequipments}, that a $\ast$-module $A$ from $X$ to $Y$ extends to a functor $\Pi(X)\times\Pi(Y)\rightarrow \cat$. This is a morphism $\Pi(X) \rightarrow \Pi(Y)$ in $\twoprof$. We can deduce that a $\ast$-equipment $\AA$ is a lax monad in $\twoprof$ on the object $\Pi(A_0)$. Furthermore, a functor $F_0 : X \rightarrow Y$ gives rise to a functor $\Pi(F_0) : \Pi(X) \rightarrow \Pi(Y)$, which itself determines a morphism $\Pi(F_0)^\o : \Pi(Y) \rightarrow \Pi(X)$ of $\twoprof$. Functors of $\ast$-equipments can be described as l-upmaps of l-monads of the form $(\Pi(F_0), F^{\Pi(F_0)}, \kappa^{\Pi(F_0)})$. However, a natural transformation $t : F_0 \rightarrow G_0$ does not extend to a transformation $\Pi(F_0) \rightarrow \Pi(G_0)$ or to a 2-cell $\Pi(F_0)^\o \rightarrow \Pi(G_0)^\o$ of $\twoprof$.
\end{remark}

\section{Further comparisons with virtual double categories}\label{further}
As stated before, our equipments are close to the virtual double categories of \cite{CrSh10}, and can replace them in practice. In particular, one can reproduce for equipments different constructions on virtual double categories. An example is the construction $\MMod$ of \cite{CrSh10}.

\begin{definition}
Suppose that $\AA$ is an equipment whose categories of vectors have coequilizers which are preserved by compositions. Then, the equipment $\MMod(\AA)$ is defined to consist of:
\begin{itemize}
\item The category of monoids $\mon(\AA)$ as the category of scalars.
\item For each pair of monoids the category of vectors $\MMod(\AA)((x, a), (y, b))$ defined to have objects (\textbf{modules}) vectors $m : x \vectarr y$ together with 2-cells $ma \Rightarrow a$ and $bm \Rightarrow b$ compatible with the multiplication and units of $(x, a)$ and $(y, b)$, and morphisms (\textbf{module maps}) 2-cells $m \Rightarrow m'$ compatible with the left and the right actions.
\item The composition structure defined by the usual module composition using coequalizers in the vector categories of $\AA$.   
\end{itemize}
\end{definition}

Using $\MMod$ one obtains important new classes of equipments from the existing one. These classes of equipments themselves are used to capture important examples of operads as $T$-monoids. In particular, this is how one arrives at symmetric multicategories. Following \cite{CrSh10}, define $\MMod(\MMatV) = \PProf(V)$. The objects of $\PProf(V)$ are $V$-categories, its scalar arrows are $V$-functors, and its vectors are $V$-profunctors. $T$-monoids in $(\TT_0, \PProf(\set))$, where $\TT_0$ is the free symmetric strict monoidal category monad are ``enhanced'' symmetric multicategories. ``Normalizing'' the latter gives symmetric multicategories. For the details of this example, as well as other examples, see \cite{CrSh10}.

Interesting from the point of view of this paper is how $\MMod$ can be expressed in monad theoretic terms of Section~\ref{laxmonads}. Given an equipment $\AA$, let for the time being $M_0$ denote $\mon(\AA)$, and let $M$ denote the module of modules, so that $\MMod(\AA) = (M_0, M)$. Then, what we have is a morphism $M_0 \rightarrow A_0$ in $\cat$ , taking a monoid in $\AA$ to its underlying object, together with an obvious 2-cell in $\mod$
 
\[
\xymatrix{
M_0 \ar[rr]^{M}_{\;}="1"  \ar[d]&& M_0  \ar[d]\\
A_0 \ar[rr]_{A}&& A_0. \ar@{=>} "1"+(0, -3.5); "1"+(0, -8.5)
}
\]

\noindent Whether it is possible to describe $(M_0, M)$ internally to a tricategory (perhaps through some universal properties) remains to be seen. 

\smallskip Whether virtual double categories are lax monads in any tricategory is still to be examined. Possibly, this may be achieved by replacing morphisms and 2-cells of $\mod$ from pseudofunctors $X^\op\times Y \rightarrow \cat$ and pseudonatural transformations between them to some kind of morphisms $X^\op\times Y \rightarrow \prof$ to the bicategory of profunctors $\prof$ and transformations between them.

\bibliographystyle{alpha}

\end{document}